\documentclass[a4paper, 10pt]{article}
\usepackage{amsmath} \usepackage{amsfonts} \usepackage{amssymb}
\usepackage{latexsym}
\usepackage[english]{babel}
\usepackage{amscd}
\usepackage[dvips,final]{graphics}
\usepackage{fancyhdr}
\usepackage{locengli,math}
\usepackage{graphicx,pst-plot,psfig,pstricks,pst-node,pst-text,pst-tree,pst-char,pst-coil}

\newcommand{\LastUpdate}{10 04 2002}
\newcommand{\Ito}{\textrm{Ito}}

\begin{document}

\begin{center}
\textbf{\LARGE{\textsf{Coassociativity breaking and oriented graphs}}}

\vskip1cm
\parbox[t]{14cm}{\large{
Philippe {\sc Leroux}}\\
\vskip4mm
{\footnotesize
\baselineskip=5mm
Institut de Recherche
Math\'ematique, Universit\'e de Rennes I and CNRS UMR 6625\\
Campus de Beaulieu, 35042 Rennes Cedex, France, leroux@maths.univ-rennes1.fr}}
\end{center}

\vskip1cm
{\small
\centerline{\LastUpdate}
\vskip1cm
\baselineskip=5mm
\noindent
{\bf Abstract:}
With each coassociative coalgebra, we associate an oriented graph.
The coproduct $\Delta$, obeying the coassociativity equation $(\Delta \otimes id)\Delta = (id \otimes \Delta)\Delta$
is then viewed as a physical propagator which can convey information. We notice that
such a coproduct is non local. To recover locality we have
to break the coassociativity of
the coproduct and restore it, in some sense by introducing
another coproduct $\tilde{\Delta}$. The coassociativity equation becomes
the coassociativity breaking equation
$(\tilde{\Delta} \otimes id)\Delta = (id \otimes \Delta)\tilde{\Delta}$.
A coalgebra equipped with such coproducts, will be called a $L$-coalgebra. We prove then
that any oriented graph whose paths may be equipped with a probability measure,
can be seen as derived from such an algebraic formalism. The aim of this
article is to study the consequences of such a physical viewpoint into the algebraic formalism, especially for a unital algebra
equipped with its flower graph. We will show also the consequences of the notion of curvature of a 1-cochain
introduced by Quillen in this framework and a common point with the associative product of
this algebra,
homomorphisms, Leibnitz derivatives and Ito derivatives. As the concept of Ito derivative
appears naturally in this framework it will be studied throughout this article. In an attempt to adapt
what was done in the case of cyclic cocycles to the Ito case, we construct a di-superalgebra
(notions due to Quillen and Loday) from
the curvature of an Ito derivative. We show also that (pre)-dialgebras are in one to one
correspondence
with dendriform algebras equipped with one associative product. We generalise our concept
of $L$-coalgebras to $L$-coalgebras of degree $n$ and show that $M_2(k)$ and the quaternions
algebra $\mathbb{H}$ are (Markov) $L$-Hopf algebras of degree 2.
We go further by introducing
the notion of probabilistic algebraic product, notion which comes naturally from graph
theory and coalgebra (matrix product in relation with the graph of $Sl(2)_q$, wedge product
in relation with the oriented triangle graph).

\newpage
\tableofcontents
\newpage

\section{Introduction}
After introducing notation and definitions in part 2, we show how an oriented graph can emerge
from a coassociative coalgebra. To interpret physically this association,
we see the graph equipped with
its coproduct $\Delta$ as a physical model of space-time with a propagator $\Delta$. We notice
that such a coproduct, or propagator is non local and that the
mathematical concept of associativity
entails the lost of the physical concept of locality. However we show that this non locality
can be controlled by a (semi)-groupoid argument, (which could explain a phase transition between
locality and non locality). Indeed, groupoids need the notion of associativity expressed in
the definition of the concatenation (notion of locality) product.
To restore locality in such (coassociative) graphs, we need to break the coassociativity by
introducing  two coproducts (propagators) $\Delta$ and $\tilde{\Delta}$. They modelize respectively
our conceptions of future and past, not present
in the non local setting since $\Delta = \tilde{\Delta}$.
The r\^ole of part 3 is to study the consequences of such coalgebras, so-called $L$-coalgebras, equipped
with two coproducts $\Delta$ and $\tilde{\Delta}$. We show that any graph, whose paths may be equipped with a probability measure,
can be naturally embedded into a (Markov)-$L$-coalgebra. Several examples are given,
particularly the degenerate case, i.e. the coassociative coalgebra, which can be seen as
an (Ito)- $L$-coalgebra i.e. the two new coproducts constructed from the coproduct $\Delta$ become Ito derivatives
when $\Delta$ is a unital homomorphism. The important concept of Ito derivative will be studied
throughout this article. Part 5 is devoted to the study of $L$-coalgebra over
an associative algebra. Every unital associative algebra $A$, presents naturally two coproducts
$\delta$ and $\tilde{\delta}$ which embed $A$ into a $L$-bialgebra, so-called the flower graph.
These coproducts allow us to recover several concepts of algebra theory including the Hochschild
complex by reading periodic orbits on the flower graph. We go further by using Quillen's ideas on the curvature
of a 1-cochain to prove that, the product $m$ of the algebra $A$, homomorphisms, Leibnitz derivatives
and Ito derivatives, from $A$ to $A$, all verify a same equation. Motivated by this remark we establish a one-to-one
mapping between the set of Ito derivatives and the set of homomorphisms. This relation is explained
at the end of part 5 by interpreting Ito derivatives into an algebraic setting, showing
the lack of distributivity of the new product defined from Ito map and the associative product $m$.
We point then a relation between the third Reidemeister movement in knot theory and this setting. In an attempt to adapt
what was done in the case of cyclic cocycles to the Ito case, we construct a di-superalgebra
(notions due to Quillen and Loday) from
the curvature of an Ito derivative. The integral calculus so obtained yieds cyclic cocycles and vanishes on Leibnitz bracket.
We show also that (pre)-dialgebras are in one to one
correspondence
with dendriform algebras equipped with one associative product. We construct also non local
forms and a non local bundle from a coassociative coproduct $\Delta$ of a coassociative coalgebra $C$ and its embedding into
its natural Ito $L$-coalgebra, in studying the
deformation of an element $a \in A$, living on a petal of the flower graph of $A$ (the basis) from
its lift on its associated non local fiber via $\Delta$ i.e.
we study $\Delta(a) - \Delta_f(a) :=\Delta(a) - (a \otimes 1 + 1 \otimes a)$. This framework do not disturb
the primitive elements of $C$.

Part 6 studies the consequences of some ideas from part 5 to any (Markov) $L$-coalgebra. We generalise
the notion of $L$-coalgebra to $L$-coalgebra of degree $n$.
For instance by creating virtual petals on any (Markov) $L$-coalgebra, we show that
any  (Markov) $L$-bialgebra of degree $n>1$ can be equipped with two Ito-derivatives. We show also that
coproducts of a $L$-bialgebra defines
Ito differentials or a Leibnitz differentials. Two examples of Markov $L$-Hopf algebras of degree 2
end this part, the quaternion one
and $M_2(k)$ via the Pauli matrices.

In part 3, we notice that the usual iterates of a completely positive map come from combinatorics
on a complete graph equipped with a coproduct \footnote{The idea to view coproduct, coming from a
Hopf-algebra, to deal with combinatorics is present in the works of Rota,
\cite{Rota}.} (i.e. a $L$-coalgebra). We generalise this setting
by defining a new composition law, coming from a  $L$-coalgebra. The
iterates, computed with this new law, of a given completely positive map, will still give a semigroup of completely positive
maps.

Part 7 ends this article by showing that any coproduct
can embed any polynomial vector space
into a polynomial algebra. This observation comes from the graph
of $Sl(2)_q$ which embeds a polynomial vector space into an algebra isomorphic to $M_2(k)$.
The generalisation of this idea to any graph, and particularly equipped with a probability
measure leads to the natural concept of probabilistic algebraic product and the notion of random polynomial algebra.
\newpage
\section{Graphs and coassociative coalgebras}
\subsection{Definitions and notation}
\normalsize \textsc{ }
We denote by $k$, the field $\mathbb{R}$ or $\mathbb{C}$ and
consider only unital algebras.
We will use frequently the following symbolic notation: $\sum_{a} a_1 \otimes a_2 \equiv a_1 \otimes a_2 $. For
example, $\Delta(a) = a \otimes a + b \otimes c $ will be denoted by $\Delta(a) = a_{1} \otimes a_{2}$. Moreover all mappings
considered here are linear.
\begin{defi}{}
For every $n>1$, let $ V_1, V_2, ... ,V_n $ be $n$-vector spaces over $k$, we define the {\it{transposition map}} $\tau$ by,
$$V_1 \otimes V_2 \otimes ... \otimes V_n \xrightarrow{\tau} V_n \otimes V_1 \otimes ... \otimes V_{n-1}$$
$$x_1 \otimes x_2 \otimes ... \otimes x_n \mapsto x_n \otimes x_1 \otimes ... \otimes x_{n-1}.$$
\end{defi}
\begin{defi}{}
A $k$-vector space $(A,m,\eta)$  equipped with a product $m: A \otimes A \xrightarrow{} A$ verifying
$m(m \otimes id) =m(id \otimes m)$ (associativity) and a unit map $\eta: k \xrightarrow{} A,  \ \lambda \mapsto \lambda 1_A$
is called an {\it{unital associative algebra}}.
\end{defi}
By reversing the arrows of the previous definition we obtain:
\begin{defi}{}
A {\it{coassociative coalgebra}} over $k$ \cite{Sweedler} \cite{Majid} is an object $(C, \Delta, \epsilon )$
such that $C$ is a $k$-vector space. The counit map $\epsilon: C \xrightarrow{} k$ and the coproduct map $\Delta: C \xrightarrow{} C \otimes C$ verify:
\begin{enumerate}
\item{The coassociativity equation: $ (\Delta \otimes id)\Delta = (id \otimes \Delta)\Delta $.}
\item{The counit equation: $ (id \otimes \epsilon)\Delta = id = (\epsilon \otimes id)\Delta$.}
\end{enumerate}
\end{defi}
\Rk
Coassociativity means that
\begin{equation*}
\begin{CD}
C @>\Delta>> C^{ \otimes 2} \\
@V{\Delta}VV		@VV{\Delta \otimes id}V \\
C^{ \otimes 2} @>{id \otimes \Delta }>> C^{ \otimes 3}
\end{CD}
\end{equation*}
is commutative.
\begin{defi}{}
A \textit{bialgebra} $( C, m, \eta, \Delta, \epsilon, k )$ over $k$ is a vector space such that
$(C, \Delta, \epsilon)$ is a coalgebra and $(C,m, \eta)$ is
an algebra. Moreover the coproduct and counit are unital homomorphisms, i.e.
for all $x,y \in C$,
\begin{eqnarray*}
\epsilon(x)\epsilon(y) = \epsilon(xy) \ &\textrm{and}& \ \epsilon(1_C) = 1_k, \\
\Delta(x)\Delta(y) = \Delta(xy) \ &\textrm{and}& \ \Delta(1_C) = 1_C \otimes 1_C,
\end{eqnarray*}
\end{defi}
The homomorphism property means that the following diagram
\[
\begin{array}{c@{\hskip 1cm}c@{\hskip 1cm}c}
\rnode{a}{C \otimes C} & \rnode{b}{C} & \rnode{c}{C \otimes C}\\[1cm]
\rnode{d}{C \otimes C \otimes C \otimes C} & \ & \rnode{e}{C \otimes C \otimes C \otimes C}
\end{array}
\psset{nodesep=3pt}
\everypsbox{\small}
\ncline{->}{a}{b}\Mput[b]{m}
\ncline{->}{b}{c}\Mput[b]{\Delta}
\ncline{->}{a}{d}\Mput[r]{\Delta \otimes \Delta}
\ncline{->}{d}{e}\Mput[b]{id \otimes \tau \otimes id}
\ncline{->}{e}{c}\Mput[l]{m \otimes m}
\]
commutes
\begin{defi}{}
A {\it{Hopf algebra}} $( H, m, \eta, \Delta, \epsilon, S, k )$ is a bialgebra with
an antipode $S: H \xrightarrow{} H$ verifying:
$$m(S \otimes id)\Delta = m(id \otimes S)\Delta = \eta \epsilon.$$
\end{defi}
If an antipode exists, it is unique and is an unital antialgebra map and an anticoalgebra map, i.e. for all $x,y \in H$,
\begin{eqnarray*}
S(xy) &=& S(y)S(x) \ \textrm{and} \ S(1)=1 \ \textrm{(unital antialgebra map)},\\
(S \otimes S) \Delta x &=& \tau \Delta S(x) \ \textrm{(anticoalgebra map)},\\
\epsilon (S(x)) &=& \epsilon (x).
\end{eqnarray*}
\subsection{Emergence of the notion of graph}
\label{s1}
The main idea is to associate with a tensor product $ a \otimes b $, an arrow
$a \xrightarrow{} b $
and with an associative coalgebra a graph. For convenience we recall the definition of a graph.
\begin{defi}{[Oriented Graph]}
An \textit{oriented graph} $G$ is a quadruple \cite{Petritis}, $(G_{0},G_{1},s,t)$ where
$ G_{0}$ is the denumerable vertex set of $G$ and $G_{1}$ its denumerable arrow set. Moreover
$s,t$ are two mappings, $G_{1} \xrightarrow{s,t} G_{0}$, called source and terminus. A vertex $v$
is called a source if $t^{-1}(\{v \})$ is empty or a sink if $s^{-1}(\{v\})$ is  empty. A graph $G$ is
said to be row (resp. locally) finite if for each vertex $v$, $s^{-1}(\{v\})$ (resp. $t^{-1}(\{v\})$) is finite.
\end{defi}
To deal with random walks on $G$, assumed row and locally finite, we need a family of
probability measures, indexed by the
vertex set.
Let $v \in G_{0}$. We define $F_{v}=\{a \in G_{1}, s(a)=v\}$ and a family of probability
vectors on
$ F_{v}$ by:
$$F_{v} \xrightarrow{P_{v}} [0,1] \ \textrm{such that} \ \sum_{a \in F_{v}} P_{v}(a)=1.$$
This probability vector allows to define a family of probability measures on the $\sigma$-algebra
of subsets of $ F_{v}$. A simple random walk evolving on $G$ is a $G_0$-valued Markov process $(X_n)_{n \in \mathbb{N}}$
with transition probabilities $\mathbb{P}(X_{n+1}=u  \vert X_n=v)= P_v(a)$ if $s(a) =v$ and $t(a)=u$ or $0$ otherwise.
($\mathbb{P}$ is here the unique Markov probability measure on $G_0 ^\mathbb{N}$.)
\Rk
The case of non-oriented graphs can be dealt in this framework by imposing that for each arrow
$a \in G_1$, such that $s(a)=u$ and $t(a)=v$, there exists a unique $\bar{a} \in G_1$ with
$s(\bar{a})=v$ and $t(\bar{a})=u$. We then identify $a$ with $\bar{a}$. Should this identification
be omitted the graph is oriented, the condition of existence of $\bar{a}$ meaning that every arrow
has an inverse.
\Rk
\textbf{[Generalisation of the definition of a graph]} As we shall see in part three, the whole graph structure
can be encoded in its coproducts. Therefore, provided that the coproducts are well defined, we can
enlarge the concept of graph. For instance we can omit the denumerability condition in the definition of a graph
given above.
\begin{defi}{}
Let $C$ be a coassociative coalgebra. We denote $G(C)$ the graph associated with $C$. Its vertex set is
$C$ and its arrow set, the set of those tensor products appearing in the definition of the coproduct.
\end{defi}
Let see two examples \cite{Majid}:
\begin{exam}{}
We consider the coalgebra generated by $C=( 1, X, g, g^{-1} )$ with the following relations:
$$
\Delta X = X \otimes 1 + g \otimes X, (\textrm{thus there will be two arrows:} \ X \xrightarrow{} 1 \ \textrm{and} \ g \xrightarrow{} X),  \ \
\Delta g = g\otimes g, \ \
\Delta 1 = 1\otimes 1. \ \
$$
$$
\epsilon (X) = 0, \ \
\epsilon( g) = 1, \ \
\epsilon( 1) = 1.
$$
With this coalgebra structure we associate the following graph (left corner):
\begin{center}
\includegraphics*[width=5cm]{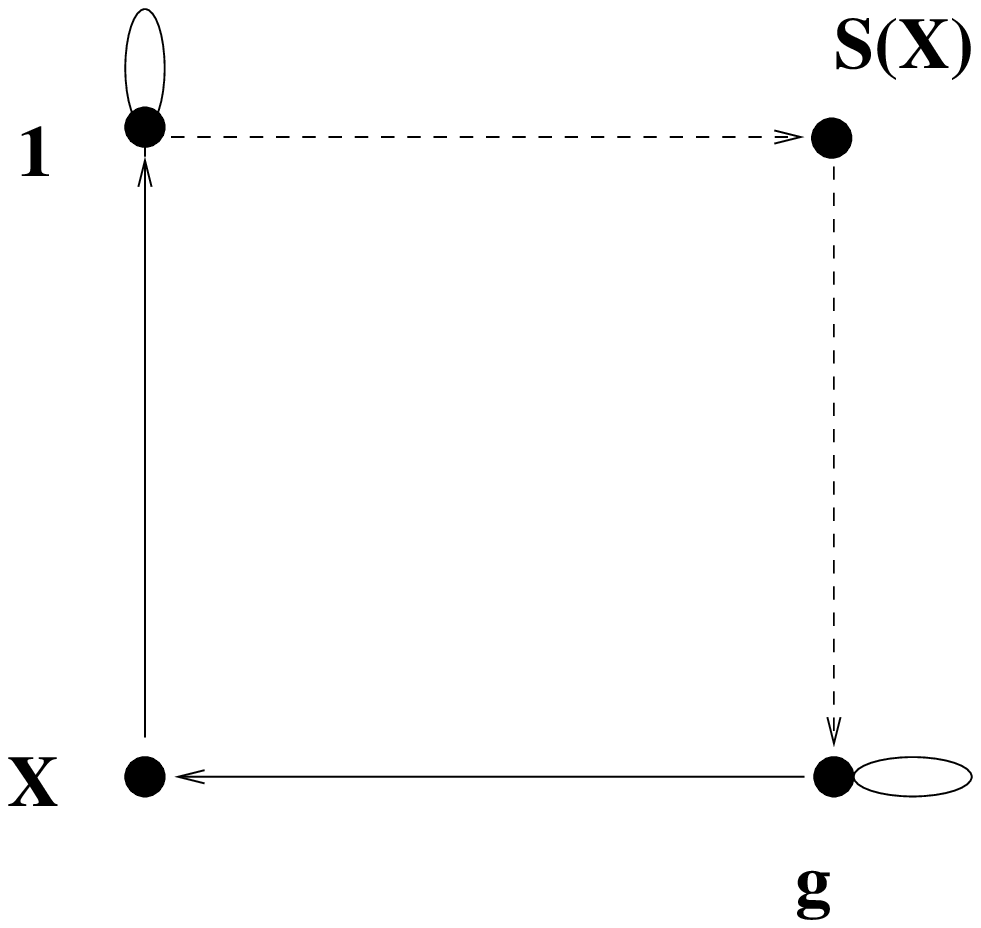}
\end{center}
The following algebraic relations turn $C$ into a bialgebra.
$$
gg^{-1} = 1 = g^{-1}g, \ \
Xg = qgX, \ \
Xg^{-1} = q^{-1}g^{-1}X.
$$
where $q$ is a fixed invertible element of the field $k$.
The following antipode structure turns $C$ into a Hopf algebra.
$$
S(X) = -g^{-1}X, \ \
S(g) = g^{-1}, \ \
S(g^{-1}) = g.
$$
The coproducts associated with the antipode action are:
$$
\Delta g^{-1} = g^{-1} \otimes g^{-1}, \ \
\Delta S(X) = 1 \otimes S(X) + S(X) \otimes S(g), \ \
\epsilon (g^{-1}) = 1, \ \
\epsilon(S(X)) = 0.
$$

The left corner of the picture represents the coalgebra $( 1, X, g ) $ without its antipode part
$( g^{-1}, S(X), 1 )$, represented in the right corner.

We notice that the antipode $S$ describes a kind of folding of the graph. Here $g$ is glued with
its algebraic inverse $g^{-1}=S(g)$, $1=S(1)$ is glued with itself and $X$ is glued with $S(X)$.
Moreover when we fold the manifold as we just described, we observe that the arrow emerging from,
for example $X$, and those emerging from his algebraic gluing $S(X)$ are in opposite sense.
The antipode $S$ realises what we could call an algebraic gluing, a geometrical gluing for the graph
in a compatible way with the algebraic product.
Moreover if we imagine the graph in a dynamical way, we can say that the antipode realises a time reversal.
That is just the physical interpretation of the anticoalgebra map denomination when we see a Hopf algebra as
an oriented graph.
\end{exam}
\begin{exam}{[$Sl(2)_q$]}
We define the well-known coalgebra structure:
$$
\Delta a = a \otimes a + b \otimes c, \ \
\Delta b = a \otimes b + b \otimes d, \ \
\Delta c = d \otimes c + c \otimes a, \ \
\Delta d = d \otimes d + c \otimes b.
$$
Here is the picture of $G(Sl(2)_q)$:
\begin{center}
\includegraphics*[width=5cm]{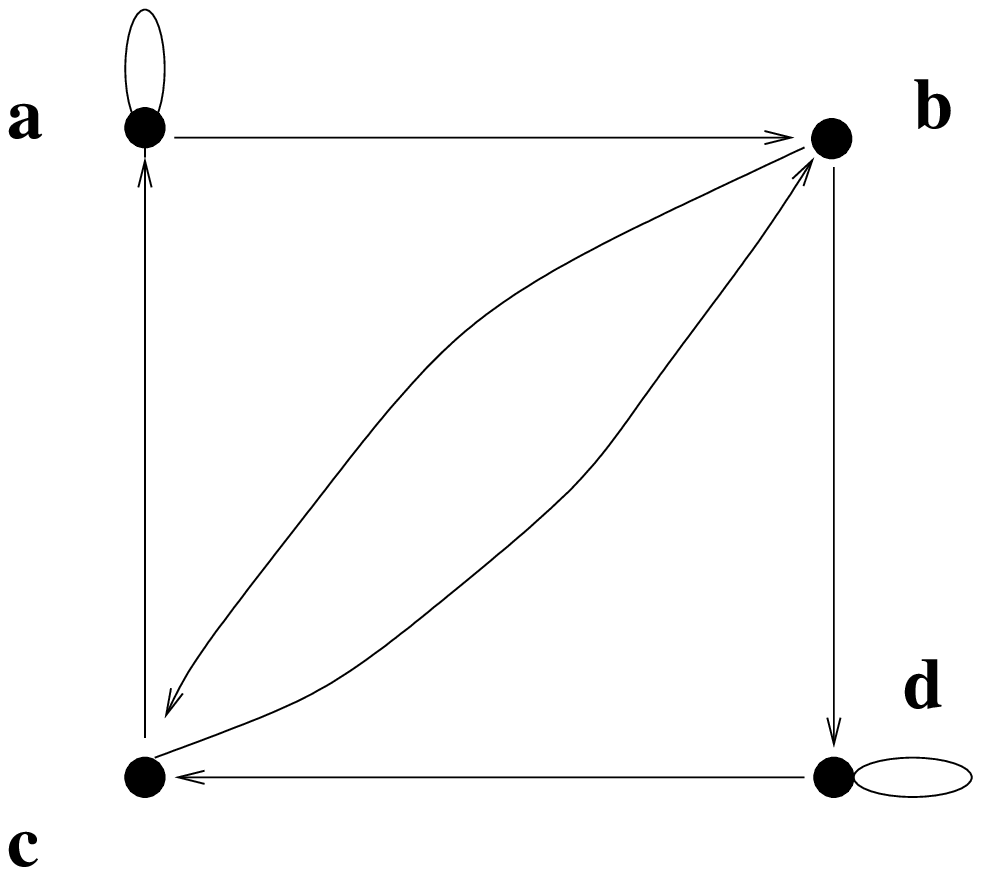}
\end{center}
Another physical relevant remark is the non locality of this coproduct. To assure the coassociativity of
the coproduct we need to delocalise the effect of the coproduct or what a physicist would call
a propagator on the oriented graph $G$.
For the sake of an example, notice that when a signal is emitted from $a$ and evolves using the coproduct
structure $\Delta$, the only physical accessible vertices ---if we take validity of the locality principle
for granted---are either $a$ or $b$.
Here only $a$ must be chosen because $a \otimes b $ is not represented in the definition of the
coproduct of $a$. Nevertheless the information arrives at $a$ but
in a completely delocalised way some information has to start from $b$ to arrive at $c$ without
possible communication from $a$ to $b$! (We recall that here we live in a world such
that the only way to communicate is by $\Delta$.)

Thus to satisfy the mathematical property of associativity of the product $m$,
which becomes coassociativity of a coalgebra when we reverse the direction of the arrows of commutative
diagrams we lose, if we see the coalgebra as an oriented graph and the coproduct as a propagator, the physical
principle of locality.
\end{exam}
\Rk
Let $C=(x_1, \ldots, x_n)$. To produce a coassociative coalgebra with counit from $C$, place $(x_1, \ldots, x_n)$
in an arbitrary matrix $T \in M_{m}(C)$ with $m^2 \geqslant n$. Then $T \bar{\otimes} T$, where $(T \bar{\otimes} T)_{ij} := \sum_k T_{ik} \otimes T_{kj} := \Delta T_{ij}$ will give a coproduct.
For instance the first example of a coassociative coalgebra we met is given by:
\[
 T = \begin{pmatrix}
 1 & 0\\
X & g
\end{pmatrix} , \ \
T \bar{\otimes} T = \begin{pmatrix}
1 \otimes 1 & 0\\
X \otimes 1 + g \otimes X &  g \otimes g
\end{pmatrix}
=\begin{pmatrix}
\Delta 1 & 0\\
\Delta X & \Delta g\
\end{pmatrix}.
\]
Nevertheless there exist
other ways to produce coassociativity all relying on the matrix product. Instead of computing
$T \bar{\otimes} T$, we compute $T \bar{\otimes} U$, where $U$ is a well chosen matrix with coefficients
from
the $(x_1, \ldots, x_n)$.
Let us see an example on the two points case: $C=(a,b)$.
If we choose:
\[ T_1 = \begin{pmatrix}
 a & 0\\
b & 0
\end{pmatrix} ,
\]
and compute $T_1 \bar{\otimes} T_1$. We will obtain
$\Delta(a)= a \otimes a, \ \ \Delta(b)= b \otimes a $ and a possible coassociative
coalgebra structure on $C$. (Another possibility is to choose a loop on each vertex, that is
\[ T_2 = \begin{pmatrix}
 a & 0\\
0 & b
\end{pmatrix},
\]
and compute $T_2 \bar{\otimes} T_2$.)
Yet if we choose:
\[ T_3 = \begin{pmatrix}
 a & b\\
b & 0
\end{pmatrix} , \ \
U_3 = \begin{pmatrix}
 b & 0\\
a & 0
\end{pmatrix},
\]
and compute $T_3 \bar{\otimes} U_3$, we find that the coproduct defined by $\Delta(a)= a \otimes b + b \otimes a, \ \ \Delta(b)= b \otimes b $
is coassociative and has still a counit given by $\epsilon(a)=0, \ \epsilon(b)=1$.
Similarly,
\[ T_4 = \begin{pmatrix}
 a & a\\
b & b
\end{pmatrix} , \ \
U_4 = \begin{pmatrix}
 b & 0\\
a & 0
\end{pmatrix},
\]
give a coproduct defined by $\Delta(a)= a \otimes b + a \otimes a, \ \ \Delta(b)= b \otimes b + b \otimes a$.
Yet it has no counit but only a right counit, given by $\epsilon(a)=1, \ \ \epsilon(b)=0$ and thus verifying $(id \otimes \epsilon)\Delta = id$. No left counit exists.
Another coassociative coproduct on $C$ can be found by computing $T_5 \bar{\otimes} U_5$ with,
\[ T_5 = \begin{pmatrix}
 a & b\\
b & a
\end{pmatrix} , \ \
U_5 = \begin{pmatrix}
 b & 0\\
a & 0
\end{pmatrix}.
\]
We have $\Delta(a) = a \otimes a + b \otimes b, \ \Delta(b) = a \otimes b + b \otimes a $. The
counit is defined by $\epsilon(a)=1,  \ \epsilon(b)=0$. The graph associated with this coalgebra
is the same as the example just before. (Remark that $\Delta_1(b) \equiv \Delta(a), \ \epsilon_1(b)=1$ and $\Delta_1(a) \equiv \Delta(b), \ \epsilon_1(a)=0$ works also.)
Interesting questions arise. Let $G$ be a graph. Does it exist a method to know if a non local
structure, that is a coassociative coproduct can be put on its vertex set, and how many such non local structures exist?
Is the matrix product the only way to obtain such coassociative coproduct ?
In the following we do not address these questions. However when we wish to put a coassociative coalgebra
structure on a $k$-vector space $C$ we shall use a matrix $T$ and compute $T \bar{\otimes}T$. In this case
there is unicity between the graph and the coassociative coproduct on it.
We stop here considerations on these open questions by the following theorem and its physical meaning.
\begin{defi}{[Groupoid]}
A {\it{groupoid}} \cite{Connes} consists of a set $G$, a distinguished subset $G^{(0)} \subset G$, two maps
$t,s: G \xrightarrow{} G^{(0)} $ and a law of composition
$$ \circ: G^{(2)}=\{(\gamma_1,\gamma_2) \in G \times G; \ s(\gamma_1)=t(\gamma_2) \} \xrightarrow{} G$$
such that:
\begin{enumerate}
\item{$s(\gamma_1 \circ \gamma_2) = s(\gamma_2), \ t(\gamma_1 \circ \gamma_2) = t(\gamma_1), \ \forall(\gamma_1, \gamma_2) \in G^{(2)},$}
\item{$\forall x \in G^{(0)}, \ s(x)=t(x)=x; \ \forall \gamma \in G, \ \gamma \circ s(\gamma) = \gamma, \ t(\gamma) \circ \gamma = \gamma, $}
\item{$(\gamma_1 \circ \gamma_2) \circ \gamma_3 = \gamma_1 \circ (\gamma_2 \circ \gamma_3), $}
\item{Each $\gamma$ has a two-sided inverse $\gamma^{-1}$, with $ \gamma\gamma^{-1}=t(\gamma), \ \gamma^{-1}\gamma= s(\gamma) .$}
\end{enumerate}
\end{defi}
\begin{theo}
Let $G$ be a set verifying all the properties of the definition of a groupoid, but maybe the last one \footnote{For instance an oriented graph, or the spectrum of an atom.}.
Let $C$ be a $k$-vector space. We consider
$\tilde{C} = \{ c_\gamma= \alpha(c,\gamma), \gamma \in G \}$ the $k$-vector space obtained by an action
$\alpha: C \times G \xrightarrow{} C$. Then $\tilde{C}$ can be embbeded into a coassociative
coalgebra.
\end{theo}
\Proof
Fix $c_\gamma \in \tilde{C} $ and define $\Delta(c_\gamma) = \sum_{\gamma_1 \circ \gamma_2 = \gamma} c_{\gamma_1} \otimes c_{\gamma_2}.$
We observe that any $\gamma \in G$ can be decomposed in three parts, i.e. it exists $\gamma_i$ such that $\gamma = \gamma_1 \circ \gamma_2 \circ \gamma_3$.
By definition of the coproduct and the associativity of the product $\circ$ we get,
\begin{eqnarray*}
\sum_{\gamma_1 \circ \gamma_2 = \gamma} c_{\gamma_1} \otimes \Delta c_{\gamma_2} &=& \sum_{\gamma_1 \circ \gamma_2 = \gamma} \sum_{\gamma_1' \circ \gamma_2' = \gamma_2} c_{\gamma_1} \otimes (c_{\gamma_1'} \otimes c_{\gamma_2'})
= \sum_{\gamma_1 \circ (\gamma_1' \circ \gamma_2') = \gamma} c_{\gamma_1} \otimes (c_{\gamma_1'} \otimes c_{\gamma_2'}),  \\
\sum_{\gamma_1 \circ \gamma_2 = \gamma} \Delta c_{\gamma_1} \otimes c_{\gamma_2} &=& \sum_{\gamma_1 \circ \gamma_2 = \gamma} \sum_{\gamma_1'' \circ \gamma_2'' = \gamma_1}
(c_{\gamma_1''} \otimes c_{\gamma_2''}) \otimes c_{\gamma_2}= \sum_{(\gamma_1'' \circ \gamma_2'') \circ \gamma_2 = \gamma} (c_{\gamma_1''} \otimes c_{\gamma_2''}) \otimes c_{\gamma_2},
\end{eqnarray*}
proving that $(id \otimes \Delta)\Delta=(\Delta \otimes id)\Delta$ since the sums involved are over all possible decompositions of $\gamma$ in three parts.
As $\gamma = \gamma \circ s(\gamma)=t(\gamma) \circ \gamma$, we define $\epsilon(\gamma)= 0$, if $\gamma \in G \setminus G^{(0)}$
and $\epsilon(\gamma)= 1$ otherwise, i.e. if $\gamma \in G^{(0)} $. We have
$$\Delta(c_{\gamma}) = \sum_{\gamma_1 \circ \gamma_2 = \gamma; \ \gamma \in G \setminus G^{(0)}} c_{\gamma_1} \otimes c_{\gamma_2} + c_{t(\gamma)} \otimes c_{\gamma} + c_{\gamma} \otimes c_{s(\gamma)},$$
thus $(id \otimes \epsilon)\Delta = (\epsilon \otimes id)\Delta = id$.
\eproof
\Rk
Let $G$ be a groupoid.
Let $C$ be an associative algebra. We consider $\tilde{C} = \{ c_\gamma= \alpha(c,\gamma), \gamma \in G \}$, as defined above.
Then $\tilde{C}$ can be embbeded into an associative algebra by
defining the convolution product as follows:
$$(a*b)_{\gamma}= \sum_{\gamma_1 \circ \gamma_2 = \gamma} a_{\gamma_1} b_{\gamma_2}, \ \forall a,b \in C.$$
As $\circ$ is associative, so is the convolution product. $(\tilde{C}, *)$ is an associative algebra.
It is not in general a bialgebra.
\begin{exam}{}
Let $T$ be a matrix. We recover the coassociativity of the $\bar{\otimes}$ notation by noticing that
$\Delta T_{(ij)} = \sum_{(ik) \circ (kj) = (ij)} T_{(ik)} \otimes T_{(kj)}$.
\end{exam}
\Rk
If we take for granted that some graphs, i.e. some models of space-time, equipped
with a coassociative coproduct, behave in a non local way, we have to find a process which
describes the phase transition between physics governed by the locality principle and physics
governed by the coassociative principle. An attempt to explain this phase transition could lie in
the choice of the (semi)-groupoid $G$, of the previous theorem. If $G$ is a graph and acts on
itself we obtain a coassociative coproduct on the path algebra of $G$. We remain in the
realm of locality because paths are macroscopic objects. The coproduct can be seen as a
deconcatenation of the path. If $G$ acts on vertices of a graph, seen as microscopic objects,
via a matrix $T$ and the
$\bar{\otimes}$ product, we obtain non local coproducts.
\begin{defi}{}
A coalgebra is {\it{cocommutative}} if $\Delta = \tau \Delta$ .
\end{defi}
\begin{prop}
If C is a cocommutative coalgebra then its graph $G(C)$ will be non-oriented.
\end{prop}
\Proof
Let $a,b,x \in C$ and suppose that the term $ a \otimes b $ appears in the description of $\Delta x$. The same must be true
for $ b \otimes a $, since $\Delta = \tau \Delta$.
There is an arrow emerging from $a$ to $b$ and from $b$ to $a$. We have just proved that the graph
is bi-oriented. By identifying the arrow emerging from $a$ to $b$ to that one from $b$ to $a$, we
obtain a non-oriented graph.
\eproof

We have seen that we can associate a graph with an associative coalgebra. It is tempting to try to
associate with a graph an associative coalgebra \footnote{There exists such a construction, but defined
on the paths of a graph and not on its vertex set \cite{rosso}.}.
We equipped $G_0$ (identified with $G$ in the following) with its free $k$-vector space. In the case of a
coassociative coalgebra $C$ we notice that $\Delta$ is defined on the vertex of the graph and maps
$G(C)$ to $ G(C) \otimes G(C) $. (In the following we identify $C$ to $G(C)$.)
Can we associate a coproduct on the vertex of a graph $G$ such that the algebraic setting so constucted
defines a coassociative coalgebra ?
The answer is no in general because as we have seen, the coassociative coproduct is non local.
Nevertheless we can find such a coalgebra structure. The price to pay to restore locality is
to break the coassociativity and restore a kind of coassociativity by introducing two propagators \cite{Lerdresden}:
\begin{enumerate}
\item{the right coproduct $\Delta $:  $ G \xrightarrow{} G \otimes G ,$ }
\item{the left coproduct  $\tilde{\Delta}$: $ G \xrightarrow{} G \otimes G .$ }
\end{enumerate}
These propagators will verify the coassociativity breaking equation:
\begin{equation}
\label{coabreak}
(\tilde{\Delta} \otimes id)\Delta = (id \otimes \Delta)\tilde{\Delta},
\end{equation}
instead of the coassociativity equation:
\begin{equation}
\label{coa}
(\Delta \otimes id)\Delta = (id \otimes \Delta)\Delta.
\end{equation}

We postpone a speculative physical interpretation of the meaning of the two equations in part 3
and focus on the restoration of the locality.

\section{Coassociativity breaking}
\subsection{The paths of a random walk on an oriented graph: an algebraic framework}
\label{troisun}
Before expressing the algebraic setting, we start with the example of a graph $G$.
To include random walks on $G$ we allowed the arrows to carry a probability.
We define for a graph $G$ with no sinks and no sources:
\begin{enumerate}
\label{rdd}
\item{ The right coproduct (the notion of future): $\Delta: G_{0} \xrightarrow{} G_{0} \otimes G_{0}
$
$$
\label{eqdr}
v \mapsto \Delta(v)=\sum_{ a \in G_{1} \atop s(a) = v} P_{v}(a) \ v \otimes t(a).
$$
    }
\item{ The left coproduct (the notion of past): $\tilde{\Delta}: G_{0} \xrightarrow{\tilde{\Delta}} G_{0} \otimes G_{0}
$
$$
\label{eqdl}
v \mapsto \tilde{\Delta}(v)= \dfrac{1}{^\#t^{-1}(\{v\})} \sum_{ a \in G_{1} \atop  t(a) = v}  s(a) \otimes v.
$$
     }
\item{ The right counit: $\label{eqcr} \epsilon: G_{0} \xrightarrow{\epsilon} k, \ \ v \mapsto 1.$}
\item{ The left counit: $\label{eqcl} \tilde{\epsilon}: G_{0} \xrightarrow{\tilde{\epsilon}} k, \ \ v \mapsto 1.$}
\end{enumerate}
\Rk
Here there is a dissymmetry in the definition of the right and left coproduct because what we consider as
relevant for a vertex $v$ is the probability associated with the arrows emerging from $v$,
i.e.\ its future and not the probability associated with an arrow arriving at $v$ which only depends on
another vertex.
\Rk
Recall that we wish to put a graph on mathematical structures such as algebras,
coalgebras and so forth, that is we seek for a geometrical link between elements belonging to such
structures. Thus if the graph has sinks and loops, we connect the identity element of the structure
to be considered as the future (resp.\ past) of a sink (resp.\ source) and
place a loop on the identity element, $I$, that is we impose:
$\tilde{\Delta} I = I \otimes I, \ \ \Delta I = I \otimes I.$
\Rk
With this setting the coassociativity breaking equation:
$(\tilde{\Delta} \otimes 1)\Delta = (1 \otimes \Delta)\tilde{\Delta}$
is realised \footnote{If at Planck scale, space time is a graph, the notion of time has to emerge from the very definition of the graph, i.e.
from the arrows. As arrows are related to coproducts, that is to dynamics
the space $G$ and its coproduct(s), (i.e. a $L$-coalgebra), carry themselves
a notion of space-time.

Equations \ref{coabreak} and \ref{coa} have then an interesting and speculative physical interpretation.
They do not introduce a notion of time, already carried by the very definition of the graph itself, but
a notion of future, represented by $\Delta$ and a notion of past,
represented by $\tilde{\Delta}$ in the first case and not in the second case, because $\Delta=\tilde{\Delta}$.
Thus a coassociative graph means a graph $G$ with dynamics $\Delta$ obeying equation \ref{coa},
where the human notion of future and past do not exist.
On the contrary, what we shall call a Markov $L$-coalgebra, that is a graph $G$ with its dynamics
induced by $\Delta$ and $\tilde{\Delta}$, carry naturally such notions, which emerge by considering
only the physical concept of locality.

To summarize, locality concept implies notions of future and past. Coassociativity
concept (non locality) implies the lost of these notions.

We suggest an interpretation of quantum entanglement, a subject of actuality in quantum information
theory. Would it be possible that when two particles are entangled the notion of future and past
disappear (not the notion of time as we have just seen)?

If the answer is positive we would avoid all the puzzling questions about the possible transmission of
information, and at what speed, (because there is not). We would understand as well why such particles are
in perfect correlations because correlations means the preexistence of a notion of future and past,
inexisting in the world of two such entangled particles.

As for the notion of causality, another important concept of physics, we notice that it is intrinsic
to the notion of dynamics.}.
\subsection{Axioms}
\label{troisdeux}
\begin{defi}{[$L$-coalgebra]}
A \textit{$L$-coalgebra} with counits is a vector space over $k$ equipped with its right and left coproducts $(C, \Delta, \tilde{\Delta}, \epsilon, \tilde{\epsilon})$
such that the following equations hold:
\begin{enumerate}
\item{The coassociativity breaking equation: $ (\tilde{\Delta} \otimes id)\Delta = (id \otimes \Delta)\tilde{\Delta}.$}
\item{The right counit equation: $ (id \otimes \epsilon)\Delta = id.$}
\item{The left counit equation: $ ( \tilde{\epsilon} \otimes id)\tilde{\Delta} = id. $}
\end{enumerate}
\end{defi}
\Rk
The coassociativity breaking equation means that
\begin{equation*}
        \begin{CD}
        C @>\Delta>> C^{ \otimes 2} \\
        @V{\tilde{\Delta}}VV		@VV{\tilde{\Delta} \otimes id}V \\
        C^{ \otimes 2} @>{id \otimes \Delta }>> C^{ \otimes 3}
        \end{CD}
        \end{equation*}
is commutative.
The right counit equation means that
\begin{equation*}
        \begin{CD}
        C @>\Delta>> C^{ \otimes 2} \\
	@V{id}VV		@VV{ id \otimes \epsilon}V \\
        C @= C \otimes k
        \end{CD}
        \end{equation*}
is commutative, and the left counit equation means that
\begin{equation*}
        \begin{CD}
        C @>\tilde{\Delta}>> C^{ \otimes 2} \\
	@V{id}VV		@VV{ \tilde{\epsilon} \otimes id}V \\
        C @= k \otimes C
        \end{CD}
        \end{equation*}
is commutative.

To discriminate between the different types of $L$-coalgebras, we define:
\begin{defi}{[Markov $L$-coalgebra]}
A {\it{Markov $L$-coalgebra}} $C$ is a $L$-coalgebra such that for all $x \in C $, $\Delta x = x \otimes x_1 $
and $ \tilde{\Delta}x = x_0 \otimes x $ for some $x_0,x_1 \in C$, (in symbolic notation).
\end{defi}
Such a structure reproduces locally what we have in mind when we speak about random walks on
a graph.
\begin{theo}
We can always associate with an oriented graph $G$, equipped with a probability structure,
a unique Markov $L$-coalgebra with counits.
\end{theo}
\Proof
The proof is carried by verifying the items \ref{eqdr}, \ref{eqdl}, \ref{eqcr}, \ref{eqcl} of
subsection \ref{troisun}.
For the bijection part, we note that all the geometric and probabilistic information are coded
into the coproducts. Therefore two different graphs lead necessarily to two different coproducts and conversely.
\eproof \\
\Rk
If there is not a probability structure, either we choose the equiprobability, i.e.\ for each $a \in
s^{-1}(\{v\}), \ \ P_v(a) = \dfrac{1}{^\#s^{-1}(\{v\})} $
or we choose to affect 1 to each arrow. In this case we have a $L$-coalgebra without right counit.
However, since we are more interested in coproducts than counits, we will choose the second
choice and suppress the term $\dfrac{1}{^\#t^{-1}(\{v\})}$ of the definition of the left coproduct as well.

As we associate with each tensor product an arrow, each $L$-coalgebra will be represented by a unique graph, that which is defined
from the coproducts. However a graph, as a geometric object, can be a support for several $L$-coalgebras. For instance the graph
associated with $Sl(2)_q$ is the support for $Sl(2)_q$, as a coassociative coalgebra and could be the support as well for the Markov
$L$-coalgebra, without counits, naturally associated with it.
\begin{defi}{}
We define the set of the coproducts (propagators) of a graph $G$ by
$\textsf{Prop}(G)=\{(\Delta,\tilde{\Delta}); \Delta, \\ \tilde{\Delta}:  G \xrightarrow{} G^{\otimes 2};  (\tilde{\Delta} \otimes id)\Delta =
(id \otimes \Delta)\tilde{\Delta}\}$.
\end{defi}
\begin{prop}
The co-products of a Markov $L$-coalgebra are the generators of all possible (random) walks.
\end{prop}
\Proof
The sequence
$ \Delta_1 \equiv \Delta, \Delta_2 = 1 \otimes \Delta, \Delta_3 = 1 \otimes 1 \otimes \Delta, \ldots $
generates all possible (random) walks starting at any vertex. Similarly, The sequence of powers of
$\tilde{\Delta}$, generates all the possible (random) walks arriving at a given vertex.
\eproof
\begin{defi}{[the $\star$ product]}
The $\star$ product is defined by $\Delta_n \star \Delta_m = \Delta_{n+m}$.
\end{defi}
\Rk
The degenerate case $(\Delta=\tilde{\Delta}, \ \epsilon=\tilde{\epsilon})$ corresponds
to a coassociative coalgebra.
\Rk
The previous proposition allows us to generalise the concept of walk on a graph to any
$L$-coalgebra because only the coproducts are involved. As an example the walk of say, a wave
function from $a$, vertex of the graph associated with $Sl(2)_q$ is described by the sequence
$(a, \Delta(a), \Delta_2(a), \Delta_3(a), \ldots)$. Obviously this walk is a deterministic (and
non local) walk
because no notion of probability is used in this coassociative coalgebra.
\begin{defi}{[$L$-bialgebra]}
$C$ is a \textit{$L$-bialgebra} (with counits) over $k$ if it is a $L$-coalgebra (with counits) over $k$ and an unital algebra.
Moreover its coproducts and counits are
homomorphisms.
\end{defi}
Here we choose to study only $L$-bialgebra associated with an associative product $m$.
\begin{defi}{[$L$-Hopf algebra]}
$C$ is a \textit{$L$-Hopf algebra} (with counits) if it is a $L$-bialgebra (with counits) equipped with right and left antipodes
$S, \ \tilde{S}: C \xrightarrow{} C $, such that:
$$
m(id \otimes S)\Delta = \eta \epsilon, \ \
m(\tilde{S} \otimes id)\tilde{\Delta} = \eta \tilde{\epsilon}.
$$
\end{defi}
\Rk
This structure seems to play a less important r\^ole than in the coassociative case. We notice that
the r\^ole of the antipode $S$ in a Hopf algebra is to produce both a gluing of the graph and a
reversal of the arrows. Therefore it acts as a time reversal in a compatible way with the underlying algebraic
structure. However we have:
\begin{prop}
Let $G$ a Markov $L$-coalgebra. If the arrows of $G$ are reversed, then the graph so obtained is still
a Markov $L$-coalgebra.
\end{prop}
\Proof
Let $G_r$ be the graph $G$ whose arrows have been reversed. By denoting $\Delta_r = \tau \tilde{\Delta}$,
$\tilde{\Delta}_r = \tau \Delta $, $ \epsilon_r = \tilde{\epsilon} $, $ \tilde{\epsilon}_r = \epsilon $,
its coproducts and counits, the proof is now obvious.
\eproof
\subsubsection{About $L$-cocommutativity}
Let $C$ be a coassociative coalgebra with coproduct $\Delta$.
In \cite{Quillen} Quillen introduced in the case of a bicomodule $M$ over $C$, that is
a vector space equipped with left and right comodule structures which commute, i.e.
$(\Delta_l \otimes id)\Delta_r = (id \otimes \Delta_r)\Delta_l$ the notion of cocommutator subspace
defined by:
$$ M^{\natural} = \ker \{\Delta_l - \tau \Delta_r : M \xrightarrow{} C \otimes M \}$$
The adaptation of this concept in the case of a Markov $L$-coalgebra $G$ gives the following
interesting result:
\begin{defi}{}
We denote $ G^{\natural} = \ker \{\Delta - \tau \tilde{\Delta} : G \xrightarrow{} G \otimes G \}$
the \textit{$L$-cocommutator} subspace of $G$.
\end{defi}
\Rk
If we do not consider probability
on the Markov $L$-coalgebra $G$,
the $L$-cocommutator subspace of $G$ is the set of vertices $v$ of the graph whose the arrows emerging and arriving at $v$ can be identified
to give locally a non oriented graph.
\subsection{Examples of non degenerate $L$-coalgebras}
\begin{exam}{[The Cuntz-Krieger algebra]}
In \cite{CK} \cite{Pask} a $C^{*}$-algebra, called a {\it{Cuntz-Krieger algebra}}, is associated with a graph $G$.
If $G$ is a row-finite (i.e. $\forall v \in G_0, \ s^{-1}(\{v\})$ is finite) oriented graph, a Cuntz-Krieger
$G$-family consists of a set $\{ P_v : v \in G_0 \}$ of mutually orthogonal projections and
a set $\{ S_e : e \in G_1 \}$ of partial isometries satisfying,
$$ \forall (e,v) \in G_1 \times G_0, \  S ^{*}_e  S^{}_e=P_{t(e)}, \  P_v=\sum_{e: s(e)=v} S^{}_e S^*_e .$$
\begin{prop}
A Cuntz-Krieger algebra $CK$ associated with a graph
without sinks and loops and whose vertex set is finite is a Markov $L$-bialgebra.
\end{prop}
\Proof
By defining
$$\Delta(P_v) = \sum_{v_1 \in t(s^{-1}(\{v\}))} P_v \otimes P_{v_1} \ \textrm{and} \
\tilde{\Delta}(P_v) = \sum_{v_0 \in s(t^{-1}(\{v\}))} P_{v_0} \otimes P_v,$$ we embed $CK$ into
a Markov $L$-coalgebra. There is no counit.
Thanks to the mutual orthogonality property of the projectors, we get for instance,
$$\Delta(P_v) \Delta(P_{v'}) = \sum_{v_1 \in t(s^{-1}(v))}\sum_{v_1' \in t(s^{-1}(v'))}
P_v P_{v'} \otimes P_{v_1}P_{v_1'} = \delta(v,v') \Delta(P_v) = \Delta(P_v P_{v'}),$$ where
$\delta$ is the Kronecker symbol.
Since the vertex set is finite, $CK$ has an identity element $\sum_{v \in G_0} P_v :=I$. In general
we do not have $\Delta(I) = I \otimes I$.
\eproof
\end{exam}
\begin{exam}{[Unital algebra]}
Let $A$ be a unital algebra. $A$ carries a non trivial $L$-bialgebra called the flower graph
with coproducts $\delta(a) = a \otimes 1$ and $\tilde{\delta}= 1 \otimes a$.
It is the subject of part 5.
\end{exam}
\begin{exam}{[The degenerate case]}
Any coassociative coalgebra $C$ with coproduct $\Delta$ can be embedded into a $L$-coalgebra structure
in a non degenerate way.

A generalisation of an idea, applied by R.L. Hudson \cite{Hudson1} \cite{Hudson2}, used in the case of a $C$-bi-comodule
where $C$ is a coassociative coalgebra, into the present setting is given below. We define:
\begin{enumerate}
\item{$ C \xrightarrow{\overrightarrow{d}} C \otimes C \ \textrm{such that} \ \overrightarrow{d}(x) = \Delta x -x \otimes 1.  $
 }
\item{$ C \xrightarrow{\overleftarrow{d}} C \otimes C \ \textrm{such that} \ \overleftarrow{d}(x) =  \Delta x - 1 \otimes x.  $
       }
\end{enumerate}
\begin{theo}[Hudson \cite{Hudson2}]
$\overleftarrow{d},\overrightarrow{d}$, turns a coassociative coalgebra into a non degenerate $L$-coalgebra.
\end{theo}
\Proof
With the notation above we must prove the coassociativity breaking equation. In adopting symbolic notation
fix $x \in C$ such that $\Delta x = x_1 \otimes x_2$.
$$ x \xrightarrow{\overleftarrow{d}} x_1 \otimes x_2 - 1 \otimes x  \xrightarrow{id \otimes \overrightarrow{d}}
x_1 \otimes \Delta x_2 - \Delta x \otimes 1 - 1 \otimes \Delta x +  1 \otimes  x \otimes 1, $$
$$ x \xrightarrow{\overrightarrow{d}} x_1 \otimes x_2 - x \otimes 1  \xrightarrow{\overleftarrow{d} \otimes id}
\Delta x_1 \otimes x_2 - \Delta x \otimes 1 - 1 \otimes \Delta x +  1 \otimes x \otimes 1. $$
Hence we prove that the following diagram,
\begin{equation*}
        \begin{CD}
        C @>\overrightarrow{d}>> C^{ \otimes 2} \\
        @V{\overleftarrow{d}}VV		@VV{\overleftarrow{d} \otimes 1}V \\
        C^{ \otimes 2} @>{1 \otimes \overrightarrow{d} }>> C^{ \otimes 3}
        \end{CD}
\end{equation*}
commutes.
Moreover $ \overleftarrow{d}$ is obviously not equal to $\overrightarrow{d}$ on the whole coalgebra.
\eproof \\
Here is an example of an $L$-coalgebra which in general cannot be thought in terms of a Markov $L$-coalgebra.
A priori there would not exist counits.
\begin{defi}{[Ito $L$-coalgebra]}
Let $C$ be a coassociative coalgebra over $k$. We call $(C, \overleftarrow{d}, \ \overrightarrow{d})$, an Ito $L$-coalgebra \footnote{
Here appears, in this non local setting, two operators $ (\overleftarrow{d}, \ \overrightarrow{d})$ which could be interpreted as
a generalisation of notion of past and future.}
over $k$.
\end{defi}
Let $C$ be a coassociative coalgebra. As we will see in part 5, the coproducts of the flower graph, $\delta(x) := x \otimes 1$ and $\tilde{\delta}(x) := 1 \otimes x$ for all
$x \in C$, carry a natural chirality. That is why we define:
\begin{defi}{[Chiral bimodule]}
Let $C$ be a coassociative coalgebra. As we have seen, $(C, \overleftarrow{d},\overrightarrow{d})$ is an Ito $L$-coalgebra.
We embed $C ^{\otimes 2}$
into a {\it{$C$-chiral bimodule}} by defining only on the image of the
coproducts of the Ito $L$-coalgebra, $\overleftarrow{d}, \overrightarrow{d}$ the following products: \\
Let $c, x, y \in C$,
\begin{eqnarray*}
x.\overleftarrow{d}(c) &=& \tilde{\delta}(x)\overleftarrow{d}(c); \ \ \ \ \overleftarrow{d}(c).y = \overleftarrow{d}(c)\tilde{\delta}(y), \\
x.\overrightarrow{d}(c) &=& \delta(x)\overrightarrow{d}(c); \ \ \ \ \overrightarrow{d}(c).y = \overrightarrow{d}(c)\delta(y).
\end{eqnarray*}
\Rk
The notion of chiral bimodule is well defined because if $\overrightarrow{d}(c)=\overleftarrow{d}(c)$
then $c$ must be equal to a multiple of identity.
\end{defi}
\begin{theo}[Hudson \cite{Hudson2}]
Let $C$ be a bialgebra. We embed $C ^{\otimes 2}$ into its {\it{$C$-chiral bimodule}}.
If $\Delta$ is a unital homomorphism, then $\overleftarrow{d},\overrightarrow{d}$ are Ito derivatives.
\end{theo}
\Proof
Let $x,y \in C$, recall that $\overrightarrow{d}(x) = \Delta(x) - \delta(x)$ and $\overleftarrow{d}(x) = \Delta(x) - \tilde{\delta}(x)$.
We have $\overleftarrow{d}(1) =0= \overrightarrow{d}(1)$. Moreover,
$$\overrightarrow{d}(x)\overrightarrow{d}(y) = \Delta(xy) + xy \otimes 1 - \Delta(x)(y \otimes 1) - (x \otimes 1)\Delta(y)$$
We insert: $$ 0 = xy \otimes 1 - xy \otimes 1 $$ and regrouping the terms we find:
$$ \overrightarrow{d}(xy)= \overrightarrow{d}(x)\overrightarrow{d}(y) + \overrightarrow{d}(x)\delta(y) + \delta(x)\overrightarrow{d}(y).$$
Similarly,
$$ \overleftarrow{d}(xy)= \overleftarrow{d}(x)\overleftarrow{d}(y) + \overleftarrow{d}(x)\tilde{\delta}(y) + \tilde{\delta}(x)\overleftarrow{d}(y).$$
\eproof
\end{exam}
\Rk
If $C$ is a coassociative coalgebra with $\Delta(x) = x \otimes 1 + 1 \otimes x$, i.e. $x$ primitive,
then $( \overleftarrow{d} \otimes 1)\overrightarrow{d}(x)=0$ for all $x \in\ C$.
\begin{exam}{[$C$-bicomodule, with $C$ a coassociative coalgebra]}
Let $B$ be a $C$-bicomodule i.e. there exist $\delta$ and $\tilde{\delta}$ such that the following
diagram commute:
\begin{equation*}
        \begin{CD}
        B @>\delta>> B \otimes C  \\
        @V{\tilde{\delta}}VV		@VV{\tilde{\delta} \otimes 1}V \\
        C \otimes B @>{1 \otimes \delta }>> C \otimes B \otimes C
        \end{CD}
\end{equation*}
We define $A = C \otimes B \otimes C$. The coproducts
$\delta,  \tilde{\delta}$ induce coproducts $\Delta,  \tilde{\Delta}$ such that the following diagram commute:
\begin{equation*}
        \begin{CD}
        A @>\Delta>> A \otimes A  \\
        @V{\tilde{\Delta}}VV		@VV{\tilde{\Delta} \otimes 1}V \\
        A \otimes A @>{1 \otimes \Delta }>> A \otimes A \otimes A
        \end{CD}
\end{equation*}
Where :
\begin{equation*}
        \begin{CD}
        u \otimes b \otimes v @>\Delta>> ( u \otimes \delta(b) ) \otimes (v \otimes 1_B \otimes 1_C)  \\
        @V{\tilde{\Delta}}VV		@VV{\tilde{\Delta} \otimes 1}V \\
        ( 1_C \otimes 1_B \otimes u ) \otimes ( \tilde{\delta}(b) \otimes v) @>{1 \otimes \Delta }>> ( 1_C \otimes 1_B
         \otimes u ) \otimes (\tilde{\delta} \otimes id )\delta(b)) \otimes (v \otimes 1_B \otimes 1_C)
        \end{CD}
\end{equation*}
\end{exam}
\begin{exam}{[Tensor product]}
Let $C$, $D$ be $L$-coalgebras. Define the right coproduct $\Delta_{C \otimes D}$ to be the composite:
$$ C \otimes D \xrightarrow{\Delta_C \otimes \Delta_D} (C \otimes C) \otimes  (D \otimes D)
\xrightarrow{1_C \otimes \tau \otimes 1_D} (C \otimes D) \otimes (C \otimes D) $$
and the left coproduct $\tilde{\Delta}_{C \otimes D}$ by:
$$C \otimes D \xrightarrow{\tilde{\Delta}_C \otimes \tilde{\Delta}_D} (C \otimes C) \otimes  (D \otimes D)
\xrightarrow{1_C \otimes \tau \otimes 1_D} (C \otimes D) \otimes (C \otimes D) $$
For the counits part we define:\\
$\epsilon_{C \otimes D}$ as
$ C \otimes D \xrightarrow{\epsilon_C \otimes \epsilon_D} k \otimes k \simeq k$
and $\tilde{\epsilon}_{C \otimes D}$ as
$ C \otimes D \xrightarrow{\tilde{\epsilon}_C \otimes \tilde{\epsilon}_D} k \otimes k \simeq k$\\
With this setting $C \otimes D$ becomes a $L$-coalgebra over $k$.
\end{exam}
\subsection{R\'esum\'e}
For the convenience of the reader we summarize the different categories set forth here.
\begin{center}
\[
\begin{array}{c@{\hskip 1cm}c}
\rnode{a}{L-\textrm{coalgebra}}\\[1cm]
\rnode{b}{\textrm{Non degenerate (Ito)} \ L-\textrm{coalgebra}} & \rnode{c}{\textrm{Coassociative} \ \textrm{coalgebra, (bialgebra)}}\\[1cm]
\rnode{d}{\textrm{Markov} \ L-\textrm{coalgebra}} & \rnode{e}{\textrm{Unital} \ \textrm{algebra}}
\end{array}
\psset{nodesep=3pt}
\ncline{->}{a}{b}
\ncline{->}{a}{c}
\ncline[linestyle=dotted]{<-}{b}{c}
\ncline{->}{b}{d}
\ncline[linestyle=dotted]{<-}{d}{e}
\]
\end{center}
\section{Semantics and completely positive semigroups}
The aim of this section is to construct a completely positive semigroup,
which can be driven from a Markov $L$-coalgebra.
Let $\mathcal{H}$ be a separable Hilbert space. We denote
$B(\mathcal{H})$ the space of bounded operators on $\mathcal{H}$, let $A$ be such an operator,
$A^{*}$ will denote its adjoint. By a theorem due to Stinespring \cite{Stinespring}, a completely positive semigroup on $B(\mathcal{H})$
is a linear mapping
$$ B(\mathcal{H}) \xrightarrow{\Phi} B(\mathcal{H}) \ \ \ \textrm{defined by},$$
$$ \rho \mapsto \Phi(\rho) = \sum_{i=1}^{n} A_i \rho A_i ^{*}, \ \ \textrm{where} \ A_i \in B(\mathcal{H}).$$
\Rk
We choose the Hilbert framework, because we have Quantum
Mechanics in mind.

Suppose we wish to study the iterates of $\Phi$, i.e.
the sequence $(\Phi, \Phi^2, \Phi^3, \ldots, \Phi^n,\ldots)$, where $\Phi^2 \equiv \Phi \circ \Phi$ and so on.
One of the aim of this part is to show that such a process can be recovered from random walks (combinatorics) on
special graphs, and to generalise the operation $\circ$ to any other graphs.
\subsection{The $\circ$ operation}
In the following we do not consider probability on graphs, i.e. we impose 1 over each arrow.
\begin{defi}{}
A graph $G$ is said {\it{complete}} if
any vertex of $G$ is connected by an arrow to any other vertex of $G$, included itself.
\end{defi}
\Rk
In this case: $G=G^{\natural}$.
\Rk
Let $n > 0$.
In the following we consider only $G = (A_1, \ldots,A_n)$ a subset of $B(\mathcal{H})$.
The free semigroup generated by $G$ with the product $m$ of $B(\mathcal{H})$
is denoted by $(\tilde{G},m)$.
\begin{defi}{}
We define
$$R : (\tilde{G},m) \xrightarrow{} \textsf{Hom}(B(\mathcal{H}) )$$
$$ X \mapsto X(\cdot)X^{*}. $$
\end{defi}
\begin{theo}
Let $G = (A_1, \ldots,A_n)$ be in $B(\mathcal{H})$ and consider the complete
graph having $G$ as its vertex set with the natural
Markov $L$-coalgebra associated with it.
The action of the right coproduct $\Delta$ and the product $m$ of $B(\mathcal{H})$
over $(A_1, \ldots,A_n)$ generates $\Phi^2$.
\end{theo}
\Rk
$R$ can be now naturally extended by linearity.
\Proof
It suffices to remark that $\Phi^2(\cdot) = \sum_{i,j=1}^{n} A_jA_i(\cdot)(A_jA_i)^{*}$. Noticing that
$\Delta(A_i) = A_i \otimes \sum_{i=1}^{n} A_j$ we see that $\sum_{i=1}^{n} m\Delta(A_i)$ gives the same
words. Thus $\Phi^2(\cdot) = \sum_{i=1}^{n} R(m\Delta(A_i))(\cdot) = \sum_{i,j=1}^{n} R(A_iA_j)(\cdot)$.
\eproof \\
\begin{coro}
The set of paths on such a graph $G$ is encoded into the sequence $(\Phi, \Phi^2, \Phi^3, \ldots, \Phi^k,\ldots).$
\end{coro}
\Proof
In the previous part we observed that all paths emerging from a given vertex could be obtained by the right
coproduct thanks to the sequence $(\Delta_1 \equiv \Delta, \Delta_2 = 1 \otimes \Delta, \Delta_3 = 1 \otimes 1 \otimes \Delta, \ldots $).
Iterating the previous construction, we observe that $\Phi^k$ can be constructed by $\sum_{i=1}^{n} m\Delta_k(A_i)$, since
$\Phi^k(\cdot) = \sum_{i=1}^{n} R(m\Delta_k(A_i))(\cdot)$.
\eproof \\
The operation $\circ$ for completely positive map not only generates the completely positive
semigroup $(\Phi^k)_{k \in \mathbb{N}}$ (because $\Phi^{k+l}=\Phi^{k} \circ \Phi^{l} $)
but also generates all the paths on the complete graph with $n$ vertices via $\Delta$.
The idea is to replace $\circ$ by $\circ_G$ that is to generalise this procedure to any Markov
$L$-coalgebra.
\begin{defi}{}
Let $G$ be a Markov $L$-coalgebra generated by $n$ vertices $(A_1, \ldots,A_n) \in B(\mathcal{H})$
equipped with this coproduct $\Delta_G$. Let $\Psi(\cdot):
B(\mathcal{H}) \xrightarrow{} B(\mathcal{H}), $
$ \rho \mapsto \Psi(\rho) = \sum_{i=1}^{n} A_i \rho A_i ^{*}.$
We define for all $k > 1, \ \Psi^{\circ_G k}(\cdot)= \sum_{i=1}^{n} R(m\Delta_{(k), G}(A_i))(\cdot)$.
\end{defi}
\begin{theo}
All the linear mappings $(\Psi, \Psi^{\circ_G 2}, \Psi^{\circ_G 3}, \ldots, \Psi^{\circ_G n}, \ldots)$
are completely positive.
\end{theo}
\Proof
Obvious by the Stinespring theorem.
\eproof \\

As for $\circ$, we wish that $\circ_G$ generates a semigroup.
\begin{defi}{}
We define $\Psi^{\circ_G k} \circ_G \Psi^{\circ_G l}(\cdot) = \sum_{i=1}^{n} R(m(\Delta_{(k), G} \star \Delta_{(l), G}(A_i))(\cdot)$,
where the product \footnote{ To simplify exposition, we do not consider probability measures on a graph. We see here that the introduction
of such probability measures do not affect the present setting.} $\star$ is defined in \ref{troisdeux}.
\end{defi}
\Rk
With such a definition the sequence $(\Psi^{\circ_G k})_{k \in \mathbb{N}}$ is a completely positive
semigroup driven by the graph $G$. Moreover if $G$ is complete, we recover the $\circ$ product.
\Rk
All what said depend only on the right coproduct of the Markov $L$-coalgebra, and so can be extended
easily to any $L$-coalgebra. Only the concreteness is lost.
\begin{exam}{}
As an example, consider the graph associated with the first coassociative coalgebra
we met in the second part.
By definition we have:
$$\Psi(\cdot)= 1(\cdot)1 + X(\cdot)X^{*} + g(\cdot)g^{*}.$$
Applying $\Delta$ to $(1,X,g)$ we get,
$$\Psi^{\circ_G 2}(\cdot)= m(1 \otimes 1)(\cdot)m(1 \otimes 1) + m(X \otimes 1 + g \otimes X)(\cdot)m(X \otimes 1 + g \otimes X)^{*} + m(g \otimes g)(\cdot)m(g \otimes g)^{*},$$
that is, $\Psi^{\circ_G 2}(\cdot)= 1(\cdot)1 + (X  + gX)(\cdot)(X  + gX)^{*} + g^2(\cdot)(g^2)^{*}$
and so forth by applying the power of $\Delta$.
\end{exam}
\Rk
The reader can now notice the interest to link sinks and sources to the identity element.
\begin{exam}{[contractive completely positive $n$-tuple]}
In \cite{Popescu} Popescu is interested in a contractive $n$-tuple, that is
$n$ operators $(T_1, \ldots, T_n)$ from $B(\mathcal{H})$ such that $\sum_{i=1}^n T_iT_i^* \leq I$.
Such a $n$-tuple generates, for the $\circ$ operation, a contractive completely positive semigroup.
Suppose now we consider a quantum channel, $\Psi$, that is $n$-tuple whose $\sum_{i=1}^n T_iT_i^* = I$
and we study the semigroup of the power of this channel generated by the $\circ_G$ operation,
---with $G$, a Markov $L$-coalgebra with $n$ vertices--- how such a semigroup behaves ?

If $G$ is not a complete graph, $\Psi^{\circ_G k}$ for $k > 1$ is not unital any longer, in fact
$\Psi^{\circ_G k}(I) < I$.
Moreover if $G$ has an attractor $A$, that is if all the paths of the graph $G$ will converge
to $A$ the limit, when $k \xrightarrow{} \infty$, of $\Psi^{\circ_G k}$ will be only composed by
operators defining the attractor $A$.
\end{exam}
\section{The flower graph}
\subsection{Introduction and motivation}
In \cite{Quillen} Quillen defines the notion of curvature of a linear map relative
to the Hochschild homology of an associative algebra $A$ with product $m$.
Here we study only the case of a unital algebra. The unit element will be noted 1 or $I$.
For the convenience of the reader we remind briefly the boundary operators $b$ and $b'$. \\
\noindent
\textbf{Notation}: For every $n \in \mathbb{N}$ we denote $(a_1,a_2,...,a_n)$ the tensor product: $a_1 \otimes a_2 \otimes...\otimes a_n$
and for $n > 1$, we define
$b': A^{\otimes n} \xrightarrow{} A^{\otimes (n-1)},$
$$
b'(a_1, ... , a_n) = \sum_{i=1}^{n-1} (-1)^{i-1}(a_1,...,a_ia_{i+1},...,a_n).
$$
Since $A$ is unital, the $b'$-complex,
$$ ... \xrightarrow{b'} A^{\otimes 3} \xrightarrow{b'} A^{\otimes 2} \xrightarrow{b'} A \xrightarrow{} 0, $$
is exact.
The {\it{differential}} $b: A^{\otimes n} \xrightarrow{} A^{\otimes (n-1)}$ is usually defined by:
$$b(a_1, ... , a_n) = b'(a_1, ... , a_n) + (-1)^{n-1}(a_na_1, ... , a_{n-1}).$$
The Hochschild homology $H(A,A)$, when $A$ is unital, is usually computed from the $b$-complex:
$$ ... \xrightarrow{b} A^{\otimes 3} \xrightarrow{b} A^{\otimes 2} \xrightarrow{b} A \xrightarrow{} 0. $$
\begin{defi}{}
Let $A$ a unital algebra, we denote $F(A) = k \oplus A  \oplus A^{\otimes 2} \oplus A^{\otimes 3} \ldots$
the {\it{Fock space}} of $A$.
\end{defi}
The coproduct usually used on $F(A)$ is
$ \Delta : F(A) \xrightarrow{} F(A) \otimes F(A), $
$$\Delta(a_1, ... , a_n) = \sum_{i=0}^{n} (a_1,...,a_i) \otimes (a_{i+1}, ... , a_n).$$
The operator $\Delta$ acts as a deconcatenation map if we think of $(a_1, ... , a_n)$ as path over a graph.
By a $n$-cochain on $A$ we mean a multilinear function $f(a_1, ... , a_n)$ with values in some
vector space $V$ or equivalently a linear map from $A^{\otimes n}$ to $V$.
These cochains form a complex $\textsf{Hom}(F(A), V)$, where the differential is
$\hat{\delta}(f) = -(-1)^{n}fb', \ \  f \in \textsf{Hom}(A^{\otimes n}, V).$
Suppose $L$ is a unital algebra, the complex of cochains $\textsf{Hom}(F(A), L)$ has a product defined by,
$ fg=m(f \otimes g)\Delta$
where $ m: L \otimes L \xrightarrow{} L $ is the product of $L$.
If $f$ and $g$ have respectively degrees $p$ and $q$, we define the associative product,
$$ (fg)(a_1, ... , a_n)= (-1)^{pq} f(a_1, ... , a_p)g(a_{p+1}, ... , a_{p+q}).$$
As an important example, we have:
\begin{defi}{[Curvature of a 1-cochain \cite{Quillen}]} Let $\rho$ be a 1-cochain, that is a linear map from $A$ to $L$.
We can view $\rho$
as a connection form and construct its {\it{curvature}}: $\omega = \hat{\delta} \rho + \rho^2$, which will be a
2-cochain. Then, the {\it{curvature}},
$$\omega(a_1,a_2) = (\delta \rho + \rho^2)(a_1,a_2)=\rho(a_1a_2)-\rho(a_1)\rho(a_2),$$
quantifies how close $\rho$ is to a homomorphism.
\end{defi}
\begin{theo}
A unital algebra $A$ is a Markov $L$-bialgebra.
\end{theo}
\Proof
We define below all the operations showing that $A$ is in fact a Markov $L$-bialgebra.
The right coproduct $\delta: A \xrightarrow{} A^{\otimes 2} $ is defined by $a \mapsto a \otimes 1$ and
the left coproduct $\tilde{\delta}: A \xrightarrow{} A^{\otimes 2} $ by $a \xrightarrow{} 1 \otimes a$,
showing thus the Markov property.
For the right and left counit $\epsilon,\tilde{\epsilon}: A \xrightarrow{} k $ choose
any unital algebra maps.
We have for all $a \in A$,
$$(\tilde{\delta} \otimes 1)\delta a = (1 \otimes \delta)\tilde{\delta} a = 1 \otimes a \otimes 1 \ \
\textrm{and} \ \
(\tilde{\epsilon} \otimes 1)\tilde{\delta} a = a =(1 \otimes \epsilon)\delta a. $$
The $L$-bialgebra property is given by,
$$\delta(ab) = ab \otimes 1 = (a \otimes 1)(b \otimes 1) = \delta(a)\delta(b), $$
$$\tilde{\delta}(ab) = 1 \otimes ab = (1 \otimes a)(1 \otimes b) = \tilde{\delta}(a)\tilde{\delta}(b). $$
Moreover, $\delta(1) = 1 \otimes 1 = \tilde{\delta}(1)$.
With the choice we made for the counits we have automatically the algebra map property.
\Rk
We can embed $A$ into a $L$-Hopf algebra by choosing any unital map. It is of no use for the following.
\begin{defi}{[Flower graph]}
We call such a Markov $L$-coalgebra a \textit{flower graph} because it is the concatenation of petals:
\begin{center}
\includegraphics*[width=5cm]{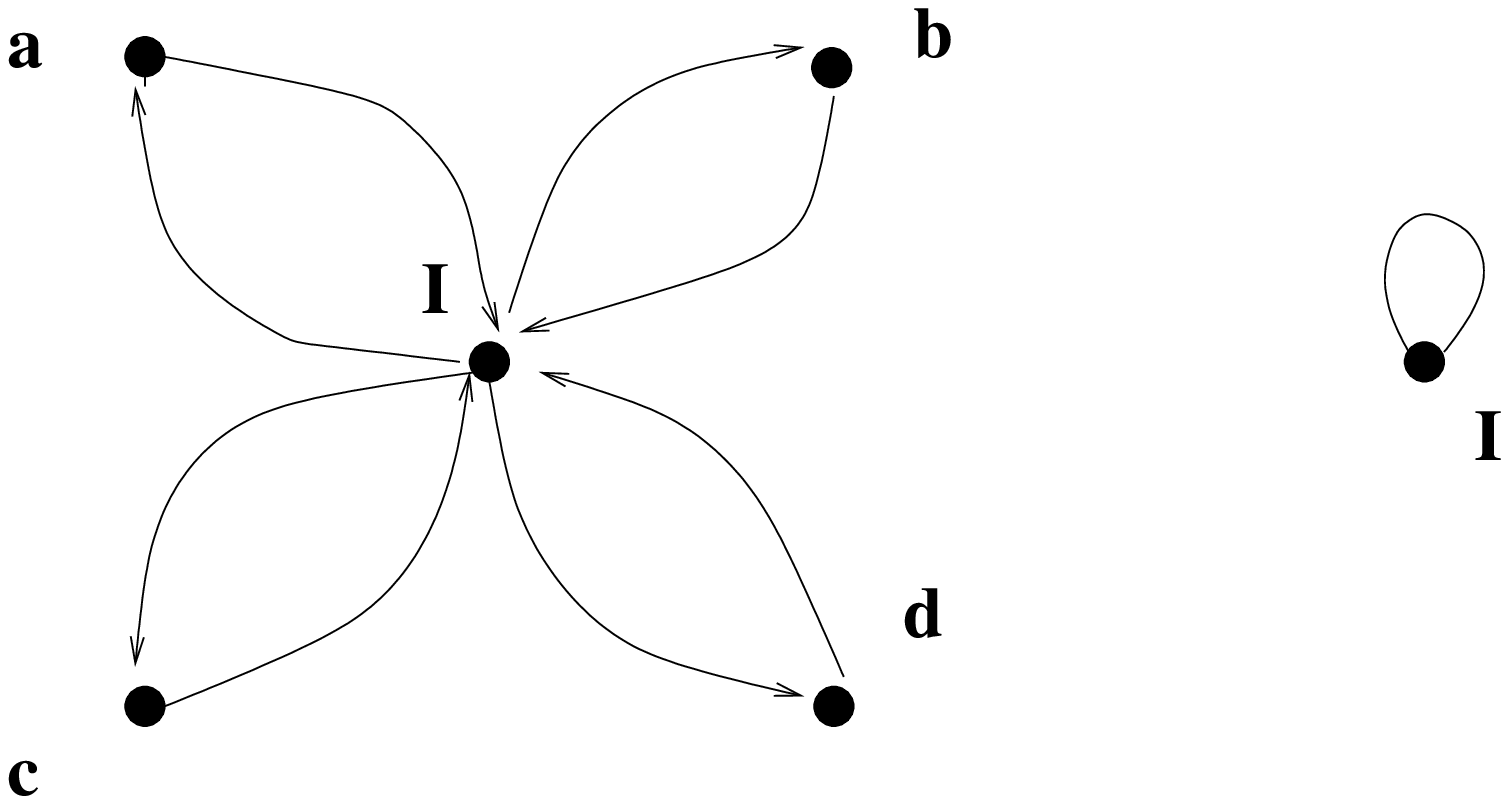}
\end{center}
\end{defi}
We think that the flower graph \footnote{As the coproducts define completely a graph,
we omit the denumerability condition of the definition of a graph.} is a fundamental object for a unital algebra, i.e.
with the product $m$ we can rediscover the fundamental concepts such as differential, Ito derivative,
commutator, $b$ and $b'$ complex and so on.
We start from an unital associative algebra $A$ equipped with its flower graph.
\begin{defi}{[pattern]}
An element $(a_1, \ldots , a_n) \in A^{\otimes n}$ defines a periodic orbit on the flower graph, we denote the periodic orbit,
$(\ldots ,I,a_1,I,a_2,I,a_3,I, \ldots ,I,a_n,I,a_1,I, \ldots,I,a_n,\ldots )$, where $I$, for convenience,
notes the neutral element 1.
A \textit{pattern} is denoted by the ordered set $[I,a_1,I,a_2,I,a_3,I,\ldots,I,a_n]$, repeated infinitly often
it generates the periodic orbit.
We call left border of the pattern the symbol $a_n$ and the right border the symbol $I$.
It can be represented by $\ldots a_n,[I,a_1,I,a_2,I,a_3,I, \ldots,I,a_n],I \ldots$
\end{defi}
To get information contained into a periodic orbit it suffices to read either its pattern or
if we want more information, its pattern and its border.
\begin{defi}{[Reading map]}
$\mathcal{L}$ is a reading map for the unital algebra $A$ if it obeys the following conditions:
\begin{enumerate}
\item {$\mathcal{L}$ starts always after the symbol $I$, ($I$ plays the r\^ole of a stop codon).}
\item {$ \mathcal{L} : A^{\otimes 3} \xrightarrow{} A $ (the reading is done three by three).
	The simplest function we can consider must use the sole information we have
        about the algebra that is $m$. We can choose: $\mathcal{L} = m(m \otimes id) = m(id \otimes m)$.}
\item {Each shift of the reading map provokes the apparition of a minus sign.}
\item {$\mathcal{L}$ must explore either all the pattern or the pattern and its border if we want more information.}
\end{enumerate}
\end{defi}
\Rk
The reading must be done three by three because the flower graph is the concatenation of petals,
and to cover such a petal one must start from let say $I$, then to $a_i$ for some $i$ then to $I$. It demands
three letters to write down.
\Rk
The petal of the flower graph arises naturally when one considers a unital algebra $A$ because every
element $a$ from $A$ can be written as: $a = 1 \cdot (a \cdot 1) = (1 \cdot a) \cdot 1$.
\begin{prop}
Reading the pattern of a periodic orbit produces the complex with boundary $b'$ and
reading the pattern and its border produces the complex with boundary $b$.
\end{prop}
\Proof
We proceed by recurrence.
As $b$ and $b'$,
the reading map is not defined on periodic orbits of period 1 because there is not enough information to
read 3 by 3.\\
Let $(a_1,a_2) \in A^{\otimes 2}$ be a periodic orbit of period 2.
We have the sequence : $\ldots a_2,[I,a_1,I,a_2],I \ldots$
For the moment we focus only on the pattern.
We start the reading after the stop codon $I$
and we find $\ldots a_2,[I, \mathcal{L}(a_1,I,a_2)],I \ldots$ that is $a_1a_2$.
Remark that $b'(a_1,a_2)=a_1a_2$.
Now if we want more information we can read the pattern and the border. Yet a problem appears.
We would have to write $\ldots a_2,[I,a_1,I, \mathcal{L}(a_2],I \ldots$ with a minus sign, but then we cannot read three by three any more.
The only thing we can do to read is to use the left boarder, by doing so we shift the pattern too behind.
So we have: $\ldots [ \mathcal{L}(a_2,I,a_1)],I,a_2,I \ldots$ and we find $-a_2a_1$.
The complete reading gives: $a_1a_2-a_2a_1$. This is equal to $b(a_1, a_2)$.\\
Let $(a_1,a_2,a_3) \in A^{\otimes 3}$ be a periodic orbit of period 3,
that is we have the sequence: $$\ldots a_3,[I,a_1,I,a_2,I,a_3],I \ldots$$
We focus on the reading of the pattern.\\
\textit{First step}:
$\ldots a_3,[I, \mathcal{L}(a_1,I,a_2),I,a_3],I \ldots$ gives $a_1a_2 \otimes a_3 $.\\
\textit{Second step}:
$\ldots a_3,[I,a_1,I, \mathcal{L}(a_2,I,a_3)],I \ldots$ gives $a_1 \otimes a_2a_3 $  with a minus sign.
We have finished to read the pattern. We obtain $a_1a_2 \otimes a_3 $  $- a_1 \otimes a_2a_3 $.
This is equal to $b'(a_1,a_2,a_3)$.
If we want more information we can read the boarder too.\\
\textit{Third step}:
$\ldots[ \mathcal{L}(a_3,I,a_1),I,a_2],I,a_3,I \ldots$ gives $a_3a_1 \otimes a_2 $ with a minus sign. But
at the step before we had a minus sign too, so we get a plus sign.
The complete reading gives: $a_1a_2 \otimes a_3 - a_1 \otimes a_2a_3 + a_3a_1 \otimes a_2$.
This is equal to $b(a_1,a_2,a_3)$.
By repeating this process by an obvious recurrence, we recover $b'$ and $b$.
We can interpret the $b'$ and $b$ complex of a unital algebra only from its Markov $L$-coalgebra,
that is by reading periodic orbits from the flower graph.
\eproof

From the coproducts of a unital algebra, we remark that we can recover
the well-known cocommutative coassociative coproduct $\Delta$
since
for $a$ different from 1 we get:
$\Delta a = a \otimes 1 + 1 \otimes a = \delta(a) + \tilde{\delta}(a)$
and
$\Delta 1 = 1 \otimes 1 = \delta(1) = \tilde{\delta}(1).$
As well, if
$\Omega$ denotes the (non commutative) differentials over $A$, the differential,
$ d : A \xrightarrow{} \Omega, \ \ a \mapsto da = a \otimes 1 - 1 \otimes a$
can be decomposed into $d=\delta-\tilde{\delta}$.
If
we define $ [\cdot, \cdot]_{\otimes}: \  A^{\otimes 2} \xrightarrow{} A^{\otimes 2}$ by
$a \otimes b \mapsto [a,b]_{\otimes} = a \otimes b - b \otimes a$,
we observe that $[a , b]_{\otimes} = \delta(a)\tilde{\delta}(b) - \tilde{\delta}(a)\delta(b)$.
However we shall define $[\cdot, \cdot]_{\otimes} = \delta(\cdot)\tilde{\delta}(\cdot) - \delta(\cdot)\tilde{\delta}(\cdot)\tau(\cdot)$
to keep the null property of the commutator on an element invariant by $\tau$.
\subsection{Curvature of a 1-cochain}
\label{curv}
Let $A, B$ two unital algebras.
Consider $\rho$ as a 1-cochain from $A$ to $B$,
with curvature $\omega$. Let $G_A$ be a graph whose vertex set is a subset of the algebra $A$.
We embed $G_A$ with its natural
Markov $L$-coalgebra. How such a  Markov $L$-coalgebra can be conveyed in $B$ by the 1-cochain $\rho$ \cite{Lersemi}?
To answer this question we denote by $\textsf{Prop}_M(A)$, (resp. $\textsf{Prop}_M(B)$) the set of Markovian coproducts which can equipped any graph $G_A$, (resp. $G_B$)
whose vertex set is a subset from the algebra $A$, (resp. $B$).
\begin{theo}
The curvature $\omega: A^{\otimes 2} \xrightarrow{} B$ induces a mapping $\omega_{*}:  \textsf{Prop}_M(A)
\xrightarrow{} \textsf{Prop}_M(B)$.
\end{theo}
\Proof
Let $x,y \in A$,
such that
$\Delta x = x \otimes x_1, \
\Delta y = y \otimes y_1, \
\tilde{\Delta}x = x_0 \otimes x, \
\tilde{\Delta}y = y_0 \otimes y,$ in symbolic notation.
Define the right coproduct: $\Delta_{\omega} : B \xrightarrow{} B^{\otimes 2} $ by
$\Delta_{\omega}\omega(x,y) = \omega(x,y) \otimes \omega(x_1,y_1)$,
and the left coproduct:
$\tilde{\Delta}_{\omega} : B \xrightarrow{} B^{\otimes 2} $ by
$\tilde{\Delta}_{\omega}\omega(x,y) = \omega(x_0,y_0) \otimes \omega(x,y)$
They obey the coassociativity breaking equation. Thus $( \Delta_{\omega},\tilde{\Delta}_{\omega} ) \in \textsf{Prop}_M(B)$.
\eproof
\begin{coro}
Let $A, B$ be as above.
The curvature of the 1-cochain $\rho$ sends the Markov $L$-coalgebra $G_A$ into a Markov $L$-coalgebra.
\end{coro}
\Proof
The theorem above gives us the coproducts.
For the new right and left counits we take $\epsilon_{\omega}(\omega(x,y)) = 1$ and $\tilde{\epsilon}_{\omega}(\omega(x,y)) = 1$.
By checking the equations (see subsection \ref{rdd}) for
the counits we are faced with $ \ldots \otimes \sum_{i,j}P_x(i)P_y(j) \omega(x_i,y_j)$
and we know that $\sum_{i,j}P_x(i)P_y(j)=1$.
\eproof
\Rk
We follow Quillen \cite{Quillen}. Let $\sigma: L \xrightarrow{} V$ be a trace on the algebra $L$ with
values in the vector space $V$, i.e. a linear map vanishing on the commutator subspace $[L,L]$.
The image of the map $m-m\tau: L \otimes L \xrightarrow{} L$ giving $[L,L]$, we can say
that a trace is realy a linear map defined on the commutator quotient space:
$$ L_{\natural}:= L / [L,L] = \textrm{coker} \{ m-m\tau: L \otimes L \xrightarrow{} L \}.$$
Thus naturally associated with the bar construction of $A$ is its cocommutator subspace,
$$ F(A)^{\natural}:= \ker \{ \Delta - \tau \Delta: F(A) \xrightarrow{} F(A) \otimes F(A) \}. $$
If  $\natural: F(A)^{\natural} \xrightarrow{} F(A)$ denote the inclusion map, it is the universal cotrace in the sense that
a cotrace $L \xrightarrow{} F(A)$
is the same as a linear map with values in $F(A)^{\natural}$.
By combining the trace $\sigma$ with the universal cotrace $\natural$, Quillen defines a morphism of complexes
$$ \sigma^\natural: \textrm{Hom}(F(A),L) \xrightarrow{} \textrm{Hom}(F(A)^\natural,V), \ \ \ \ \sigma^\natural(f)=\sigma f \natural,$$
which is a trace on the DG algebra of cochains i.e. it vanishes on $[f,g]= fg - (-1)^{deg(f)deg(g)}gf$.
We recall that
$$ \sigma^\natural(\omega(a_1,a_2), \ldots, \omega(a_{2n-1},a_{2n}) = n\sigma(\omega(a_1,a_2), \ldots, \omega(a_{2n-1},a_{2n}) -
\omega(a_{2n},a_1), \ldots, \omega(a_{2n-2},a_{2n-1})),$$
is a cyclic cocycle of degree $2n-1$.
Let us denote by $m_B$ the product in $B$ and $\sigma: B \xrightarrow{} V$.
The image by $\sigma^{\natural} m_B$ of the coproducts
$\Delta_{\omega}, \ \tilde{\Delta}_{\omega}$
induced by the curvature of $\rho: A \xrightarrow{} B$ will give
cyclic cocycles of degree 3.
\begin{prop}
Let $C$ the sub-algebra of $A$ generated by $\langle 1,X_1,...,X_n \rangle$ and $\rho$ be a 1-cochain from $A$ to $A$.
Denote by $\omega$ its curvature and embed $C$ into a Markov $L$-coalgebra by fixing a graph on it.
We demand that $C$ does not belong to the image of $\omega$, (compatibility property).
Under this assumption, we can identify the induced coproducts by $\omega$ with the old ones and
the curvature $\omega$ embeds the Markov $L$-coalgebra $C$ into a $L$-bialgebra with
$\omega$ as a (non associative) product.
\end{prop}
\Proof
We start by noticing that if $m$ is a product in a $L$-bialgebra $C$ we have, if
$a,a_1,b,b_1 \in C$, $\Delta a = a \otimes a_1$ and $\Delta b = b \otimes b_1$ (in symbolic notation).
The homomorphism property of the coproduct reads:
$\Delta(m(a \otimes b )) = m(\Delta(a)\Delta(b)) $
i.e.
$$\Delta(m(a,b)) = ab \otimes a_1b_1 = m(a,b) \otimes m(a_1,b_1), \ \ \tilde{\Delta} m(a,b) = m(a_0,b_0) \otimes m(a,b).$$
However we have constructed the induced coproducts,
(identified here to the old coproducts, thanks to the compatibility property) such that:
\begin{eqnarray*}
\Delta(\omega)(a,b) &=& \Delta\omega(a,b) = \omega(a,b) \otimes \omega(a_1,b_1) =  (\omega)(a,b) \otimes (\omega)(a_1,b_1),\\
\tilde{\Delta}(\omega)(a,b) &=& \tilde{\Delta}\omega(a,b) = \omega(a_0,b_0) \otimes \omega(a,b) = (\omega)(a_0,b_0) \otimes (\omega)(a,b).
\end{eqnarray*}
Thus we reproduce the action of what we know about the ordinary product when we say
that such or such coproduct is a $m$-homomorphism, i.e.
the following diagram,
\[
\begin{array}{c@{\hskip 1cm}c@{\hskip 1cm}c}
\rnode{a}{C \otimes C} & \rnode{b}{C} & \rnode{c}{C \otimes C}\\[1cm]
\rnode{d}{C \otimes C \otimes C \otimes C} & \ & \rnode{e}{C \otimes C \otimes C \otimes C}
\end{array}
\psset{nodesep=3pt}
\everypsbox{\small}
\ncline{->}{a}{b}\Mput[b]{\omega}
\ncline{->}{b}{c}\Mput[b]{\Delta}
\ncline{->}{a}{d}\Mput[r]{\Delta \otimes \Delta}
\ncline{->}{d}{e}\Mput[b]{id \otimes \tau \otimes id}
\ncline{->}{e}{c}\Mput[l]{\omega \otimes \omega}
\]
commutes.
\eproof
\subsection{Ito derivatives}
A common point between the product map $m$ of an unital associative algebra, Leibnitz derivatives,
homomorphisms and Ito derivatives from $A$ to $A$, is that they all obey a same equation.
To explain what this equation is we need some definitions. We will follow part 3 of the
Quillen's paper \cite{Quillen}.

Let $A$ be an algebra with product map $ m: A^{\otimes 2} \xrightarrow{} A $ and unit map $ \eta: k \xrightarrow{} A $.
Let $M$ be a $A$-bimodule, that is, a vector space with left and right product maps
$ m_l: A \otimes M \xrightarrow{} M, m_r: M \otimes A \xrightarrow{} M, $
defining left and right module structures which commute that is, $m_r(m_l \otimes id) = m_l( id \otimes m_r ).$
Consider $A \otimes V \otimes A$, where $V$ is a vector space, as a $A$-bimodule with $ m_r= id \otimes id \otimes m$
and $ m_l= m \otimes id \otimes id$.
We have the following proposition:
\begin{prop}[Quillen \cite{Quillen}]
There is a one-one correspondence between linear maps $h: V \xrightarrow{} M $ and bimodule morphism \footnote{The following tilde notation has nothing to do with the tilde notation of the left coproduct of a $L$-coalgebra.}
$ \tilde{h}: A \otimes V \otimes A \xrightarrow{} M$ given by:
$$ \tilde{h} = m_r(m_l \otimes id)(id \otimes h \otimes id), h=\tilde{h}(\eta \otimes id \otimes \eta).$$
\end{prop}
Consider the exact $b'$-complex (defined at the begining of this section):
$$ ... \xrightarrow{} A^{\otimes 4} \xrightarrow{\Carre} A^{\otimes 3} \xrightarrow{m \otimes id - id \otimes m} A^{\otimes 2} \xrightarrow{m} A
\xrightarrow{} 0,$$
where $\Carre = m \otimes id \otimes id - id \otimes m \otimes id + id \otimes id \otimes m.$
The bimodule $ \Omega_A $ of the (noncommutative) differentials over $A$ is defined to be the kernel of $m$.
We say that the linear mapping $ D : A \xrightarrow{} M $, where $M$ is a bimodule,
is a {\it{Leibnitz derivative}} if the corresponding bimodule morphism $\tilde{D} : A^{\otimes 3} \xrightarrow{} M$
satisfies:
$$\tilde{D}\Carre = 0, $$
that is:
$ 0 = \tilde{D}\Carre (1_A,x,y,1_A) = \tilde{D}( x,y,1_A) - \tilde{D}(1_A,xy,1_A) + \tilde{D}(1_A,x,y),$
or $ 0 = xD(y) - D(xy) +(Dx)y .$
At the begining of this part we said that we have to recover all the basic tools used on an algebra with
the help of the product $m$ and the flower graph generated by $\delta, \tilde{\delta}$.
We observe that we can naturally find the Leibnitz derivative as stated above
but not the Ito derivative.
For the convenience of the reader we remind the following definition:
\begin{defi}{[Ito derivation]}
A linear map $d_{\Ito} : A \xrightarrow{} B$ such that:
\begin{enumerate}
\item {$d_{\Ito}(1_A) = 0_B,$}
\item {$d_{\Ito}(xy) =  d_{\Ito}(x)y + xd_{\Ito}(y) + d_{\Ito}(x)d_{\Ito}(y),$}
\end{enumerate}
is called an {\it{Ito derivative}}.
\end{defi}
We note $\omega_{\Ito}$ the curvature of such a derivative. Without loss of generality we shall only consider
$A = B$.
\begin{lemm}
Let $(x,y) \in A$.
$\omega_{\Ito}(1,x) = \omega_{\Ito}(x,1) = d_{\Ito}(x) $.
\end{lemm}
\Proof
By definition $\omega_{\Ito}(1,x) = d_{\Ito}(1.x) - d_{\Ito}(1)d_{\Ito}(x) $
and $d_{\Ito}(1) = 0$.
\eproof

We now adapt the Quillen's proposition to $V = A^{\otimes 2}$.
\begin{defi}{}
Let $\Xi$ denote the pairing isomorphism,
$ \Xi : A^{\otimes 4} \xrightarrow{} A \otimes (A^{\otimes 2}) \otimes A$,
$$ (a_1,a_2,a_3,a_4) \mapsto a_1 \otimes (a_2,a_3) \otimes a_4.$$
\end{defi}
When we started from $A$ and began to consider the repercussion of the flower graph on $A$, we found
that we could recover the differential $d$ as $d = \delta - \tilde{\delta}$.
Now if we start from $A^{\otimes 2}$ the possible actions of the coproducts of the flower graph on it can be for example
the mappings $\delta \otimes \delta$, $\tilde{\delta} \otimes \delta$, $\delta \otimes \tilde{\delta}$, $\tilde{\delta} \otimes \tilde{\delta}$.
\begin{defi}{}
The following sequence, (not exact in $A^{\otimes 4}$), allows us to define an important mapping $\Carre_* = \Xi \Carre' $,
$$ 0 \xrightarrow{} A^{\otimes 2} \xrightarrow{\Carre'} A^{\otimes 4} \xrightarrow{\Xi} A \otimes (A^{\otimes 2}) \otimes A  \xrightarrow{} 0, $$
with $\Carre' = \delta \otimes \delta - \delta \otimes \tilde{\delta} - \tilde{\delta} \otimes \delta + \tilde{\delta} \otimes \tilde{\delta}$.
\end{defi}
\Rk
Notice that $\Carre' = d \otimes d$ and $\Carre'd=0$ \footnote{These equations are not a particular
case of the flower graph. If instead,
we would have a graph with only orbits of period two, say $(\ldots a,a_0,a,a_0 \ldots)$
and if we would define $d(a) := a \otimes a_0 - a_0 \otimes a$, we would still have $\Carre'd=0$. }
. Hence,
$$ 0 \xrightarrow{} A \xrightarrow{d} A^{\otimes 2} \xrightarrow{d \otimes d} A^{\otimes 4} \xrightarrow{id \otimes id \otimes d \otimes d \pm d \otimes d \otimes id \otimes id} A^{\otimes 6} \ldots $$
is a complex built only from the action of the coproducts of the flower graph.

Recall that thanks to Quillen's proposition we denote $ \tilde{h}: A \otimes (A^{\otimes 2}) \otimes A \xrightarrow{} M$
in one to one correspondence with $h: A^{\otimes 2} \xrightarrow{} M $.
Define $\tilde{m}, \tilde{\omega}_{\gamma}, \tilde{\omega}_{\Ito}, \widetilde{Dm}$ the correspondence between
the product $m$, the curvature of an homomorphism $\gamma$, the curvature of an Ito derivative $\omega_{\Ito}$
and the product $m$ composed by a Leibnitz derivative $D$.
\begin{theo}
$(\tilde{m}, \tilde{\omega}_{\gamma}, \tilde{\omega}_{\Ito}, \widetilde{Dm})$ verify the same equation:
$$\tilde{F}\Carre_* = 0.$$
\end{theo}
\Proof
Let $ \tilde{F}: A \otimes (A^{\otimes 2}) \otimes A \xrightarrow{} M$
verifying $\tilde{F}\Carre_* = 0,$ with $\Carre_* = \Xi(\delta \otimes \delta - \delta \otimes \tilde{\delta} - \tilde{\delta} \otimes \delta + \tilde{\delta} \otimes \tilde{\delta})$
that is:
$$ 0 = \tilde{F}(x,(1,y),1) - \tilde{F}(x,(1,1),y) - \tilde{F}(1,(x,y),1) + \tilde{F}(1,(x,1),y).$$
Therefore, $ 0 = x \cdot F(1,y) \cdot 1 - x \cdot F(1,1)y - 1 \cdot F(x,y) \cdot 1 + 1 \cdot F(x,1) \cdot y$.
\begin{itemize}
\item{$F = m$ yields: $0 = xy - xy -xy +xy$.}
\item{$F = \omega_{\gamma}$ yields: $ 0 = 0 -0 -\omega_{\gamma}(x,y) +0$.}
\item{$F = \omega_{\Ito}$ yields: $ 0 = xd_{\Ito}(y) - 0 -(d_{\Ito}(xy) - d_{\Ito}(x)d_{\Ito}(y)) + d_{\Ito}(x)y$.}
\item{$F = Dm$ yields: $ 0 = xD(y) - 0 -D(xy) + D(x)y $.}
\end{itemize}
This complete the proof.
\eproof

The converse is also true.
\begin{theo}
Suppose the curvature $\omega$ of a linear map $\rho$ verifies $\tilde{\omega}\Carre_* = 0 $ .
\begin{itemize}
\item{If $\rho$ send the neutral element of the product to itself, then $\rho$ is a homomorphism.}
\item{If $\rho$ send the neutral element of the product to the neutral element of the
vector space structure, then $\rho$ is an Ito derivative. }
\end{itemize}
\end{theo}
\Proof
 Obvious.
\eproof

How can we produce Ito derivatives from known objects ?
\begin{lemm}
Let $\rho$ an unital homomorphism from $A$ to $A$.
The linear map $d = \rho - id$ is an Ito derivative.
\end{lemm}
\Proof
Let $x,y \in A$. We have:
$$ dxdy = (\rho (x) - x )(\rho (y) - y ) = \rho (xy) -(xy -xy) -x\rho (y) -\rho (x)y + xy,$$
that is: $ dxdy = d(xy) - xdy -d(x)y $.
\eproof
\begin{lemm}
Let $d$ an Ito derivative from $A$ to $A$.
The linear map $\rho = d + id$ is a homomorphism.
\end{lemm}
\Proof
$d(x)d(y) = d(xy) -xd(y) -d(x)y $
and
$ \rho (x) \rho (y) = d(x)d(y) + xd(y) + d(x)y + xy = d(xy) + xy = \rho (xy).$
\eproof

We have just proved the following:
\begin{theo}
Let $A$ be an unital algebra.
The set of Ito derivatives from $A$ to $A$ is in one to one with the set of homomorphisms from $A$ to $A$ \footnote{At the end of this part
we will give an algebraic interpretation of this theorem.}.
\end{theo}
\begin{theo}
Let $\rho$ an unital linear map from $A$ to $A$, with curvature $\omega_{\rho}$.
Decompose $\rho = \zeta + id $ where $\zeta$ is a linear map sending $1_A$ to $0_A$ with curvature
$\omega_{\zeta}$. We have: $\tilde{\omega}_{\rho} \Carre_* = \tilde{\omega}_{\zeta}\Carre_* $.
\end{theo}
\Proof
Let $x,y \in A$.
\begin{eqnarray*}
\tilde{\omega}_{\rho} \Carre_*(x,y) = -\omega_{\rho}(x,y) &=& -\rho(xy) + \rho(x)\rho(y) \\
&=& -((\zeta(xy) + xy) - (\zeta(x) + x)(\zeta(y) + y)) \\
&=& -(\zeta(xy) - \zeta(x)y - x\zeta(y) - \zeta(x)\zeta(y)) \\
&=& x \cdot \omega_\zeta  (1,y) \cdot 1 - x \cdot \omega_\zeta (1,1) \cdot y - 1 \cdot \omega_\zeta (x,y) \cdot 1 + 1 \cdot \omega_\zeta (x,1) \cdot y \\
&=&\tilde{\omega}_{\zeta}\Carre_*(x,y),
\end{eqnarray*}
which completes the proof.
\eproof
\Rk
Let $\textsf{Hom}_{1_A}(A)$ be the set of unital linear maps from $A$ to $A$ and $\textsf{Hom}_{0_A}(A)$ be the set of linear
maps which send $1_A$ to $0_A$.
The one to one mapping $\Psi$:
$$ \textsf{Hom}_{0_A}(A) \xrightarrow{\Psi} \textsf{Hom}_{1_A}(A) $$
$$ \zeta \mapsto \rho, $$
can be viewed as a gauge transformation in $\textsf{Hom}(A)$ which leaves the curvature of the maps involved,
invariant by $\Carre_*$.
\subsubsection{Free module and Hochschild complex}
\label{fn}
Let $A$ be a unital algebra.
By establishing a link between $A^{\otimes 2}$ and its free bimodule
$ A \otimes (A^{\otimes 2}) \otimes A$ we showed a common point between
the curvature of Ito derivatives and Leibnitz derivatives. We go further by showing the usefulness of
the Bianchi identity applied to the curvature of an Ito derivative.
We recall that  if $\rho$ is a 1-cochain, we define \cite{Quillen}
its curvature by $\omega = \hat{\delta} \rho + \rho^2$
where $ \hat{\delta}$ is related to the Hochschild boundary $b'$. The Bianchi identity
reads $\hat{\delta}\omega = -[\rho, \omega] $ with
$$[\rho, \omega](a_1,a_2,a_3) = \rho(a_1) \omega(a_2,a_3) - \omega(a_1,a_2)\rho(a_3)
= \omega(a_1a_2,a_3) - \omega(a_1,a_2a_3) \ \ \forall a_i \in A.$$
If $\rho$ is an Ito derivative, we have $ \omega(a_0,a_1) = a_0\rho(a_1)
+\rho(a_0)a_1$.
We denote $ A \otimes (A ^{\otimes n}) \otimes A = \hat{A} ^{\otimes n}$.
To explain the following theorem, recall we showed $\tilde{ \omega} \Carre_{*} = 0$. Is it possible to
construct multilinear maps $f_n: A^{\otimes n} \xrightarrow{} A$, operators $\Carre_{*n}: A^{\otimes n} \xrightarrow{} \hat{A} ^{\otimes n}$
such that $\tilde{f_n} \Carre_{*n} = 0$ ?
In the following we denote $\tilde{A} ^{\otimes n} := \textsf{Im}(\Carre_{*n})$.
\begin{theo}
Let $\rho$ be an Ito derivative, with curvature $\omega$.
We have the following complex between the Hochschild complex with boundary $b'$ and its
associated free bi-module:
\begin{center}
\[
\begin{array}{c@{\hskip 1cm}c@{\hskip 1cm}c@{\hskip 1cm}c@{\hskip 1cm}c@{\hskip 1cm}c}
\rnode{a1}{ } & \rnode{a2}{0} & \rnode{a3}{0} & \rnode{a4}{0} & \rnode{a5}{0} & \rnode{a6}{ } \\[1cm]
\rnode{b1}{0} & \rnode{b2}{A} & \rnode{b3}{A{^{\otimes 2}}} & \rnode{b4}{A{^{\otimes 3}}} & \rnode{b5}{A{^{\otimes 4}}} & \rnode{b6}{\ldots} \\[1cm]
\rnode{c1}{0} & \rnode{c2}{\tilde{A}} & \rnode{c3}{\tilde{A}{^{\otimes 2}}} & \rnode{c4}{\tilde{A}{^{\otimes 3}}} & \rnode{c5}{\tilde{A}{^{\otimes 4}}} & \rnode{c6}{\ldots} \\[1cm]
\rnode{d1}{ } & \rnode{d2}{0} & \rnode{d3}{0} & \rnode{d4}{0} & \rnode{d5}{0} & \rnode{d6}{ }
\end{array}
\psset{nodesep=3pt}
\ncline{->}{b2}{b1}
\ncline{->}{b3}{b2}\Mput[b]{b'_2}
\ncline{->}{b4}{b3}\Mput[b]{b'_3}
\ncline{->}{b5}{b4}\Mput[b]{b'_4}
\ncline{->}{c2}{c1}
\ncline{->}{c3}{c2}\Mput[t]{b'_{*2}}
\ncline{->}{c4}{c3}\Mput[t]{b'_{*3}}
\ncline{->}{c5}{c4}\Mput[t]{b'_{*4}}
\ncline{->}{a2}{b2}
\ncline{->}{b2}{c2}\Mput[l]{ }
\ncline{->}{c2}{d2}
\ncline{->}{a3}{b3}
\ncline{->}{b3}{c3}\Mput[l]{\Carre_{*2}}
\ncline{->}{c3}{d3}
\ncline{->}{a4}{b4}
\ncline{->}{b4}{c4}\Mput[l]{\Carre_{*3}}
\ncline{->}{c4}{d4}
\ncline{->}{a5}{b5}
\ncline{->}{b5}{c5}\Mput[l]{\Carre_{*4}}
\ncline{->}{c5}{d5}
\]
\end{center}
where for all $n$, $b'_n := b'$.
We denote for all $n>1$, $\Xi_n: A ^{\otimes (n+2)} \xrightarrow{} \tilde{A} ^{\otimes n}
\ \ (a_0, \ldots, a_{n+1}) \mapsto a_0 \otimes (a_1 \ldots a_n) \otimes a_{n+1}$
and $\Carre_{*n} := \Xi_n \circ (\underbrace{d \otimes id \ldots id \otimes d}_
{n \ \ \textrm{terms}})$
Let us denote $f_2 = \omega$, $f_3 = \hat{\delta}\omega$, and for all $n>3$
\begin{eqnarray*}
f_{2n} &=& \hat{\delta}\omega \rho^{2n -3} + \rho^{2}\hat{\delta}\omega\rho^{2n -5} + \rho^{4}\hat{\delta}\omega \rho^{2n -7} + \rho^{6}\hat{\delta}\omega \rho^{2n -9}
+ \ldots + \rho^{2(n - 2)}\hat{\delta}\omega\rho + \rho^{2(n-1)}\omega, \\
f_{2n + 1} &=& \hat{\delta}\omega \rho^{(2n +1) -3} + \rho^{2}\hat{\delta}\omega \rho^{(2n+1) -5}
+ \rho^{4}\hat{\delta}\omega \rho^{(2n+1) -7}
+ \rho^{6}\hat{\delta}\omega\rho^{(2n+1) -9} + \ldots + \rho^{(2n +1) -3} \hat{\delta}\omega.
\end{eqnarray*}
Then $\forall n >1, \  \tilde{f}_n \Carre_{*n} = 0$. Moreover $\hat{\delta}f_{2n} = f_{2n + 1}$
and $\hat{\delta}f_{2n +1} =0$.
\end{theo}
\Rk
For $n=2$ it is already proved. All the proposition comes from the following remark:
\begin{eqnarray*}
a_0 \rho(a_1)\rho(a_2) &=& \omega(a_0,a_1)\rho(a_2) \\
& & - \rho(a_0) \omega(a_1,a_2) \\
& & + \rho(a_0) \rho(a_1)a_2.
\end{eqnarray*}
This equality shows the usefulness of the Bianchi identity when its computed on an Ito derivative.
The idea is then to find an operator, here $\Carre_{*3}$, and $f_3$, here $\widetilde{\hat{\delta} \omega}$
such that $\widetilde{\hat{\delta} \omega}\Carre_{*3}(a_0, a_1, a_2) = 0$.
\Proof
We now fix $n>2$ and remark that:
\begin{eqnarray*}
\label{yyy}
a_0 \rho(a_1) \ldots \rho(a_{2n-1}) &=& \hat{\delta}\omega(a_0,a_1,a_2)\rho(a_3) \ldots \rho(a_{2n - 1})\\
& & + \rho(a_0)\rho(a_1) \hat{\delta}\omega(a_2,a_3,a_4)\rho(a_5) \ldots \rho(a_{2n -1}) + \ldots \\
& & + \rho(a_0) \ldots \rho(a_{2n - 3}) \omega(a_{2n - 2},a_{2n - 1}) \ \ \ \ \ \ \ \ \ \ \ \ \ \ \ \ \ \ \ \ (4)\\
& & + \rho(a_0) \ldots \rho(a_{2n - 2}) a_{2n - 1}.
\end{eqnarray*}
which explains the definition of $f_{2n}$.
Moreover,
\begin{eqnarray*}
\Carre_{*n} &=& \Xi_n \circ (d \otimes id \ldots id \otimes d) \\
&=& + a_0, (1, a_1, \ldots, a_{2n-2},a_{2n-1}), 1 \ \ \ \ (5) \\
& & - a_0, (1, a_1, \ldots, a_{2n-2}, 1), a_{2n-1} \ \ \ \ (6) \\
& & - 1, (a_0, a_1, \ldots, a_{2n-2},a_{2n-1}), 1  \ \ \ \ (7) \\
& & + 1, (a_0, a_1, \ldots, a_{2n-2}, 1), a_{2n-1} \ \ \ \ (8) \\
\end{eqnarray*}
Since $\rho(1) = 0$, we note that $\widetilde{\hat{\delta}\rho^{2n -3}}$
yields zero when applied on equations $(6)$ and $(8)$.
When applied to $(5)$ we obtain $ + a_0\omega(1,a_1)\rho(a_2) \ldots \rho(a_{2n -1})$. \\
On $(7)$ we get
$ - \hat{\delta}\omega(a_0,a_1,a_2)\rho(a_3) \ldots \rho(a_{2n -1})$.
The other terms of the definition of $f_{2n}$ only apply on equation $(7)$ and give
the other terms of the sum (\ref{yyy}), except the last term
$ \widetilde{\rho^{2n -2}\omega}$ which, when applied on equations $(7)$ and $(8)$, yields the two last terms
of the sum (\ref{yyy}).
We obtain $\tilde{f}_{2n} \Carre_{*2n}(a_0, \ldots, a_{2n-1}) = 0 $ i.e. the sum (\ref{yyy}).
The same remark is used to prove the odd case.
To prove $\hat{\delta}f_{2n}  = f_{2n +1}$ we remark that for all
$n>2$, $f_{2n + 1} = \hat{\delta}(\rho^{2n})$, hence $\hat{\delta}f_{2n +1} =0$. Moreover,
$f_{2n} = \hat{\delta}f_{2n -1} \rho + \rho^{2n-2}\omega$.
\begin{eqnarray*}
\hat{\delta}f_{2n} &=& \hat{\delta}(\hat{\delta}f_{2n -1} \rho + \rho^{2n-2}\omega) \\
&=& ((-1)^{2n-1} \hat{\delta}f_{2n -1} \hat{\delta}\rho) + (\hat{\delta}\rho^{2n-2}\omega + \rho^{2n-2}\hat{\delta}\omega) \\
&=& (-\hat{\delta}\rho^{2n-2} (\omega - \rho^2) + (\hat{\delta}\rho^{2n-2}\omega + \rho^{2n-2}\hat{\delta}\rho^2) \\
&=& \hat{\delta}\rho^{2n-2}\rho^2 + \rho^{2n-2}\hat{\delta}\rho^2 \\
&=& \hat{\delta}\rho^{2n} = f_{2n +1}.
\end{eqnarray*}
\Rk
We must be careful about sign. Recall that
$$(f.g)(a_1, \ldots, a_{p+q}) = (-1)^{pq}f(a_1,\ldots,a_p)g(a_{p+1}, \ldots a_{p+q}).$$
In our case $f_n$ are defined up to a sign, without importance for the result. We restore
the right sign by noticing that
for all $n>0$, $f'_{4n} \equiv -f_{4n}$ and $f'_{4n +1} \equiv -f_{4n +1}$. All the
other $f'_{4n+l} \equiv f_{4n+l}$, with $l=2 \ \textrm{or} \ 3$, are correctely defined.
\Rk
Nothing has been said about $b'_{*n}$.
We define $b'_{*n}:= \Carre_{*(n-1)}(m \otimes id \ldots id \otimes m)(id \otimes b_{n} \otimes id) J$,
where, as usual $m$ denotes the product of $A$. The aim of the projection $J$ is to select the equation $(6)$ among the
four possibilities $(5), (6), (7), (8)$
(because for the sequel of this remark the four equations carry the same information) that is,
\begin{eqnarray*}
 \Carre_{*n}(a_0, \ldots, a_n) &=&
 \Xi_n \circ (d \otimes id \ldots id \otimes d)(a_0, \ldots, a_n) \\ &\xrightarrow{J}&
\Xi_n \circ (\delta \otimes id \ldots id \otimes \tilde{\delta})(a_0, \ldots, a_n) \\
 & & = (a_0, (1, a_1, \ldots, a_{n-2}, 1), a_{n-1}).
\end{eqnarray*}
We recall that, $\delta$ and $\tilde{\delta}$ are the coproducts of the unital algebra $A$.
Now if we prove that
$$(m \otimes id \ldots id \otimes m)(id \otimes b_{n} \otimes id)
(a_0, (1, a_1, \ldots, a_{n-2}, 1), a_{n-1}) =
b_{n}(a_0, a_1, \ldots, a_{n-2}, a_{n-1}),$$
the commutativity of the sequence above will be proved and since
$$ b'_{*n}(b'_{*(n+1)}\Carre_{*(n+1)}) = b'_{*n}(\Carre_{*n}b'_{(n+1)}) = (\Carre_{*(n-1)}b'_{n})b'_{(n+1)} = 0.$$
showing that $b'_{*n}b'_{*(n+1)} = 0$. The sequence will be a complex as claimed in the theorem.
However by definition,
$$ b'_n = m \otimes id \ldots \otimes id - id \otimes m \ldots \otimes id + \ldots
+ (-1)^{n+1} \ id \otimes id \ldots \otimes m. $$
Hence,
\begin{eqnarray*}
(m \otimes id \ldots id \otimes m)(id \otimes b'_n \otimes id) &=&  \\
& & m(id \otimes m) \otimes id \ldots \otimes id \otimes m(id \otimes id) \\
& & - m(id \otimes id) \otimes m \ldots \otimes id \otimes m(id \otimes id) \\
& & + \ldots + (-1)^{n+1}
\ m(id \otimes id) \otimes id \otimes id \ldots id \otimes m(m \otimes id).
\end{eqnarray*}
However,
\begin{itemize}
\item {$m(id \otimes m)\delta(a_0) \otimes a_1 = m(a_0 \otimes a_1),$}
\item {$(id \otimes m) a_{(n-2)} \otimes \tilde{\delta} a_{(n-1)} = a_{(n-2)} \otimes a_{(n-1)}, $}
\end{itemize}
proving that
$$(m(id \otimes m) \otimes id \ldots \otimes id \otimes m
(id \otimes id))(\delta \otimes id \ldots id \otimes \tilde{\delta})(a_0, \ldots, a_n) =
(m \otimes id \ldots \otimes id)(a_0, \ldots, a_n).$$
\begin{itemize}
\item {$m(id \otimes id)\delta(a_0) = a_0,$}
\item {$(id \otimes m) a_{(n-2)} \otimes \tilde{\delta} a_{(n-1)} = a_{(n-2)} \otimes a_{(n-1)}, $}
\end{itemize}
proving that
$$- (m(id \otimes id) \otimes m \ldots \otimes id \otimes m
(id \otimes id))(\delta \otimes id \ldots id \otimes \tilde{\delta})(a_0, \ldots, a_n) = -
(id \otimes m \ldots \otimes id)(a_0, \ldots, a_n),$$
and the equality between all the other terms of the sum except the last one. Yet,
\begin{itemize}
\item{$m(m \otimes id) a_{(n-2)} \otimes \tilde{\delta}a_{(n-1)} = m( a_{(n-2)} \otimes a_{(n-1)}),$}
\end{itemize}
proving that
$
((-1)^{n+1}
\ m(id \otimes id) \otimes id \otimes id \ldots id \otimes m(m \otimes id))
(\delta \otimes id \ldots id \otimes \tilde{\delta})(a_0, \ldots, a_n)
= ((-1)^{n+1} \ id \otimes id \ldots \otimes m)(a_0, \ldots, a_n).$
This concludes the proof.
\eproof
\Rk
Related to the Hochschild complex with boundary $b'$ is
the following complex,
$$ 0 \xrightarrow{} \textsf{Hom}(\tilde{A}, A)
\xrightarrow{\hat{\delta}} \textsf{Hom}(\tilde{A} ^{\otimes 2}, A)
\xrightarrow{\hat{\delta}} \textsf{Hom}(\tilde{A} ^{\otimes 3}, A)
\xrightarrow{\hat{\delta}} \textsf{Hom}(\tilde{A} ^{\otimes 4}, A) \xrightarrow{\hat{\delta}} \ldots $$
It is worth noticing that the equations above admit Ito derivatives $\rho$ but also
Leibnitz derivative $D$
if we associate formally with the curvature $\omega$ the bilinear map $Dm$.
\subsubsection{An Ito graded differential algebra}
Let $A$ be an unital algebra.
This short subsection is an attempt to adapt what was done in the case of cyclic cocycles to the Ito case. Connes \cite{Connes}
defines an operator $F$ such that $F^2=I$ and constructs a Leibnitz
graded differential algebra $\Omega^*$ from the $a_0[F,a_1] \ldots [F,a_n]$. The
differential $d$ acts as $a_0[F,a_1] \ldots [F,a_n] \mapsto [F,a_0][F,a_1] \ldots [F,a_n]$
and the product is based on the Leibnitz property, (where we denote $D(a)=[F,a]$),
\begin{eqnarray*}
a_0 D(a_1) \ldots D(a_k)a_{k+1} &=& \\
 & & (-1)^k \sum_{j=1}^k (-1)^j a_0 D(a_1) \ldots D(a_{i-1})D(a_{i}a_{i+1})D(a_{i+2}) \ldots D(a_{k+1}) \\
&+& (-1)^k a_0 a_1 D(a_2) \ldots D(a_{k+1}), \ \ \ \forall a_j \in A.
\end{eqnarray*}
The product of two forms is then associative thanks to the Leibnitz property.

Recall that if $\rho: A \xrightarrow{} A$ is an Ito map, i.e.
$\rho(I)=0, \ \ \rho(ab)= \rho(a)b + a\rho(b) + \rho(a)\rho(b)$, with curvature $\omega$, we get
$\omega(a_1,a_2) = a_1 \rho(a_2) + \rho(a_1)a_2$.
This property is nearly the same that the Leibnitz one. An idea would be to replace
$D(a)$ by $\rho(a)$ and construct the space of forms from,
$$ \eta = a_0 \rho(a_1) \ldots \rho(a_k), \ \ \  \forall a_j \in A.$$
There are two drawbacks in this naive framework. The first one is that
$$\rho(ab)= \rho(a)b + a\rho(b) + \rho(a)\rho(b),$$
mixing 1-forms and 2-forms.
However, since $\omega(a,I) = \omega(I,a) = \rho(a)$, we get rid of this obstacle
by denoting $\rho(a) \equiv \omega(a)$. Hence
we define for $k > 0$, $\Omega^k$ the linear span of the operators
$$ \eta = a_0 \omega(a_1) \ldots \omega(a_k), \ \ \ \forall  a_j \in A.$$
For $k=0$, $\Omega^0 = A$.
Then the $k$-vector space $\Omega^*$ is defined as $\Omega^* := \bigoplus \Omega^k$.
As in the Leibnitz case we remark,
\begin{eqnarray*}
a_0 \omega(a_1) \ldots \omega(a_k)a_{k+1} &=& \\
 & & (-1)^k \sum_{j=1}^k (-1)^j a_0 \omega(a_1) \ldots \omega(a_{i-1})\omega(a_{i},a_{i+1})\omega(a_{i+2}) \ldots \omega(a_{k+1}) \\
&+& (-1)^k a_0 a_1 \omega(a_2) \ldots \omega(a_{k+1}), \ \ \ \forall a_j \in A.
\end{eqnarray*}
Thus $\Omega^*$ can be embedded into an algebra structure.
Here is the second drawback.
The product
of two forms is no longer associative. Nevertheless, we give the following result.

If $\eta_1 \in \Omega^{k_1}$ and $\eta_2 \in \Omega^{k_2}$ then $\eta_1 \eta_2 \in \Omega^{k_1 + k_2}$
and, for all $k$, $\Omega^{k}$ is a $A$-bimodule.
The differential $d: \Omega^* \xrightarrow{} \Omega^*$ is defined as follow:
$$ d(a_0 \omega(a_1) \ldots \omega(a_k)) =  \omega(a_0) \omega(a_1) \ldots \omega(a_k).$$
\begin{prop}
By construction $d^2=0$ and
$$ d(\eta_1 \eta_2) = d(\eta_1)\eta_2 + (-1)^{k_1} \eta_1 d(\eta_2) + (-1)^{k_1} d(\eta_1) d(\eta_2), \ \ \forall \eta_j \in \Omega^{k_j}. $$
\end{prop}
\Proof
$d^2 =0$ is straightforward since $\omega(I) = 0$.
Let $\eta_1 = a_0 \omega(a_1) \ldots \omega(a_{k_1})$ and $\eta_2 = b_0 \omega(b_1) \ldots \omega(b_{k_2})$. For convenience
we rename for all j, $b_j$ in $a_{k_1 + j + 1}$ so that $\eta_2 = a_{k_1 + 1} \omega(a_{k_1 + 2}) \ldots \omega(a_{k_1 +k_2 + 1})$. Then,
\begin{eqnarray*}
d(\eta_1 \eta_2) &=& (-1)^{k_1} \sum_{j=1}^{k_1} (-1)^j \omega(a_0) \rho(a_1) \ldots \omega(a_{i-1})\omega(a_{i},a_{i+1})\omega(a_{i+2}) \ldots \omega(a_{k_1}) \omega(a_{k_1 + 1}) \ldots \omega(a_{k_1 +k_2 + 1}), \\
&+& (-1)^{k_1} \omega(a_0 a_1) \omega(a_2) \ldots \omega(a_{k_1}) \omega(a_{k_1+ 1}) \ldots \omega(a_{k_1 +k_2 + 1}) \\
\end{eqnarray*}
Yet by definition, $\omega(a_0 a_1)= \omega(a_0, a_1) -\omega(a_0)\omega(a_1)$ and
\begin{eqnarray*}
d(\eta_1 \eta_2) &=& (-1)^{k_1} \sum_{j=0}^{k_1} (-1)^j \omega(a_0) \omega(a_1) \ldots \omega(a_{i-1})\omega(a_{i},a_{i+1})\omega(a_{i+2}) \ldots \omega(a_{k_1}) \omega(a_{k_1 + 1}) \ldots \omega(a_{k_1 +k_2 + 1}), \\
&+& (-1)^{k_1} \omega(a_0)\omega(a_1) \omega(a_2) \ldots \omega(a_{k_1}) \omega(a_{k_1+ 1}) \ldots \omega(a_{k_1 +k_2 + 1}). \\
\end{eqnarray*}
However,
\begin{eqnarray*}
d(\eta_1) \eta_2 &=& \omega(a_0)\omega(a_1) \omega(a_2) \ldots \omega(a_{k_1}) a_{k_1 + 1} \omega(a_{k_1 + 2}) \ldots \omega(a_{k_1 +k_2 + 1}) \\
&=& (-1)^{k_1} \sum_{j=0}^{k_1} (-1)^j \omega(a_0) \omega(a_1) \ldots \omega(a_{i-1})\omega(a_{i},a_{i+1})\omega(a_{i+2}) \ldots \omega(a_{k_1 +k_2 + 1}), \\
&-& (-1)^{k_1} a_0\omega(a_1) \omega(a_2) \ldots \omega(a_{k_1}) a_{k_1 + 1} \omega(a_{k_1 + 2}) \ldots \omega(a_{k_1 +k_2 + 1}).\\ \\
(-1)^{k_1} \eta_1 d(\eta_2) &=&  (-1)^{k_1} a_0 \omega(a_1) \omega(a_2) \ldots \omega(a_{k_1}) \omega(a_{k_1 + 1}) \omega(a_{k_1 + 2}) \ldots \omega(a_{k_1 +k_2 + 1}).\\
(-1)^{k_1} d(\eta_1) d(\eta_2) &=& (-1)^{k_1} \omega(a_0) \omega(a_1) \omega(a_2) \ldots \omega(a_{k_1}) \omega( a_{k_1 + 1}) \omega(a_{k_1 + 2}) \ldots \omega(a_{k_1 +k_2 + 1}).
\end{eqnarray*}
\eproof
\Rk
We could define then a graded curvature $\omega_d(\eta_1,\eta_2) :=d(\eta_1\eta_2) -(-1)^{k_1}d(\eta_1)d(\eta_2)$.

To recover an associative product, we have to modify the present setting.
\subsubsection{Curvature of an Ito map and (super)-dialgebra}

In part \ref{curv} we have seen that the curvature of a 1-cochain played
the r\^ole of a product. We use this remark to construct from an Ito map, an
anti-$\mathbb{Z}_2$ graded algebra of non commutative forms. Then we show this algebra
is a dialgebra \cite{Loday} and that the products of this dialgebra embed it into a di-superalgebra.
Moreover we show that any (pre)-dialgebra is a dendriform algebra.
In the sequel $(-1)^x$ will mean $(-1)^{\deg(x)}$. Let us start by some definitions.
\begin{defi}{}
A {\it{super-algebra}} \cite{Quillensuper} is a $k$-vector space $S = S_+ \oplus S_-$ of even and odd elements, belonging
respectively to $S_+$ and $S_-$, equipped with an associative product which respect this $\mathbb{Z}_2$ grading,
i.e.
$aa' \in S_{+}$ iff $a$ and $a'$ are both even or both odd and $aa' \in S_{-}$ otherwise.

An {\it{anti-superalgebra}} could then be $As = As_+ \oplus As_{-}$
a $k$-vector space of even and odd elements belonging
respectively to $As_+$ and $As_{-}$, equipped with an associative product such that
$aa' \in As_{-}$ iff $a$ and $a'$ are both even or both odd and $aa' \in As_{+}$ otherwise.
\end{defi}
Here is an example constructed from an associative unital algebra $A$,
with unit $I$ and the curvature $\omega$ of an Ito map $\rho$. Recall that $\omega(I,I)=0$ and $\omega(I,a)=\omega(a,I)=\rho(a).$
Let $\Omega^* = \bigoplus \Omega^k$, with for all $k > 0$, $\Omega^k$ be
the $A$-bimodule constructed over the linear span of the operators:
$$ a_0 \underbrace{\omega(a_1,a_2) \ldots \omega(a_{2k-1},a_{2k})}_{\textrm{k}}a_{2k+1}. $$
For $k=0$, $\Omega^0 = A$.
The product $\star$ is defined from the curvature,
\begin{eqnarray*}
(a_0 \omega(a_1,a_2) \ldots \omega(a_{2k-1},a_{2k})a_{2k+1}) \star
(b_0 \omega(b_1,b_2) \ldots \omega(b_{2l-1},b_{2l})b_{2l+1}) = \\ a_0 \omega(a_1,a_2) \ldots
\omega(a_{2k+1},b_0) \ldots \omega(b_{2k-1},b_{2k})b_{2l+1}.
\end{eqnarray*}
\Rk
To preserve the associativity of the product $\star$,
the product between $0$-forms and other forms is not defined, except for the identity
element $I$.
Notice that the product embeds two forms of degree $k$ and $l$ into a form of degree $k+l+1$.
\begin{defi}{}
The differential $d:\Omega^* \xrightarrow{} \Omega^*$ is defined for 0-forms as $d(a) := \omega(I,a) = I \star a = a \star I$ and
for forms of higher order by
$d(a_0 \omega(a_1,a_2) \ldots \omega(a_{2k-1},a_{2k})a_{2k+1})
= \\ \omega(I,a_0) \omega(a_1,a_2) \ldots \omega(a_{2k-1},a_{2k})a_{2k+1} + (-1)^k
a_0 \omega(a_1,a_2) \ldots \omega(a_{2k-1},a_{2k})\omega(a_{2k+1},I)$.
\end{defi}
\begin{prop}
The operator $d$ verifies for all $x, y \in \Omega^{\deg(x)} \times \Omega^{\deg(y)}$,
$$ d^2=0 \ \ \textrm{and} \ \ \  d(x \star y) = d(x) \star y + (-1)^{x+1}x \star d(y)=d(x) \star y -(-1)^{x}x \star d(y).$$
\end{prop}
\Proof
Recall that $\omega(I,I)=0$ and denote $k =\deg(x)$
and $x=a_0 \omega(a_1,a_2) \ldots \omega(a_{2k-1},a_{2k})a_{2k+1}$. We have,
$d^2(x)=
d(I\omega(I,a_0) \omega(a_1,a_2) \ldots \omega(a_{2k-1},a_{2k})a_{2k+1} + (-1)^k
a_0 \omega(a_1,a_2) \ldots \omega(a_{2k-1},a_{2k}) \omega(a_{2k+1},I)I)=
+ (-1)^{k+1}I\omega(I,a_0) \ldots \omega(a_{2k+1},I)
+ (-1)^k \omega(I,a_0)  \ldots \omega(a_{2k+1},I)I =0.$ \\
The remaining property follows by straightforward computations.
\eproof
\begin{defi}{[Dialgebra \cite{Loday}]}
A {\it{dialgebra}} $D$ is a $k$-vector space equipped with two associative operations $\dashv$ and
$\vdash$, called respectively right and left products, satisfying 3 more axioms:
\begin{enumerate}
\item{$x \dashv (y \dashv z) = x \dashv (y \vdash z),$}
\item{$(x \vdash y) \dashv z = x \vdash (y \dashv z),$}
\item{$(x \dashv y) \vdash z = (x \vdash y) \vdash z.$}
\end{enumerate}
\end{defi}
Based on this idea we define,
\begin{defi}{}
A \textit{pre-dialgebra} of type I, (respectively of type III), is
a $k$-vector space
equipped with two associative products verifying all the conditions of a dialgebra
but maybe the last one (respectively the first one).
\end{defi}
\begin{theo}
Defining $x \dashv y := -x \star d(y)$ and
$x \vdash y := (-1)^{x+1}d(x) \star y$, $\Omega^*$ is embedded
into a (non unital) dialgebra.
\end{theo}
\Proof
Let $x,y,z \in \Omega^*$. By $xy$ we mean $x \star y$.
The associativity is straightforward,
$x \dashv (y \dashv z) = x \dashv ydz = xd(ydz)=xdydz$ and
$(x \dashv y) \dashv z = xdy \dashv z = xdydz$. For the other product,
$x \vdash (y \vdash z) = x \vdash (-1)^{y+1}d(y)z = (-1)^{x+1}(-1)^{y+1}dxd(y)z =(-1)^{x+y}dxd(y)z$ and
$(x \vdash y) \vdash z = (-1)^{x+1}(dx)y \vdash z = (-1)^{x+1}(-1)^{(x+(y+1) + 1)+1}d((dx)y)z =
(-1)^{x+1}(-1)^{((x+1)+y+1)+1}(-1)^{x+2}dxdyz= (-1)^{x+y}dxd(y)z $. Here we must be careful
with minus sign. $d(x)y$ is a $x+y+1+1$ form because $d$ maps forms of degree $k$ into
forms of degree $k+1$ as does the product itself too.
The end of the proof is left to the reader.
\eproof
\Rk
In our case, we have $(-1)^{x+1}d(x \star y)= x \vdash y - x \dashv y$. The fact that $d^2=0$ must be closely related
to the relation between the two products $(\vdash, \dashv)$. Instead of
defining the differential $d$ and showing that $d^2=0$, we could start from a dialgebra
and define the differential as $x \succ y = x \vdash y - x \dashv y$. The
axioms of a dialgebra entail that for all $(a,b) \in D$,
$$ b \dashv (x \vdash y - x \dashv y) =0, \ \ \ \ \  (x \vdash y - x \dashv y) \vdash a =0.$$
Hence the vector-space of elements $z :=x \vdash y - x \dashv y$ is a right-nilpotent
space for the law $\vdash$ and a left nilpotent space for the law $\dashv$.
\begin{defi}{[Leibnitz algebra \cite{Loday}]}
By definition, a {\it{Leibnitz algebra}} is an algebra equipped with a bracket satisfying
the Leibnitz identity:
$$ [[x,y],z] = [[x,z],y] + [x,[y,z]].$$
\end{defi}
\Rk
We recall that $[\cdot,\cdot]$ is a non associative product whose associativity defect
can be controlled by the Leibnitz identity.
\Rk
In the case of a dialgebra, $[x,y]_L= x \dashv y - y \vdash x$ defines a Leibnitz bracket \cite{Loday}.
This bracket verifies also,
\begin{itemize}
\item {$[x,y \dashv z]_L = y \vdash [x,z]_L + [x,y]_L \dashv z = [x,y \vdash z]_L,$}
\item {$[x \dashv y, z]_L = x \dashv [y,z]_L + [x,z]_L \dashv y, $}
\item {$[x \vdash y, z]_L = x \vdash [y,z]_L + [x,z]_L \vdash y. $}
\end{itemize}
In our case $[x,y]_L= (-1)^{y}d(y) \star x - x \star d(y)$.
\begin{defi}{[Dendriform algebra \cite{Loday}]}
A {\it{dendriform algebra}} $E$ is a $k$-vector space equipped with
two binary operations,
$$ \prec \ , \ \succ: E \otimes E \xrightarrow{} E,$$
which satisfy the following axioms:
\begin{enumerate}
\item {$(a \prec b) \prec c = a \prec (b \prec c) + a \prec (b \succ c),$ }
\item {$(a \succ b) \prec c = a \succ (b \prec c),$ }
\item {$(a \prec b) \succ c + (a \succ b) \succ c = a \succ (b \succ c),$ }
\end{enumerate}
for all elements $a,b,c \in E$. $E$ is said commutative iff $a \prec b = b \succ a$. Thanks to
a theorem of Loday, such dendriform algebras are Zinbiel algebras. A Zinbiel algebra is
an algebra whose associativity defect of $\cdot$ product is controlled by:
$$ (x \cdot y) \cdot z - x \cdot (y \cdot z) = x \cdot (z \cdot y).$$
\end{defi}
\Rk
If we define $a*b := a \prec b + a \succ b$, then this new product is associative.
\begin{theo}
\label{hth}
Let $(D, \dashv, \vdash)$ be a pre-dialgebra of type I.
The relations,
$$ a \prec b = a \dashv b, \ \ \ \ \ \ \ \textrm{and} \ \ \ \ \ \ \ a \succ b = a \vdash b - a \dashv b,$$
embed $D$ into a (non commutative) dendriform algebra.
Similarly, if $(D, \dashv, \vdash)$ is a pre-dialgebra of type III,
the relations,
$$ a \succ b = a \vdash b, \ \ \ \ \ \ \ \textrm{and} \ \ \ \ \ \ \ a \prec b = a \dashv b - a \vdash b ,$$
embed $D$ into a (non commutative) dendriform algebra.
Conversely, any dendriform algebra with $\prec$ associative is a pre-dialgebra of type I and
any dendriform algebra
with $\succ$ associative is a pre-dialgebra of type III.
\end{theo}
\Rk
Notice that
the case of a commutative dendriform (or Zinbiel) algebra does not fit the hypothesis of this
theorem since the
sole product is by definition not associative. We do have two different products. Moreover, if we assume one law associative, say
$\succ \ \equiv \ \vdash$, the
second one
 $a \dashv b = a \prec b + a \succ b = a*b$ is also associative. There does exist a compatibility between the proposition and its converse.
Before giving the proof we need two auxiliary results.
\begin{lemm}
Let $(D, \dashv, \vdash)$ be a pre-dialgebra of type I.
With the relations,
$$ a \prec b = a \dashv b, \ \ \ \ \ \ \ \textrm{and} \ \ \ \ \ \ \ a \succ b = a \vdash b - a \dashv b,$$
The first,  respectively the second, axiom of a pre-dialgebra of type I is equivalent to the first , respectively the second,
axiom of a dendriform algebra with $\prec$ associative.
\end{lemm}
\Proof
$ a \prec (b \succ c) =0  \Leftrightarrow a \dashv ( b \vdash c - b \dashv c) =0$, which proves that the first axiom of
a pre-dialgebra of type I is
equivalent to the first axiom of a dendriform algebra.
Similarly, $(a \succ b) \prec c - a \succ (b \prec c) =0 \Leftrightarrow (a \vdash b) \dashv c - a \vdash (b \dashv c) =0$, since
the product $ \prec \ \equiv  \ \dashv $ is supposed to be associative.
\eproof
\begin{lemm}
Let $(D, \dashv, \vdash)$ be a pre-dialgebra of type III.
With the relations,
$$ a \succ b = a \vdash b, \ \ \ \ \ \ \ \textrm{and} \ \ \ \ \ \ \ a \prec b =  a \dashv b - a \vdash b,$$
The third,  respectively the second, axiom of a pre-dialgebra of type III is equivalent to the third, respectively the second,
axiom of a dendriform algebra with $\succ$ associative.
\end{lemm}
\Proof
The proof is the same. We verify that $(a\prec b) \succ c =0 \Leftrightarrow (a \dashv b) \vdash c - (a \vdash b) \vdash c =0$ and so on.
\eproof
\Proof (of the theorem \ref{hth})
In the case of
a pre-dialgebra of type I, the proof is completed by noticing that
the third axiom of a dendriform algebra with $\prec$ associative
is not enough to prove the third axiom of
a dialgebra. Therefore from the axioms of a dendriform algebra with $\prec$ associative, we prove only the axioms of a predialgebra of
type I, i.e. the axioms 1 and 2 of a dialgebra. (Similarly for a dendriform algebra with $\succ$ associative.)
\eproof
\begin{exam}{}
An associative algebra is a trivial dialgebra by defining $a \vdash b =ab =a \dashv b$. It
is also a trivial dendriform algebra. The second law could be described by $a \succ b = 0(ab) = 0$.
\end{exam}
The following proposition allows us to consider for a differential calculus just one product.
\begin{prop}
Let $x,y \in \Omega^*$, we have,
$ d(x \dashv y) = -dx \star dy = d(x \vdash y),  \
d[x,y]_L = d(y) \star d(x) - d(x) \star d(y)= d(x \dashv y - y \dashv x)$.
\end{prop}
\Proof
Obvious.
\eproof
\Rk
In our example the associative products $(\dashv \ , \ \vdash)$ respect the $\mathbb{Z}_2$
grading of $\Omega^*$. In addition to being a dialgebra, $\Omega^*$ is a (di)-superalgebra. If we embed $(\Omega^*,\dashv \ , \ \vdash)$
into a dendriform algebra, with for example $\succ \ \equiv  \ \vdash$ associative,  we will get $d(a \prec b) =  d(a \dashv b - a \vdash b)=0$, i.e.
$(d \prec): \Omega^* \otimes \Omega^* \xrightarrow{}\Omega^*$ will give closed forms.
Notice that in the case of a graded Leibnitz algebra, the (associative)
Fedosov product \footnote{In the case of a graded Leibnitz algebra, the product of two forms $x,y$ defined by $xy \pm (-1)^{\deg x}dxdy$
is associative. This product is the Fedosov product when the minus sign is chosen.} turns it into a superalgebra. It is also a di-algebra with
$x \dashv y := xd(y)$ and $x \vdash y := (-1)^{x}d(x)y$. It is not a
di-superalgebra but a di-anti-superalgebra.
\newcommand{\lTr}{\textsf{Tr}}
\begin{theo}
The $k$-vector space $(\Omega^*, \dashv)$ is a (di)-superalgebra. For all $x \in \Omega^*$, we define the linear map
$x \mapsto \lTr(x) = \sigma^\natural( d(x) \star I)$, where $\sigma$ is a trace on $A$ and $\natural$ is the
universal cotrace defined in \cite{Quillen}. In this case $\lTr$ is a closed
trace on $(\Omega^*, \dashv)$ and vanishes on the Leibnitz commutator, i.e. $\lTr  [x,y]_L =0$.
\end{theo}
\Proof
Let $x,y \in (\Omega^*, \dashv)$ we have $\lTr (dx)  =0$ since $d^2=0$. If
$x:=a_0 \omega(a_1,a_2) \ldots \omega(a_{2n-1},a_{2n}) a_{2n+1}$ and $y:=b_0 \omega(b_1,b_2) \ldots \omega(b_{2m-1},b_{2m}) b_{2m+1}$,
we get:
$$d(x) \star d(y) \star I = \omega(a_0, I) \ldots  \omega(a_{2n+1},I) \omega(b_0,I) \ldots \omega( b_{2m+1},I).$$
Thus $\lTr(x \dashv y) = -\tau^\natural (d(x) \star d(y) \star I)
= -(-1)^{2(x+2)2(y+2)}\tau^\natural (d(y) \star d(x) \star I) = \lTr(y \dashv x)$ since in the graded algebra defined
at the begining of this section the curvature $\omega$ is a two-cochain.
\eproof
\begin{coro}
Let $x,y \in (\Omega^*, \dashv)$,
then $\lTr(x)$ is a cyclic cocycle of degree $2(\deg(x)+2) -1$ and
$\lTr \ (x \dashv y)$,
is a cyclic cocycle of degree $2(\deg(x)+ \deg(y)+4) -1$.
\end{coro}
\Proof
This is a consequence from \cite{Quillen},  the action $\sigma^\natural$ is recalled in \ref{curv}.
\eproof
\Rk
Let $w=(a_0,(I,a_1,I,\ldots,a_{n-2}),a_{n-1})$ be a pattern of a periodical orbit on the flower graph.
Define $x([w]) = a_0 \omega(I,a_1)\omega(I,a_2) \ldots \omega(I,a_{n-2})a_{n-1} \in \Omega^*$, then
the functions $f_n$ defined in section \ref{fn} can be expressed in term of $dx([w])$.
Since $\omega(I,a_1) = \omega(a_1,I)=\rho(a_1)$ we find that,
\begin{eqnarray*}
dx([w]) &=& \omega(I,a_0) \omega(I,a_1)\omega(I,a_2) \ldots \omega(I,a_{n-2})a_{n-1} \\
& & + (-1)^{(n-2)}\omega(I,a_1)\omega(I,a_2) \ldots \omega(I,a_{n-2})\omega(I,a_{n-1}) \\
&=& \rho(a_0) \rho(a_1)\rho(a_2) \ldots \rho(a_{n-2})a_{n-1}
 +(-1)^n a_0\rho(a_1) \rho(a_2) \ldots \rho(a_{n-2}) \rho(a_{n-1}),
\end{eqnarray*}
to be compared with $\tilde{f_n}\Carre_{*n}(a_0,a_1,\ldots,a_{n-2},a_{n-1}) =0$.
Hence what we have obtained in section \ref{fn} can still be rediscover by considering periodical
orbits on the flower graph.
\subsubsection{A complex of non local forms}
\textbf{Notation:} In this subsection, $\sum_i id \otimes \ldots \otimes \textrm{EXP} \otimes \ldots id$ means
that the expression EXP is placed at position $i$.

We wish to study $\Delta -\Delta_f$, i.e. the difference between a coassociative bialgebra $C$ with coproduct $\Delta$
from the coassociative coalgebra of the flower graph $\Delta_f$ generated by the coproduct of the bialgebra.
\begin{lemm}
Let us denote $ \hat{\partial}_0 = \Delta$ and for all $n>0$, $$\hat{\partial}_n =
\sum_{i=1}^{n+1} \ (-1)^{i+1} \ \underbrace{id \otimes \ldots \otimes id \otimes \Delta \otimes id \ldots \otimes id}_{n+1 \ \ \textrm{terms}}.$$
Fix $n>0$ and $(a,a_1, \ldots, a_n) \in A$. Then,
\begin{enumerate}
\item{$\hat{\partial}_{n+1} \tilde{\delta} (a_1) \otimes a_2 \otimes \ldots \otimes a_n =  \Delta(1) \otimes a_1 \otimes a_2 \otimes \ldots \otimes a_n - 1 \otimes \hat{\partial_n}(a_1 \otimes a_2 \otimes \ldots \otimes a_n).  $}
\item{$ (-1)^n \hat{\partial}_{n+1} a_1 \otimes a_2 \otimes \ldots \otimes \delta (a_n) = (-1)^n\hat{\partial_n}(a_1 \otimes a_2 \otimes \ldots \otimes a_n) \otimes 1 + a_1 \otimes a_2 \otimes \ldots \otimes a_n \otimes \Delta(1).$}
\item{$(\tilde{\delta} \otimes id)\tilde{\delta} = \Delta(1) \otimes id = (\Delta \otimes id) \tilde{\delta} $.}
\item{$ (id \otimes \delta)\delta = id \otimes \Delta(1) = (id \otimes \Delta)\delta$. }
\item{$ (\Delta \otimes id)\delta(a)= \Delta(a) \otimes 1, \ \ \ (id \otimes \Delta)\tilde{\delta}(a) = 1 \otimes \Delta(a) $.}
\end{enumerate}
\end{lemm}
\Proof
The proof is complete by noticing that $\Delta(1) = 1 \otimes 1$ and that
the operators $\delta(a) = a \otimes 1, \ \ \tilde{\delta}(a) = 1 \otimes a$ act as shifts in
the Fock space of $A$.
\eproof
\begin{theo}
\label{th}
Recall that $\overleftarrow{d} = \Delta - \tilde{\delta}$ and $\overrightarrow{d} = \Delta - \delta$. The sequence
$$ 0 \xrightarrow{} A \xrightarrow{\Delta -\Delta_f} A^{\otimes 2} \xrightarrow{\partial_1 = \overleftarrow{d} \otimes id - id \otimes \overrightarrow{d}} A^{\otimes 3} \xrightarrow{\partial_2} A^{\otimes 4} \xrightarrow{\partial_3} \ldots$$
with $\forall n>0$,
$$\partial_n := \underbrace{\overleftarrow{d} \otimes id \otimes \ldots \otimes id}_{n+1 \ \ \textrm{terms}} + \sum_{j=2}^{n-1} \ \ (-1)^{j+1} \ id \otimes \ldots \otimes id \otimes
\Delta \otimes id \otimes \ldots \otimes id
+ (-1)^{n + 1}id \otimes id \otimes \ldots \otimes id \otimes \overrightarrow{d},$$ defines a complex.
The boundary operators verify:
\begin{enumerate}
\item {$\forall n>0 \ \ \partial_{n+1}\partial_n = 0$ and $\partial_1(\Delta -\Delta_f)=0$.}
\item{$\forall n>0 \ \ \partial_n (x_1, \ldots , x_n)$ is a multilinear map which is an Ito derivative in the first and last variables and a homomorphism in others variables.}
\end{enumerate}
\end{theo}
\Proof
We only have to prove the first item, since the second one comes from the very definition of the boundary operators.
From coassociativity coalgebra theory we know that $ \hat{\partial}_{k+1} \hat{\partial_k} = 0, \ \ \forall k \in \mathbb{N}$.
Fix $n>0$.
\begin{eqnarray*}
\partial_{n+1} \partial_n &=& (\hat{\partial}_{n+1} - \tilde{\delta} \otimes id \ldots \otimes id +(-1)^{n+1} id \otimes \ldots \otimes id \otimes \delta)
(\hat{\partial_n} \\
& & - \tilde{\delta} \otimes id \ldots \otimes id + (-1)^{n} id \otimes \ldots \otimes id \otimes \delta) \\
&=& \hat{\partial}_{n+1}\hat{\partial_n} - \hat{\partial}_{n+1}(\tilde{\delta} \otimes id \ldots \otimes id) + (-1)^{n} \hat{\partial}_{n+1}( id \otimes \ldots \otimes id \otimes \delta )\\
& &  - 1 \otimes \hat{\partial_n} + (\tilde{\delta} \otimes 1)(\tilde{\delta}) \otimes id \ldots \otimes id\\
& & +(-1)^{n+1} \hat{\partial_n} \otimes 1 - id \otimes \ldots \otimes id \otimes (id \otimes \delta)(\delta).
\end{eqnarray*}
The equality $\partial_1(\Delta -\Delta_f) = 0$ is left to the reader.
\eproof
\begin{theo}
For all $n > 0$, $\partial_n = \sum_{i=1}^{n + 1} \ (-1)^{n+1} \ id \otimes \ldots \otimes (\Delta - \Delta_f) \otimes  \ldots \otimes id$.
\end{theo}
\Proof
The proof is straightforward by noticing that $id \otimes \tilde{\delta} = \delta \otimes id$.
\eproof

To interpret the two previous theorems, we define:
\begin{defi}{[Non local bundle]}
A non local bundle $B$ consists of a flower graph which plays the r\^ole of the basis, (i.e. a unital
associative algebra embedding into its Markov $L$-coalgebra) and a coassociative coalgebra which
plays the r\^ole of a fiber space.
The projection $\pi$ is $(id \otimes \epsilon) = (\epsilon \otimes id)$ and the section map is
$\sigma = \Delta$.
\end{defi}
\Rk
We have $\pi \circ \sigma = (id \otimes \epsilon)\Delta = id$. The bundle is said non local
because of the non local aspect of the coproduct \footnote{A (local) bundle can also be constructed from a graph $G$, equipped
with a family of probability vectors, since $G$ can be embedded into its natural Markov $L$-coalgebra with counits.}.
\begin{exam}{}
Consider the graph of $Sl(2)_q$. The difference between the fiber over $a$ from the petal over
$a$ is:
\begin{center}
\includegraphics*[width=12cm]{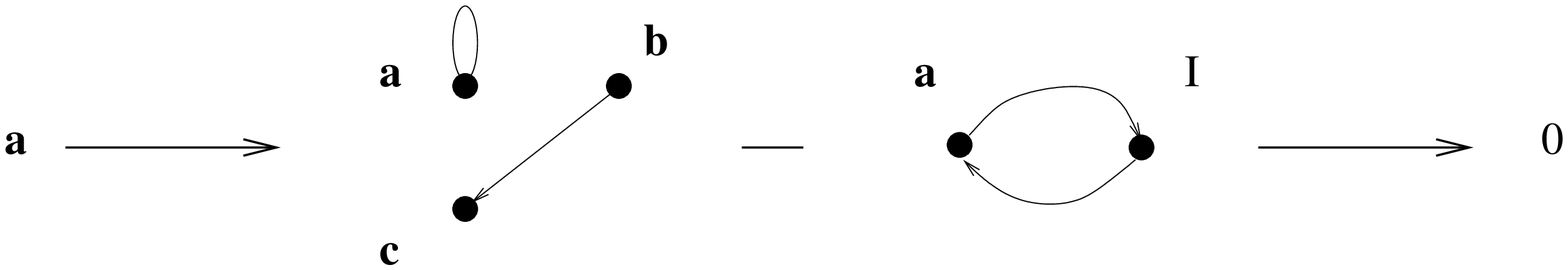}
\end{center}
Which vanishes when the operator $\partial_1$ is applied.
\end{exam}
\begin{defi}{[Primitive element]}
Let $x \in C$ where $C$ is a coassociative coalgebra. An element $x$ is said primitive if
$\Delta(x) = x \otimes 1 + 1 \otimes x = \Delta_f(x)$.
\end{defi}
\Rk
In this framework, only the primitive elements are not disturbed by the section map $\Delta$, because
a petal remains a petal.

Considering a bialgebra as a fiber space allows to interpret the two previous theorems in a physical way.
To obtain a
derivation of a function $f(x_1, \ldots, x_n)$ we must study the partial derivations of
a variable of $f$ and fix the others. The same remark holds here.
To study the deformation between the pattern of a periodic orbit $(a_1, \ldots, a_n)$ living on the
flower graph, i.e. the basis and its lift by the coproduct on the fiber space,
we have to proceed by studying that deformation on a variable by keeping the others fixed.
The previous theorem claims that the deformation of the pattern (or the string),
$(a_1, \ldots, a_n)$ can be studied in the same way as the many-variable functions.
The theorem \ref{th} claims that the boundary operators can be written in an other way such that
inside the string, i.e. $(a_2, \ldots, a_{n-1})$ the structure of the algebra is
respected thanks to the homomorphism aspect of the coproduct and
at the border of the string, i.e. $(a_1,a_n)$, where a physical
interaction might be possible, they behave
as an Ito derivative. \\
\\
We finish this section by the following remark:
\subsubsection{A connection of the Ito derivatives to the third Reidemeister movement: the distributivity defect of an Ito map}
For all $a \in A$, we can define $[a, \cdot]: b \mapsto ab-ba$.
Hence,
$$L: A \xrightarrow{} \textsf{Hom}(A,A), \ \ \ a \mapsto [a, \cdot].$$
If we define $\circ$ such
that $a \circ b := [a, b]$ we get the Leibnitz identity:
$$ (x \circ y) \circ z - x \circ (y \circ z) = (x \circ z) \circ y ,$$
which means that we can control the lack of associativity of the product $\circ$. In the case
of Ito derivatives we showed a bijection between Ito maps and homomorphisms. For any invertible
element of $A$ we can construct the map $I: A \xrightarrow{} \textsf{Hom}(A,A)$ defined by,
$$a \mapsto \rho_a(\cdot) = a(\cdot)a^{-1} - id:= a \star (\cdot).$$
If we define $\star$ such that $a \star b := \rho_a(b)$, we get the Ito identity:
$$ x \star (y \cdot z) - (x \star y)\cdot(x \star z) = (x \star y)\cdot z  + y \cdot (x \star z),$$
which means that the lack of distributivity of the product $\star$ with regard to the
product $\cdot$ can be controlled.
This remark can be used to generalise the definition of the Ito derivative concept.
For example if $A$ is equipped with a product $\rhd$ which verify the third Reidemeister movement,
i.e. $a \rhd (b \rhd c) = (a \rhd b) \rhd (a \rhd c)$ (auto-distributivity) then for all elements of $A$, the map $
I: A \xrightarrow{} \textsf{Hom}(A,A)$ defined by,
$$ a \mapsto \rho_a(\cdot) = a \rhd (\cdot) - id := a \star (\cdot),$$
sends an element $a$ into an Ito derivative since $\Psi_a (c) := a \rhd c$ is a
$\rhd$-homomorphism \cite{Lergraph}. We will have:
$$ x \star (y \rhd z) - (x \star y) \rhd (x \star z) = (x \star y)\rhd z  + y \rhd (x \star z).$$
In some sens if $h$ is a homomorphism from $A$ to $A$ for the usual product, we can say that
$h$ verifies an (auto)-distributivity
condition since for all $a,b \in A$ we have, (if we define $h(a) := h \star a$),
$$ h \star (a \cdot b) - (h \star a)\cdot(h \star b) = h(ab) - h(a)h(b).$$
That is why there is a link between Ito derivatives and homomorphisms.

We can embed the whole algebra $A( m, [\cdot,\cdot])$ into $\textsf{Hom}(A,A)$ in such a way that each point of $A$
can be seen as a Leibnitz derivative thanks to the Poisson or the Lie bracket but it is not possible
in an algebra whose product verifies the third Reidemeister movement. However we can keep
the same embedding, $A(\rhd, \star)$ into $\textsf{Hom}(A,A)$ but in such a way that each point of $A$
can be seen as an Ito derivative thanks to the $\rhd$ product defined above.

It is
a way to connect knot theory to Ito maps.
\section{Towards other botanic specimens...}
\subsection{Differentials}
The aim of this section is to study the generalisation of what was said for the flower graph
to any  $L$-coalgebra and to show that there is a one to one mapping between Leibnitz differentials and Ito differentials.
All the products considered will be associative.
Let $G$ be a $L$-coalgebra, not necessary of Markovian type.
As we will study differentials, we suppose that $G$ is a $L$-bialgebra with unital coproducts $\Delta,\tilde{\Delta}$.
\begin{defi}{}
Let $G$ be a  $L$-bialgebra,
we can equip $G \otimes G$ with a $G$-bimodule structure by:
$$ a(u \otimes v) = \Delta(a)(u \otimes v), \ \ \ \ \
(u \otimes v)a = (u \otimes v)\tilde{\Delta}(a),$$
where $a,u,v \in G$.
We denote $\widetilde{G \otimes G}$ the space so obtained.
However we can equip $G \otimes G$ with another $G$-bimodule structure
by breaking the symmetry of the above definition:
$$ a(u \otimes v) = \tilde{\Delta}(a)(u \otimes v), \ \ \ \ \
(u \otimes v)a = (u \otimes v)\tilde{\Delta}(a).$$
We denote $\widehat{U \otimes V}$ the space so obtained.
\end{defi}
\begin{defi}{}
We denote,
\begin{itemize}
\item{$\widetilde{\Omega}_G$ the $G^{\otimes 2}$-bimodule, seen as a $G$-bimodule in the sense of the first definition,
as the set of non commutative forms of a $L$-bialgebra.}
\item{$\widehat{\Omega}_G$ the $G$-bimodule of non commutative forms of
a $L$-bialgebra, in the sense of the second definition.}
\end{itemize}
\end{defi}
\begin{defi}{[Differentials]}
In the first case, the differential map is the linear surjective map
$ \bar{d}_G: G \xrightarrow{} \widetilde{\Omega}_G $ defined by $ \bar{d}_G = \Delta - \tilde{\Delta} $.
In the second case, the differential map is the linear surjective map
$ \hat{d}_G: G \xrightarrow{} \widehat{\Omega}_G $ with $ \hat{d}_G = \Delta - \tilde{\Delta}$.
\end{defi}
\begin{lemm}
$ \bar{d}_G: G \xrightarrow{} \widetilde{\Omega}_G$ is a Leibnitz differential.
\end{lemm}
\Proof
For $a,b \in G$,
we must show that $\bar{d}_G(ab) = a\bar{d}_G(b) + \bar{d}_G(a)b$.
But
\begin{eqnarray*}
a\bar{d}_G(b) + \bar{d}_G(a)b &=& \Delta(a)(\Delta(b) - \tilde{\Delta}(b)) + (\Delta(a) - \tilde{\Delta}(a))\tilde{\Delta}(b) \\
			&=& \Delta(ab) - \Delta(a)\tilde{\Delta}(b) + \Delta(a)\tilde{\Delta}(b) - \tilde{\Delta}(ab)\\
			&=& \bar{d}_G(ab).
\end{eqnarray*}
\eproof
\begin{prop}
The map $ \hat{d}_G: G \xrightarrow{} \widehat{\Omega}_G $ with $ \hat{d}_G = \Delta - \tilde{\Delta}$
defines an Ito differential.
\end{prop}
\Proof
Let $u,v \in G$.
\begin{eqnarray*}
\hat{d}_G( u) \hat{d}_G (v) &=& (\Delta u -\tilde{\Delta} u)(\Delta v -\tilde{\Delta} v) \\
			&=& \Delta (uv) - ( \tilde{\Delta} (uv) - \tilde{\Delta} (uv) ) - \Delta (u) \tilde{\Delta} (v) - \tilde{\Delta} (u) \Delta (v) + \tilde{\Delta} (uv) \\
			&=& (\Delta (uv) - \tilde{\Delta} (uv)) - ( \Delta (u) - \tilde{\Delta} (u) )\tilde{\Delta}(v) -  \tilde{\Delta}(u)( \Delta (v) - \tilde{\Delta} (v) ) \\
			&=& \hat{d}_G (uv) - \hat{d}_G (u)v - u \hat{d}_G (v).
\end{eqnarray*}
\eproof \\
Let $G$ be a graph. There can be several ways to embed $G$ into a $L$-coalgebra.
Define  $ \textsf{DER(G)}_{\textrm{Leibnitz}} = \{  \ d: G \xrightarrow{} \widetilde{\Omega}_G, d = \Delta - \tilde{\Delta}, \ (\Delta, \tilde{\Delta}) \in \textsf{Prop}(G) \} $
and $ \textsf{DER(G)}_{\Ito} = \{  \ d: G \xrightarrow{} \widehat{\Omega}_G, d = \Delta - \tilde{\Delta}, \ (\Delta, \tilde{\Delta}) \in \textsf{Prop}(G) \} $.
We have shown the following theorem:
\begin{theo}
$$\textsf{DER(G)}_{\textrm{Leibnitz}} \xrightarrow{\textrm{identity}} \textsf{DER(G)}_{\Ito} $$
$$\bar{d} \mapsto \hat{d}$$
is an one to one mapping.
That is a Leibnitz derivative is an Ito one, and conversely if we only change the way we define the $G$-bimodule
structure.
\end{theo}
\begin{exam}{}
For the flower graph of a unital algebra, the universal differential is
$\bar{d} x = dx = x \otimes 1 -1 \otimes x$, which can be viewed as an Ito derivative if we consider
$\hat{\Omega}_A$.
\end{exam}
\Rk
We could define partial derivatives along a walk $w$ by defining $\Delta_{[w]}$, resp.
$\tilde{\Delta}_{[w]}$, equal to $\Delta$, resp. $\tilde{\Delta}$, restricted to the walk $w$.
In this case the partial derivative along a walk is denoted $\partial((\cdot), w) = d_{[w]}(\cdot)$.
\subsection{Commutator}
\begin{defi}{}
As usual, we denote the underlying associative product of $a,b$ by either $ab$ or $m(a,b)$.
Define the commutator associated with a $L$-bialgebra $G$ by:
$$ [\cdot,\cdot]_G: G \times G \xrightarrow{} G^{\otimes 2}, $$
$$ (x, y) \mapsto [x , y]_G = \Delta(x)\tilde{\Delta}(y) - \Delta(y)\tilde{\Delta}(x) .$$
\end{defi}
\Rk
If we consider $G^{\otimes 2}$ equipped with its first $G$-bimodule structure, that is if we identify
$G^{\otimes 2}$ with $\tilde{G}$, we have:
$$ [x , y]_G = x(1 \otimes 1)y - y(1 \otimes 1)x = x.y-y.x, $$
where $.$ denotes the product induced by $m$ on $\tilde{G}$.
In the following we consider $G$ as $\tilde{G}$ and the differential $\tilde{d}_G$ identified with $d$.
\begin{prop}
Let $x,y,a,b \in G$. We have the following identities:
\begin{itemize}
\item{$[\cdot,\cdot]_G $ is bilinear.}
\item{$[x , x]_G = 0.$}
\item{$ [x , ya]_G = [x , y]_G .a + y.[x , a]_G - y.(dx).a$}
\item{$ [xb , y]_G = [x , y]_G .b + x.[b , y]_G + x.(dy).b$}
\end{itemize}
\end{prop}
\Proof
Let $x,y,z,a,b \in G$ and $\lambda \in k$.\\
$ [ \cdot,\cdot ]_G  $ is bilinear since
$\Delta$ and $\tilde{\Delta}$ are linear and
\begin{eqnarray*}
[x , y + \lambda z ]_{G} &=& \Delta (x) (\tilde{\Delta} (y) + \lambda \tilde{\Delta}(z)) - ( \Delta (y) + \lambda \Delta (z)) \tilde{\Delta}(x) \\
&=& \Delta (x) \tilde{\Delta} (y) + \lambda \Delta (x)\tilde{\Delta}(z) - \Delta (y) \tilde{\Delta}(x) - \lambda \Delta (z) \tilde{\Delta}(x) \\
&=& [ x, y]_G + \lambda [ x, z]_G.
\end{eqnarray*}
The same computation can be done for the other side. It is obvious that
$[ x, x]_G = 0$ while
$ [x , ya]_G = [x , y]_G .a + y.[x , a]_G - y.(dx).a$ is obtained by straightforward computation.
Here only the left hand side is computed,
\begin{eqnarray*}
[x , ya]_G &=& [\Delta (x)\tilde{\Delta}(y)\tilde{\Delta}(a) - \Delta (y)\tilde{\Delta}(x)\tilde{\Delta}(a)] + [\Delta (y)\Delta (x)\tilde{\Delta}(a) \\
& & + \Delta (y)\Delta (a)\tilde{\Delta}(x)] - [\Delta (y) ( \Delta (x) - \tilde{\Delta}(x) )\tilde{\Delta}(a)] \\
&=& \Delta (x) \tilde{\Delta} (ya) - \Delta (ya) \tilde{\Delta} (x).
\end{eqnarray*}
The same computation can be done for the last equation.
\eproof
\subsection{The arrow set }
The aim of this part is to develop some feelings about how the two Hudson's propositions can be generalised.
These propositions assert that a coassociative coalgebra can be embedded into an Ito $L$-coalgebra.
If we remove the coassociative coproduct by the right and left coproducts of a Markov $L$-coalgebra, how can we produce two
new coproducts such that, if the old ones are unital homomorphisms, the new ones become Ito derivatives ?

The answer to this question is to find, in the proof of the Hudson's propositions, into the term: $1 \otimes x \otimes 1$.
This term means that we have made one turn on the petal of the flower graph.
In a general graph, there does not exist such a petal. Hence we have to create it. To do so, we need to consider
the arrow set $G_1$ of a graph $G$, that is we must break the spherical property of a point (vertex) and
fix a direction, say $ a \xrightarrow{} b $ still denoted $a \otimes b$. By this way we can construct two
virtual petals, either $ a \xrightarrow{} b \xrightarrow{} a $ or $ b \xrightarrow{} a \xrightarrow{} b $.
For this we need two operators which map $G_1^{\otimes 2}$ into $G_1^{\otimes 3}$. Let us see how it works.
\begin{defi}{}
Let $n \in \mathbb{N} \setminus \{0\}$.
We shall say that $(Z, \Delta_n, \tilde{\Delta}_n, \epsilon_n, \tilde{\epsilon}_n)$ is a
$L$-coalgebra of degree $n$ over $k$ if $Z$ is a $k$-vector space and if it obeys the identities:
\end{defi}
\begin{enumerate}
\item{The following graph
	\begin{equation*}
        \begin{CD}
        Z^{\otimes n} @>\Delta_n>> Z^{ \otimes n+1} \\
        @V{\tilde{\Delta}_n}VV		@VV{\tilde{\Delta}_n \otimes id}V \\
        Z^{ \otimes n+1} @>{id \otimes \Delta_n }>> Z^{ \otimes n+2}
        \end{CD}
        \end{equation*}
commutes, this means: $ (\tilde{\Delta}_n \otimes id)\Delta_n = (id \otimes \Delta_n)\tilde{\Delta}_n. $
     }
\item{\begin{equation*}
	\epsilon_n: Z^{\otimes n} \xrightarrow{} Z^{\otimes n-1} \ \ \textrm{such that:} \  (id \otimes \epsilon_n)\Delta_n = id.
	\end{equation*}}
\item{\begin{equation*}
\tilde{\epsilon}_n: Z^{\otimes n} \xrightarrow{} Z^{\otimes n-1} \ \ \textrm{such that:} \  ( \tilde{\epsilon}_n \otimes id)\tilde{\Delta}_n = id.
        \end{equation*}}
\end{enumerate}
By convention we define: $Z^{\otimes 0} = k$.
\begin{prop}
A Markov $L$-coalgebra $G$ is a Markov $L$-coalgebra of degree $n$, for any $n > 0$.
\end{prop}
\Proof
Let $\Delta, \tilde{\Delta}, \epsilon, \tilde{\epsilon}$ be the coproducts and counits of $G$ and define the following operators:
\begin{eqnarray*}
\Delta_n &=& (\underbrace{id \otimes ... \otimes id}_{n-1} \otimes \Delta ), \ \ \
\tilde{\Delta}_n = ( \tilde{\Delta} \otimes \underbrace{id \otimes ... \otimes id}_{n-1} ), \\
\epsilon_n &=& (\underbrace{id \otimes ... \otimes id}_{n-1} \otimes \epsilon ), \ \ \ \ \
\tilde{\epsilon}_n = ( \tilde{\epsilon} \otimes \underbrace{id \otimes ... \otimes id}_{n-1} ).
\end{eqnarray*}
These embed $G$ into a Markov $L$-coalgebra of degree $n$.
\eproof \\
From now on we consider the special case $n=2$. The Markov $L$-coalgebra $G$ can be embedded into a
Markov $L$-coalgebra of degree 2.
\begin{defi}{}
Define:
$$ \delta_R: G^{ \otimes 2} \xrightarrow{} G^{ \otimes 3}, \ \ \ \ \delta_L: G^{ \otimes 2} \xrightarrow{} G^{ \otimes 3} $$
$$ a \otimes b \xrightarrow{} a \otimes b \otimes a, \ \ \ \ \ a \otimes b \xrightarrow{} b \otimes a \otimes b .$$
\end{defi}
\Rk
If $x = a \otimes b \in G_1$ we notice that
$ \delta_R(x) = x \otimes s(x) $ and $ \delta_L(x) = t(x) \otimes x $  where $s,t$ are the source and terminus maps from the graph $G$.
\begin{prop}
The diagram,
\begin{equation*}
        \begin{CD}
        G^{\otimes 2} @>\delta_R>> G^{ \otimes 3} \\
        @V{\delta_L}VV		@VV{\delta_L \otimes id}V \\
        G^{ \otimes 3} @>{id \otimes \delta_R }>> G^{ \otimes 4}
        \end{CD}
        \end{equation*}
commutes, i.e. $ (\delta_L \otimes id)\delta_R = (id \otimes \delta_R)\delta_L $.
Moreover $\delta_L, \delta_R$ are both homomorphisms if $G$ has an underlying associative algebra structure.
\end{prop}
\Proof
$$ a \otimes b \xrightarrow{\delta_L} b \otimes a \otimes b \xrightarrow{id \otimes \delta_R} b \otimes ( a \otimes b  \otimes a). $$
$$ a \otimes b \xrightarrow{\delta_R} a \otimes b \otimes a \xrightarrow{ \delta_L \otimes id } ( b \otimes  a \otimes b ) \otimes a .$$
If $(a,b,c,d) \in G $ then:
$$\delta_L(a \otimes b) \delta_L(c \otimes d) = ( b \otimes a \otimes b)(d \otimes c \otimes d) =
(bd \otimes ac \otimes bd) = \delta_L (ac \otimes bd) = \delta_L ((a \otimes b)(c \otimes d)). $$
The same computation is used for proving that $\delta_R$ is an homomorphism.
\begin{theo}
If we define two operators $\overleftarrow{d}_G,\overrightarrow{d}_G$ as:
$$\overrightarrow{d}_G = \Delta_2 - \delta_R, \ \ \ \ \
\overleftarrow{d}_G = \tilde{\Delta_2} - \delta_L,$$
they verify $ (\overleftarrow{d}_G \otimes id)\overrightarrow{d}_G  = (id \otimes \overrightarrow{d}_G)\overleftarrow{d}_G. $
\end{theo}
\Proof
Let $ (a, a_0, b, b_1) \in G $.
We denote in symbolic notation $\tilde{\Delta} a = a_0 \otimes a$ and $\Delta b = b \otimes b_1$.
By definition $ \tilde{\Delta}_2 ( a \otimes b ) = a_0 \otimes a \otimes b $ and $ \Delta_2 ( a \otimes b ) = a \otimes b \otimes b_1 $.
Now:
$$ a \otimes b \xrightarrow{\overleftarrow{d}_G} a_0 \otimes a \otimes b - b \otimes a \otimes b
\xrightarrow{id \otimes \overrightarrow{d}_G} a_0 \otimes a \otimes b \otimes b_1 - a_0 \otimes a \otimes b \otimes a
- b \otimes a \otimes b \otimes b_1 + b \otimes a \otimes b \otimes a $$
and
$$ a \otimes b \xrightarrow{\overrightarrow{d}_G} a \otimes b \otimes b_1 - a \otimes b \otimes a
\xrightarrow{\overleftarrow{d}_G \otimes id} a_0 \otimes a \otimes b \otimes b_1 - b \otimes a \otimes b \otimes b_1
- a_0 \otimes a \otimes b \otimes a + b \otimes a \otimes b \otimes a .$$
\eproof
\begin{defi}{[chiral $G^{\otimes 2}$-bimodule]}
Let $x,y,z \in G^{\otimes 2}$.
We can embed $G^{\otimes 3}$ into a $G^{\otimes 2}$-bimodule, called a {\it{chiral $G^{\otimes 2}$-bimodule}},
by equipping $\overleftarrow{d}_G(G^{\otimes 2})$ and $\overrightarrow{d}_G(G^{\otimes 2})$ with the following product:
\begin{eqnarray*}
x.\overleftarrow{d}_G(z) &=& \delta_L(x) \overleftarrow{d}_G(z); \ \ \ \overleftarrow{d}_G(z).y = \overleftarrow{d}_G(z) \delta_L(y), \\
x.\overrightarrow{d}_G(z) &=& \delta_R (x) \overrightarrow{d}_G(z); \ \ \ \overrightarrow{d}_G(z).y = \overrightarrow{d}_G(z)\delta_R(y).
\end{eqnarray*}
\end{defi}
\Rk
The notion of chiral bimodule is well defined because if $\overrightarrow{d}(x \otimes y)=\overleftarrow{d}(x \otimes y)$
then $(\Delta_2 - \tilde{\Delta}_2)(x \otimes y) = (\delta_R - \delta_L)(x \otimes y)$. This happens
if and only if $\Delta_2(x \otimes y) = (x \otimes y) \otimes x$ and
$\tilde{\Delta}_2(x \otimes y)= y \otimes (x \otimes y)$. However
$\overleftarrow{d}_G$ and $\overrightarrow{d}_G$ give
zero on such element $x \otimes y$.
\begin{theo}
Let $G$ is a Markov $L$-bialgebra.
We embed $G^{\otimes 3}$ into its chiral $G^{\otimes 2}$-bimodule.
Then $\overrightarrow{d}_G, \overleftarrow{d}_G$ are Ito derivatives.
\end{theo}
\Proof
For $x,y \in G^{\otimes 2}$,
we compute,
\begin{eqnarray*}
\overleftarrow{d}_G(x)\overleftarrow{d}_G(y) &=&(\tilde{\Delta}_2 (x) - \delta_L(x))(\tilde{\Delta}_2 (y) - \delta_L(y)) \\
&=& \tilde{\Delta}_2(xy) - \tilde{\Delta}_2(x)\delta_L(y) - \delta_L(x)\tilde{\Delta}_2(y) + \delta_L(x)\delta_L(y) + (\delta_L(x)\delta_L(y) - \delta_L(x)\delta_L(y))\\
&=& \overleftarrow{d}_G(xy) - \overleftarrow{d}_G(x)\delta_L(y) - \delta_L(x)\overleftarrow{d}_G(y)\\
&=& \overleftarrow{d}_G(xy) - \overleftarrow{d}_G(x).y -  x.\overleftarrow{d}_G(y).
\end{eqnarray*}
Similarly,
\begin{eqnarray*}
\overrightarrow{d}_G(x)\overrightarrow{d}_G(y) &=& \overrightarrow{d}_G(xy) - \overrightarrow{d}_G(x)\delta_R(y) - \delta_R(x)\overrightarrow{d}_G(y) \\
&=& \overrightarrow{d}_G(xy) - \overrightarrow{d}_G(x).y -  x.\overrightarrow{d}_G(y).
\end{eqnarray*}
\eproof \\
The notion of chiral $G^{\otimes 2}$-bimodule is important if we wish to get the Ito's property.
We can easily generalise the present
setting. Fix $n>1$, and generalise the definition of $\delta_R$ and $\delta_L$ as follow,
\begin{eqnarray*}
\delta_{R,n}: G^{\otimes n} &\xrightarrow{}& G^{\otimes (n+1)}, \ \ \ \ \delta_{L,n}: G^{\otimes n} \xrightarrow{} G^{\otimes (n+1)}, \\
(a_1, \ldots, a_n) &\mapsto& (a_1, \ldots, a_n) \otimes a_1, \ \ \ \ (a_1, \ldots, a_n) \mapsto a_n \otimes (a_1, \ldots, a_n).
\end{eqnarray*}
\begin{theo}
We obtain,
\begin{enumerate}
\item {$(\delta_{L,n} \otimes id)\delta_{R,n} = (id \otimes \delta_{R,n})\delta_{L,n} $ and $(\delta_{R,n}, \delta_{L,n})$
are homomorphisms.}
\item {If we define $\overrightarrow{d}_{G,n} = \Delta_n - \delta_{R,n}$ and
$\overleftarrow{d}_{G,n} = \tilde{\Delta}_n - \delta_{L,n}$,
then  $ (\overleftarrow{d}_{G,n} \otimes id)\overrightarrow{d}_{G,n}
= (id \otimes \overrightarrow{d}_{G,n})\overleftarrow{d}_{G,n} $.}
\item {$\overleftarrow{d}_{G,n} = \overrightarrow{d}_{G,n}$ if and only if they are applied to
element say, $[w]=(a_1, \ldots, a_n)$ such that $\Delta_n(a_1, \ldots, a_n)=(a_1, \ldots, a_n) \otimes a_1$
and $\tilde{\Delta}_n(a_1, \ldots, a_n)= a_n \otimes (a_1, \ldots, a_n)$. In this case,
$\overrightarrow{d}_{G,n}[w] = \overleftarrow{d}_{G,n}[w] = 0$.}
\item {$(\delta_{R,n},\delta_{L,n})$ carry a natural chirality. The notion of $G^{\otimes n}$-bimodule
is already valid and if the coproducts of the Markov $L$-coalgebra are unital homomorphisms
then $\overleftarrow{d}_{G,n}, \overrightarrow{d}_{G,n}$ are Ito derivatives.}
\end{enumerate}
\end{theo}
\Proof
Straightforward.
\eproof \\
\begin{coro}
If we note $\overrightarrow{d}_{G,n}((\cdot), w) := \Delta_{n,[w]} - \delta_{R,n}$, i.e.
$\overrightarrow{d}_{G,n}$ restricted to the path $[w]=(a_1, \ldots, a_n)$ of a graph $G$ and if
the path $[w]$ is a periodic orbit of period $n$, then
$\overrightarrow{d}_{G,n}(w, w) =0$. The same is valid for $\overleftarrow{d}_{G,n}$.
\end{coro}
\Proof
Straightforward.
\eproof \\
\subsection{Examples}
In the following examples, from a known algebra, we construct a graph so as to the algebra
embeds the Markov $L$-coalgebra, defined by the graph, into a Markov $L$-Hopf-algebra of degree 2.
\begin{exam}{[The triangle graph and quaternions]}
Here $k=\mathbb{R}$.
Recall that quaternions define the associative algebra $\mathbb{H}=( 1, i, j, k ) $ generated by the rules:
$$ ij=k, \ \ jk=i, \ \ ki=j, \ \ ii=jj=kk=-1. $$
The quaternions fit the present formalism by
considering the oriented triangle graph,\\
\begin{center}
\includegraphics*[width=5cm]{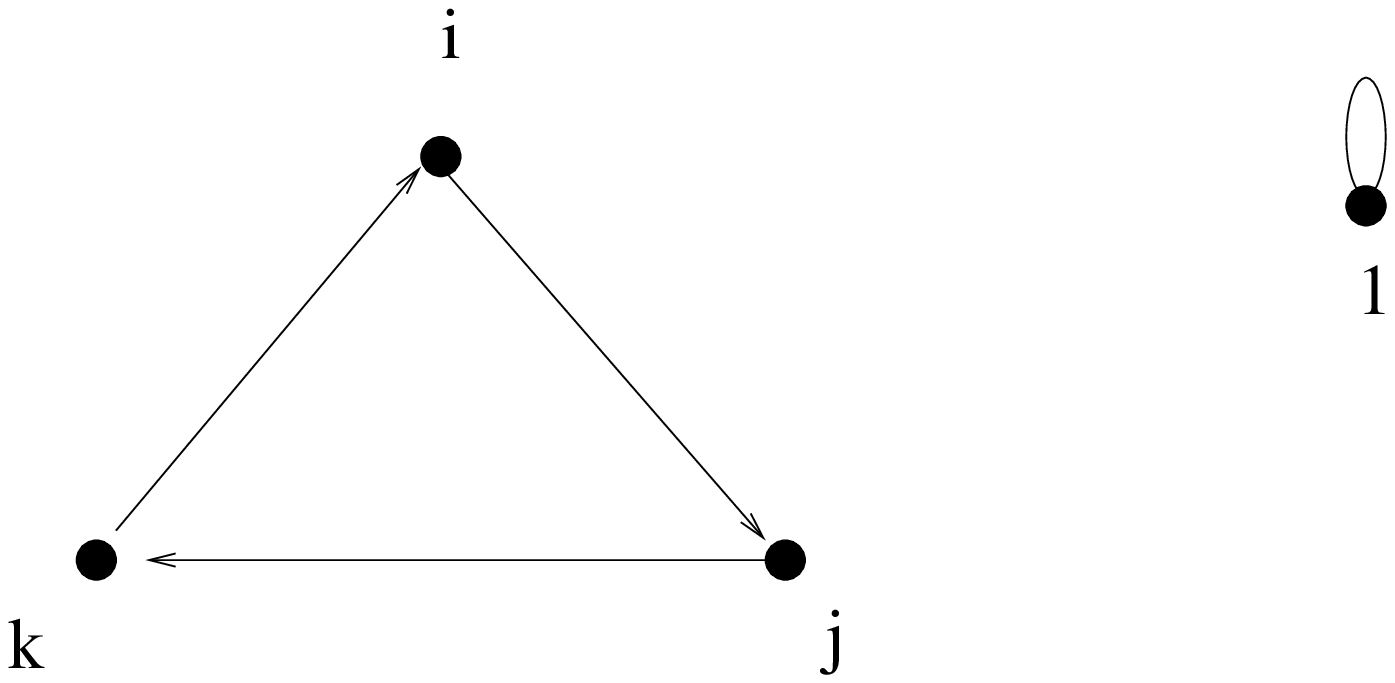}
\end{center}
Defining $x_0 \equiv i$, $x_1 \equiv j$, $x_2 \equiv k$ and adding subscripts
$\alpha, \beta \in \{0,1,2 \} \ \ \textrm{mod} \ 3$ i.e. $x_{\alpha+ \beta}
\equiv x_{\alpha+ \beta \textrm{mod} \ 3}$, we define,
$\Delta  x_{\alpha} = x_{\alpha} \otimes x_{\alpha +1}, \ \ \Delta 1 = \tilde{\Delta} 1 = 1 \otimes 1, \ \ \tilde{\Delta}x_{\alpha} = x_{\alpha -1} \otimes x_{\alpha}$
and
\begin{eqnarray*}
\Delta_2  (x_{\alpha} \otimes x_{\beta}) &=& x_{\alpha} \otimes x_{\beta} \otimes x_{\beta +1}, \ \ \ \epsilon_2(x_{\alpha} \otimes x_{\beta}) = x_{\alpha},\\
\tilde{\Delta}_2  (x_{\alpha} \otimes x_{\beta}) &=& x_{\alpha -1} \otimes x_{\alpha} \otimes x_{\beta}, \ \ \ \tilde{\epsilon}_2(x_{\alpha} \otimes x_{\beta})=x_{\beta},
\end{eqnarray*}
They embed the triangle graph into a Markov $L$-coalgebra.
\begin{theo}
The algebra of quaternions,
\begin{enumerate}
\item {embeds the triangle graph into a $L$-bialgebra of degree 2.}
\item{ Defining linear maps $S, \tilde{S}: \mathbb{H} \xrightarrow{} \mathbb{H}$
by $S(x_i)= -x_{i-1}$ and $\tilde{S}(x_{i-1})= -x_i$ for every $i \in \{ 0,1,2 \}$,
the $L$-bialgebra $\mathbb{H}$ becomes a $L$-Hopf algebra of degree 2, with $id \otimes S$ and $\tilde{S} \otimes id$ playing the r\^ole of right and left antipodes.}
\item{ The maps $S, \tilde{S}$ are unital antialgebra maps and satisfy $S\tilde{S}= id = \tilde{S}S $. They are the unique maps such that $id \otimes S$ and $\tilde{S} \otimes id$ be the
right and left antipodes
of $\mathbb{H}$ as a $L$-Hopf algebra of degree 2.}
\end{enumerate}
\end{theo}
\Proof
Let $i \in \{ 0,1,2 \}$. In the following we compute only the right coproduct part.
We show that $\Delta$ is a unital algebra map. Indeed, \\
$\Delta_2(x_{i} \otimes x_{j})\Delta_2(x_{i'} \otimes x_{j'})= (x_{i} \otimes x_{j} \otimes x_{j+1})
(x_{i'} \otimes x_{j'} \otimes x_{j'+1})= x_{i}x_{i'} \otimes x_{j}x_{j'} \otimes x_{j+1}x_{j'+1}, \ \ \textrm{and} \
\Delta_2(x_{i} \otimes x_{j})(x_{i'} \otimes x_{j'})= \Delta_2(x_{i}x_{i'} \otimes x_{j}x_{j'})= x_{i}x_{i'} \otimes x_{j}x_{j'} \otimes t(\Delta(x_{j}x_{j'})),$
where the Markov coproduct is defined by $\Delta  x_{\alpha} = x_{\alpha} \otimes x_{\alpha +1}$.
Hence we have to prove that $t(\Delta(x_{j}x_{j'})) = x_{j+1}x_{j'+1}$.
This is obvious by the following geometric proof. We suppose $j \not= j'$, $(x_{j},x_{j'})$ defines an edge of the triangle.
$(x_{j+1},x_{j'+1})$ defines the sole edge following it when we turn in a trigonometrical way.
Now we observe that up to a sign the concatenation of an edge, that is the product of its
source and its terminus give the third vertex. Hence by rotation the concatenation of
$(x_{j+1},x_{j'+1})$ will give the vertex just after. Thus up to a sign $t(\Delta(x_{j}x_{j'})) = x_{j+1}x_{j'+1}$.
The sign is easily obtained by noticing that if $(x_{j},x_{j'})$ is an arrow
of the triangle so is $(x_{j+1},x_{j'+1})$ and the sign is plus in both case when the concatenation
is realised. If the direction of $(x_{j},x_{j'})$ is in the opposite sens of an existing arrow, so
is $(x_{j+1},x_{j'+1})$ and the concatenation will give a minus sign in both cases.
In the case when $x_{j'}$ or $x_{j}$ is the identity element the proof is obvious since
there is a loop on it. The case $x_{j'} =x_{j}$ is trivial. \\
The coproducts $\Delta_2, \tilde{\Delta}_2$ are thus unital homomorphisms.
$\epsilon_2$ and $\tilde{\epsilon}_2$ are also a unital algebra map.
To prove the $L$-Hopf algebra of degree 2 part, we must prove that,
\begin{eqnarray*}
(id \otimes m) (id \otimes (id \otimes S)) \Delta_2 &=& \epsilon_2, \\
(m \otimes id) ((\tilde{S} \otimes id) \otimes id) \tilde{\Delta}_2 &=& \tilde{\epsilon}_2.
\end{eqnarray*}
This is obvious with the choice we made for the right and left antipodes.
$S$ is an antiunital map for
by definition, $-x_i=S(x_{i+1}) = S(x_{i-1}x_i)$ and $S(x_i)S(x_{i-1})=(-x_{i-1})(-x_{i-2})=(x_{i-1})(x_{i-2})=-(x_{i-2})(x_{i-1})=-(x_{i})$,
so $S(x_ix_j)=S(x_j)S(x_i)$. Moreover $S(x_ix_j) = S(x_j)S(x_i) = x_{j-1}x_{i-1}$ and
$S(x_ix_j) = -S(x_jx_i) = -S(x_i)S(x_j) = -x_{i-1}x_{j-1}= x_{j-1}x_{i-1}$ proving that $S$
is well defined.
$S$ is unital since $S(1) = S(x_i(-x_i)) = S(-x_i)S(x_i)=-(-x_{i+1}) (-x_{i+1}) = 1$.
$S$ is unique since if $S_1,S_2$ are two such right antipodes we must have  $ x_iS_1(x_{i+1}) = x_iS_2(x_{i+1})=1$
but $ x_ix_i=-1$ so $S_1(x_i) = S_2(x_i)$.
As $S_1,S_2$ are equal on the generators of the algebra they are equal. Moreover,
$S\tilde{S}(x_i) = S(-x_{i+1})=-(-x_i)=x_i.$ and $\tilde{S}S(x_i) = \tilde{S}(-x_{i-1})=-(-x_i)=x_i.$
\eproof
\Rk
$S,\tilde{S}$ are not unital anticoalgebra maps.
\Rk
As a coproduct, $\Delta$ is well defined on the oriented triangle graph, but is not an homomorphism
for the quaternion product. If it were the case, we would get, for example
$- \Delta(k) = - \Delta(ij) =\Delta(i) \Delta(j)= - ij \otimes jk = - k \otimes i $
which is true and
$\Delta(-k) = \Delta(ji) = \Delta(j) \Delta(i) = ji \otimes kj = (-k) \otimes (-i) $
which is still true. Yet we lost the $k$-linearity. Hence quaternions algebra is
an example of a Markov $L$-bialgebra of degree 2 which cannot be reduced to a
Markov $L$-bialgebra of degree 1.
\end{exam}
\begin{exam}{[The Pauli matrices]}
Here $k = \mathbb{C}$,
The Pauli matrices:
\[1_2= \begin{pmatrix}
 1 & 0\\
0 & 1
\end{pmatrix}, \ \ \
\gamma_0 = \begin{pmatrix}
 0 & 1\\
1 & 0
\end{pmatrix}, \ \ \
\gamma_1= \begin{pmatrix}
 0 & -i\\
i & 0
\end{pmatrix}, \ \ \
\gamma_2= \begin{pmatrix}
 1 & 0\\
0 & -1
\end{pmatrix},
\]
verify the algebra $\gamma_{k}\gamma_{k+1} = i\gamma_{k+2}$, $\gamma_{k}\gamma_{k}= 1_2$ and
$\gamma_{k}\gamma_{k+1} = - \gamma_{k+1}\gamma_{k}$. We know that $M_2(k)$
is the algebra generated by the Pauli matrices.
$M_2(k)$ fits the present formalism by
considering the oriented triangle graph with a loop on $1_2$ not represented here,
\[
\begin{array}{c@{\hskip 1cm}c@{\hskip 1cm}c}
\rnode{a}{ } & \rnode{b}{\gamma_{0}} & \rnode{c}{ }\\[1cm]
\rnode{d}{\gamma_{2}} & \ & \rnode{e}{\gamma_{1}}
\end{array}
\psset{nodesep=3pt}
\ncline{->}{b}{e}
\ncline{->}{e}{d}
\ncline{->}{d}{b} \ \ \ \ \ \ \ \ \ \ \ \
\begin{array}{c@{\hskip 1cm}c@{\hskip 1cm}c}
\rnode{a}{ } & \rnode{b}{i\gamma_{0}} & \rnode{c}{ }\\[1cm]
\rnode{d}{i\gamma_{2}} & \ & \rnode{e}{i\gamma_{1}}
\end{array}
\psset{nodesep=3pt}
\ncline{<-}{b}{e}
\ncline{<-}{e}{d}
\ncline{<-}{d}{b}
\]
The first one is to recall that $\gamma_{k}\gamma_{k+1} = i\gamma_{k+2}$, but it is
the second one which we are interested in because $(i\gamma_{k+1})(i\gamma_{k}) = (i\gamma_{k+2})$.
Defining $x_0 \equiv i\gamma_{0}$, $x_1 \equiv i\gamma_{1}$, $x_2 \equiv \gamma_{2}$ and adding subscripts
$\alpha, \beta \in \{0,1,2 \} \ \ \textrm{mod} 3$ i.e. $x_{\alpha+ \beta}
\equiv x_{\alpha+ \beta \textrm{mod}3}$, we define,
$\Delta  x_{\alpha} = x_{\alpha} \otimes x_{\alpha +1}, \ \ \Delta 1 = \tilde{\Delta} 1 = 1 \otimes 1, \ \ \tilde{\Delta}x_{\alpha} = x_{\alpha -1} \otimes x_{\alpha}$
and
\begin{eqnarray*}
\Delta_2  (x_{\alpha} \otimes x_{\beta}) &=& x_{\alpha} \otimes x_{\beta} \otimes x_{\beta +1}, \ \ \ \epsilon_2(x_{\alpha} \otimes x_{\beta}) = x_{\alpha},\\
\tilde{\Delta}_2  (x_{\alpha} \otimes x_{\beta}) &=& x_{\alpha -1} \otimes x_{\alpha} \otimes x_{\beta}, \ \ \ \tilde{\epsilon}_2(x_{\alpha} \otimes x_{\beta})=x_{\beta},
\end{eqnarray*}
They embed the triangle graph into a Markov $L$-coalgebra.
\begin{theo}
The algebra of Pauli matrices, i.e. $M_2(k)$,
\begin{enumerate}
\item {embeds the triangle graph into a $L$-bialgebra of degree 2.}
\item{ Defining linear maps $S, \tilde{S}: M_2(k) \xrightarrow{} M_2(k)$
by $S(x_i)= -x_{i-1}$ and $\tilde{S}(x_{i-1})= -x_i$ for every $i \in \{ 0,1,2 \}$,
the $L$-bialgebra $M_2(k)$ becomes a $L$-hopf algebra of degree 2, with $id \otimes S$ and $\tilde{S} \otimes id$ playing the r\^ole of right and left antipodes.}
\item{ The maps $S, \tilde{S}$ are unital antialgebra maps and satisfy $S\tilde{S}= id = \tilde{S}S $. They are the unique maps such that $id \otimes S$ and $\tilde{S} \otimes id$ be the
right and left antipodes
of $\mathbb{H}$ as a $L$-Hopf algebra of degree 2.}
\end{enumerate}
\end{theo}
\Proof
The proof is a corollary from the quaternion example. We only stress for instance that,
$S(x_k)= -x_{k-1}$ implies the equality $x_kS(x_{k+1})= -x_kx_{k}=-(i\gamma_k)(i\gamma_k)
=(\gamma_k)(\gamma_k)= 1_2$, usefull for computing the antipodes equalities.
\eproof
\end{exam}
\Rk
As unital associative algebras, the algebras studied in this subsection, equipped with their flower graphs, are also $L$-bialgebras of degree 1.
\section{Probabilistic algebraic products and mutation of $L$-coalgebras}
\subsection{Probabilistic algebraic products}
Let $G$ be a $k$-vector space. We equipped $G$ with a coproduct $\Delta_G: G \xrightarrow{} G \otimes G$. This coproduct
will give us a graph $G$, ($G$ can be for instance a $L$-coalgebra.).
The aim of this part is to view  a coproduct $\Delta_G$ as a product on a particular space. This idea comes
from the example of the graph of $Sl(2)_q$ and the triangle graph of the quaternions.
The most convenient mathematical tool for studying such a product is
the polynomial vector spaces.
We fix $n \in \mathbb{N}$ and denote $k \langle X_1, \ldots, X_n \rangle$ the $k$-vector space of
the polynomials constructed from the $X_i$. Let $G$ be a graph equipped with a
coproduct $\Delta_G$. We now view the $X_i$ as pointers that is as objects which will act on the
scalars from $k$. Let us see what this means.

Let $(a_1, \ldots, a_n) \in k$, $\sum_{i=1}^n a_iX_i$ will mean: $\sum_{i=1}^n a_i \lhd X_i$,
that is $X_i$ points or acts on $a_i$, thanks to the coproduct, $\Delta_G$, of the graph.
For this we equip $k \langle X_1, \ldots, X_n \rangle $ with a new (co)product $[\Delta_G]$, by defining the
following product:
$$(\sum_{i=1}^n a_i \lhd X_i)[\Delta_G](\sum_{i=1}^n b_i \lhd X_i) = \sum_{i=1}^n c_i \lhd \Delta(X_i),  $$
where the $a_i,b_i \in k$, and the $c_i \in k$ are for the moment a function of the $a_i,b_i$.
The following step is to express such a function. We simply carry the action of $\Delta(X_i)$
on the
scalars from $k$. This means, for example that if $n = 3$ and $G$ is the graph defined by $\Delta(X_1) = X_2 \otimes X_3$,
$\Delta(X_2) = X_2 \otimes X_2$ and $\Delta(X_3) = X_3 \otimes X_1$, we shall have by definition:
\begin{eqnarray*}
E &:=&(a_1 \lhd X_1 + a_2 \lhd X_2 + a_3 \lhd X_3)[\Delta](b_1 \lhd X_1 + b_2 \lhd X_2 + b_3 \lhd X_3) \\
&=& (c_1 \lhd \Delta(X_1) + c_2 \lhd \Delta(X_2) + c_3 \lhd \Delta(X_3) \\
&=& ((a_2 \otimes b_3) \lhd X_1 + (a_2 \otimes b_2) \lhd X_2 + (a_3 \otimes b_1) \lhd X_3).
\end{eqnarray*}
It is clear with this definition
that it is the coproduct of the graph $G$ which gives a product to $k \langle X_1, \ldots, X_n \rangle$.
We denote $(k \langle X_1, \ldots, X_n \rangle, [\Delta_G])$ the algebra, induced by the set $G$ equipped with a coproduct $[\Delta_G]$, (in the following, to simplify exposition
we consider only $L$-coalgebras).
\begin{theo}{\textbf{[The matrix product in $M_2(k)$]}}
Let $G$ be the graph of $Sl(2)_q$ and $\Delta$ its coproduct. Then,
$$(k \langle a, b, c ,d \rangle, [\Delta]) \simeq M_2(k)$$
\end{theo}
\Proof
We consider the graph of $Sl(2)_q$ with vertex $(a,b,c,d)$, (see part 2). We define
$X_{11} = a, X_{12} = b, X_{21} = c, X_{22} = d$ and compute:
\begin{eqnarray*}
E &:=& ( a_{11}X_{11} + a_{12}X_{12} + a_{21}X_{21} + a_{22}X_{22})(b_{11}X_{11} + b_{12}X_{12} + b_{21}X_{21} + b_{22}X_{22})\\
&=& (c_{11} \lhd \Delta(X_{11}) + c_{12} \lhd \Delta(X_{12}) + c_{21} \lhd \Delta(X_{21}) + c_{22} \lhd \Delta(X_{22})) \\
&=& ( (a_{11} \otimes b_{11} + a_{12} \otimes b_{21})X_{11} + (a_{11} \otimes b_{12} + a_{12} \otimes b_{22})X_{12} + \ldots
\end{eqnarray*}
Now we use the fact that $k \simeq k \otimes k$ to conclude.
\eproof \\
\begin{coro}
$(k \langle a, b, c ,d \rangle, [\Delta])$ is
a Markov $L$-bialgebra of degree 2.
\end{coro}
\Proof
This is a consequence from section 6.
\eproof
\Rk
It is interesting to note that the graph itself, equipped with its coproduct $\Delta$, (that is
the coassociative coalgebra) generates the matrix product in $M_2(k)$. We can
obviously extend this theorem to any dimension.
\Rk
Instead of using the field $k$ in $k \langle X_1, \ldots, X_n \rangle$ we can obviously use an associative algebra
$A$ with product $m$ and consider $A \langle X_1, \ldots, X_n \rangle$. If we compose the final result
of an operation induced by a $L$-coalgebra $G$ by the product of $m$ of the algebra $A$ we can find
for example the following corollary (we start this section by considering the field $k$ only for pedagogical reasons):
\begin{coro}
Let $G$ be the graph of $Sl(2)_q$ equipped with the coproduct $\Delta$. Then,
$$(A \langle a, b, c ,d \rangle, [\Delta]) \simeq M_2(A).$$
\end{coro}
\begin{theo}
Let $C$ be a coassociative coalgebra with $n < \infty$ elements, that is a graph $G$
with $n$ vertices
such that its coproduct $\Delta$ obeys the coassociative equation. (Recall that the coproduct
is defined from $T \bar{\otimes} T$ where $T$ is a matrix, see \ref{s1}.)\\
Then $(A \langle X_1, \ldots, X_n \rangle, [\Delta_C])$
is isomorphic to either $M_n(A)$ or to ${\bigoplus \atop{\sum n_i = n}} M_{n_i}(A)$ if $C$ is
a direct sum of subcoalgebras.
\end{theo}
\Proof
The fact that $(A \langle X_1, \ldots, X_n \rangle, [\Delta_C])$ is isomorphic to $M_n(A)$ follows immediately from the previous
corollary and the previous remarks. $C$ is the direct sum of subcoalgebras. Therefore each
subcoalgebra is isomorphic to a $M_{n_i}(A)$, subalgebra of $M_n(A)$, thus $(A \langle X_1, \ldots, X_n \rangle, [\Delta_C])$
is isomorphic to
${\bigoplus \atop{\sum n_i = n}} M_{n_i}(A)$. \\
(The matrix $T$ can always be split up in such a way
to produce ${\bigoplus \atop{\sum n_i = n}} M_{n_i}(A)$.)
\eproof \\
\begin{theo}
Let $P$ be a polynomial algebra with product $[\Delta_C]$, where $C$ is a graph equipped with a coproduct $\Delta_C$.
Asserting that $C$ is a coassociative coalgebra  is equivalent to the fact that the induced product $[\Delta_C]$ of $P$
is associative.
\end{theo}
\Proof
Obvious.
\eproof

How can we extend the previous concept to enlarge the number of different algebraic products?
We now focus on Markov $L$-coalgebra. The triangle graph of quaternions learns us that one has to
consider the position of the pointer $\Delta(X_i)$ too. For the moment, the pointer $\Delta(X_i)$
remained to the position occupied by the pointer $X_i$. It was a static point of view. The dynamical
viewpoint would be to move the pointer $\Delta(X_i)$ along, say, an orbit of a graph.
\begin{exam}{[The wedge product]}
Consider the oriented triangle graph $\equiv \triangle$, equipped with
the coproduct $\Delta_{\triangle}$ and $ (k \langle 1,X_0= i,X_1=j,X_2=k \rangle, [\Delta_\triangle])$.
We choose to put $\Delta(X_i)$ in the position $X_{i+2 \ \textrm{mod} \ 3}$
and to fix $\Delta(1)$ in the position
occupied by 1.
We have the following, where $x,y,z \in k$
\begin{eqnarray*}
E &:=& (xX_{0} + yX_{1} + zX_{2}) [\Delta_\triangle] (x'X_{0} + y'X_{1} + z'X_{2})\\
&=& ( xy'X_{2} + yz'X_{0} + zx'X_{1}).
\end{eqnarray*}
\Rk
$\Delta(X_0)$ is no longer in the position occupied by $X_0$ but in the position occupied by $X_2$.
\begin{theo}
The commutator induced by the triangle product on
$(\mathbb{R} \langle 1,X_0= i,X_1=j,X_2=k \rangle, [\Delta_\triangle])$ is isomorphic to $(\mathbb{R}^{3}, \ \wedge)$, where
$\wedge$ denotes the standard wedge product.
\end{theo}
\Proof
Let $\vec{a}=(x,y,z), \ \vec{b}=(x',y',z') \in \mathbb{R}^{3}$.
We compute:
\begin{eqnarray*}
E &:=& (xX_{0} + yX_{1} + zX_{2}) [\Delta_\triangle] (x'X_{0} + y'X_{1} + z'X_{2}) \\
& & - (x'X_{0} + y'X_{1} + z'X_{2}) [\Delta_\triangle] (xX_{0} + yX_{1} + zX_{2})\\
&=& (xy'X_{2} + yz'X_{0} + zx'X_{1}) - (x'yX_{2} + y'zX_{0} + z'xX_{1})\\
&=& (xy'- x'y)X_{2} + (yz'-y'z)X_{0} + (zx'-z'x)X_{1}\\
&=& \vec{a} \wedge \vec{b},
\end{eqnarray*}
since the field $\mathbb{R}$ is commutative. If it was not, this could be a possible generalisation
of the wedge product in the non commutative case.
\eproof
\end{exam}
\begin{exam}{}
For another important example in mathematics see \cite{Lergraph}.
\end{exam}
For the moment all the graphs and products involved were deterministic. Suppose now that we consider
a more complicated graph with probability measure on it and we choose to place $\Delta(X_i)$ along
some well-chosen orbit. We shall obtain a non deterministic algebraic product on $A \langle X_1, \ldots, X_n \rangle$,
with $A$ an associative algebra with product $m$. Such a mathematical object $(A \langle X_1, \ldots, X_n \rangle, [\Delta_G])$
will be called \textit{a random polynomial algebra}.
For the moment the author does not know any example of physical or mathematical applications of such
a concept.
As an example we consider the following Markov $L$-algebra $G$:
\begin{center}
\includegraphics*[width=5cm]{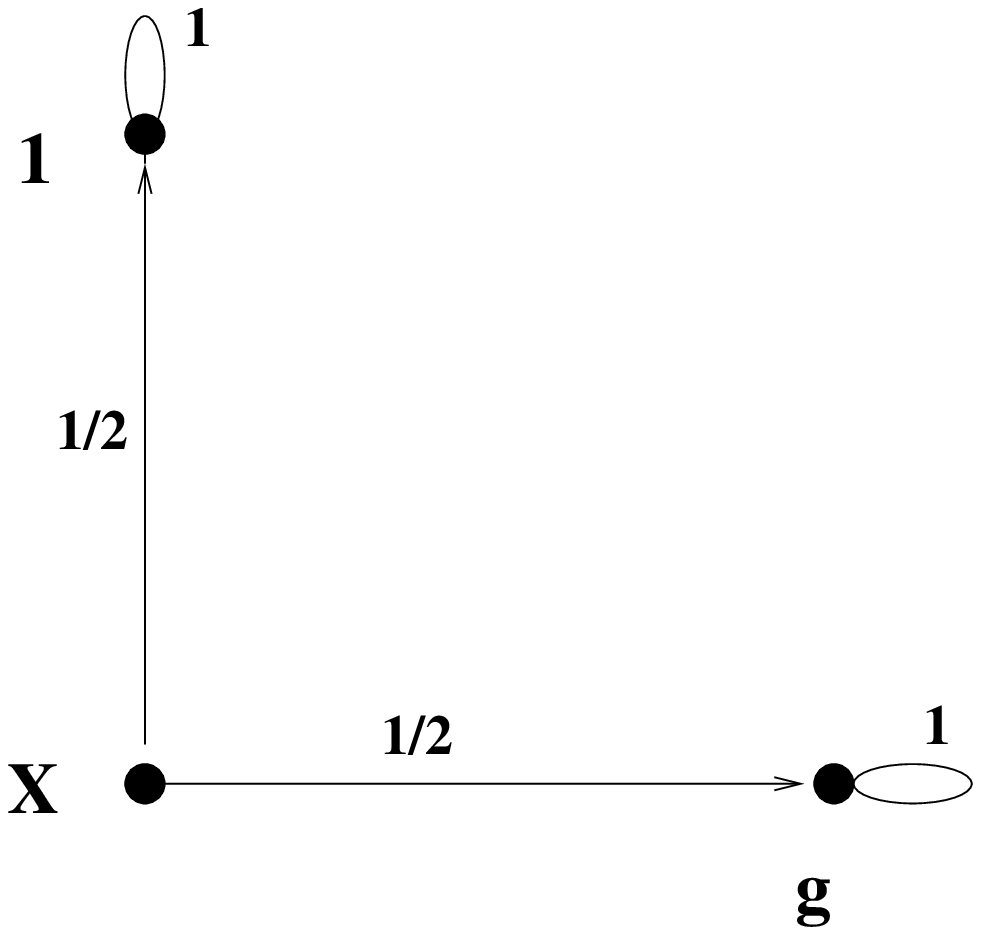}
\end{center}
and we choose the convention, called
the path convention in the sequel:
$$\lhd \Delta(\cdot) \equiv terminus(\Delta(\cdot)) = t(\Delta(\cdot)).$$
That is we choose to put the result of the operation say, $\Delta(X_i) := X_i \otimes X_j$ into
the place occupied by the pointer $\lhd X_j$.
For the moment we do not consider probability on $G$.
\begin{theo}
$(A \langle X_{11} =1, X_{12}=X , X_{22}={g} \rangle, [\Delta_G])^2 \simeq diag(M_2(A))$
where power 2 means that we consider the set $\{ a[\Delta_G]b, \ (a,b) \in (A \langle 1,X,g \rangle, [\Delta_G)] \}$.
\end{theo}
\Proof
We must only focus on $X$.
Yet, $\Delta(X) = X \otimes g + X \otimes 1$ that is the result of a product will be put on $g$ and $1$
that is in $X_{11}$ and $X_{22}$. For instance the product of two polynomials yields $(aX + bg + c1)[\Delta_G](a'X + b'g + c'1) := (ab'+bb') g + (ac'+cc')1$.
Moreover $diag(M_2(A)) \simeq (A \langle X_{11}, X_{22} \rangle, [\Delta_H])$ where $H$ is the graph with two loops
indexed by the pointers 1 and $g$.
\eproof

Now suppose we affect a probability $\frac{1}{2}$ to each arrow emerging from $X$ and a probability
1 to the two loops.
Then we have: $\Delta(X) = \frac{1}{2} X \otimes g + \frac{1}{2} X \otimes 1$. Equipped with this
graph the random walks on it, turn $(A \langle 1,X,g \rangle, [\Delta_G]) \}$ into a random algebra. To study a product
in such an algebra we must consider a probability measure $\mathbb{P}$
on $(A \langle 1,X,g \rangle , [\Delta_G])$. In fact the probability measure on paths of the random walk on $G$ is sufficient.
For example the probability to have the walk $w_1 = (X,1,1, \ldots)$ is equal to $\frac{1}{2}$ and the
probability to have $w_2 = (X,g,g, \ldots)$ is equal to $\frac{1}{2}$.
We now assign to $\mathbb{P}(a[\Delta_G]b = a[\textrm{walk:} \ w_1]b) = \frac{1}{2}$ and
$\mathbb{P}(a[\Delta_G]b = a[\textrm{walk:} \ w_2]b) = \frac{1}{2}$.
\begin{theo}
$\mathbb{P}\{(A \langle 1,X,g \rangle, [\Delta_G])^2 \simeq diag(M_2(A)) \} = 1$.
\end{theo}
\Proof
$(A \langle 1,X,g \rangle, [\Delta_G])$ already contains the sub random polynomial algebra $(A \langle 1,g \rangle, [\Delta_G])$
which is isomorphic to $diag(M_2(A))$, all we do by considering the power 2 of $(A \langle 1,X,g \rangle, [\Delta_G])$
is to eliminate $X$ to restrict ourselves to the attractors generated by the loops $1$ and $g$.
\eproof
\Rk
This theorem means that the random polynomial algebra $(A \langle 1,X,g \rangle, [\Delta_G])$ can behave as a deterministic
algebra at short term.
\Rk
The fact that we recover a deterministic algebra is due to the fact that a loop
---the $L$-coalgebra generated by a loop---is a coassociative coalgebra, in addition to the fact
that it is placed in an attractor r\^ole. This remark allows us to generalise the previous theorem by
saying that if the $L$-coalgebra generated by $G$ has an attractor $C$, where $C$ is a coassociative
coalgebra, and if we decide to choose the path convention, except on $C$ where
we choose the static one, then it exists a time $n$, possibly equal to infinity such that
the power of a random polynomial algebra equipped with the product $[\Delta_G]$ converges towards a deterministic
algebra.
\subsection{Mutation of $L$-coalgebras}
The aim of this part is twofold.
The first idea is to consider the graphs, embbeded into $L$-coalgebras,
as dynamical objects, capable of mutation
in order to consider the notion of mutation of algebraic products on polynomial algebras. The
second idea is to produce an example of sequence based on random variables.
\begin{defi}{[Mutation of $L$-coalgebras]}
If $M$ and $N$ are oriented graphs embedded into $L$-coalgebras
we denote:
$$ (M,\Delta_M)  \looparrowright (N,\Delta_N) $$
to say that the $L$-coalgebra $(M,\Delta_M)$ has undergone a mutation into the $L$-coalgebra $(N,\Delta_N)$.
\end{defi}
\begin{defi}{[Mutation of algebraic product]}
As we can associate with an $L$-coalgebra, an algebraic product we define: $ [\Delta_M] \looparrowright [\Delta_N] $
to say that the algebraic product $[\Delta_M]$ has undergone an algebraic mutation into $[\Delta_N]$.
\end{defi}
So as to be as clear as possible
we shall illustrate all these new concepts through an example.

Let $H_2$ be the graph associated with the coassociative coalgebra $SL(2)_q$, whose the vertex set is still denoted by $a,b,c,d$; $\diamondsuit_{a}$
the coassociative coalgebra represented by a loop at $a$ and $\triangle$ the $L$-Markov algebra associated with
the oriented triangle graph. For instance we choose to label the pointers of the triangle graph by $a \xrightarrow{} b \xrightarrow{} c \xrightarrow{} a$.
As we saw in the previous part, these $L$-coalgebras
define respectively the matrix product on $M_2(A)$, the standart product on $A$ and the wedge product on $A$, where $A$ is an associative algebra.
Let us consider the set $G = \{ \{H_2 \}, \  \{\diamondsuit_{a} \}, \  \{ \triangle  \} \} $ and
define a probability measure on the family of subsets of $G $,
$\mathbb{P} : \mathcal{P}(G) \xrightarrow{} [0,1]$ such that:
$$\mathbb{P}(\{\triangle\}) = \frac{\epsilon}{2}, \ \ \  \mathbb{P}(\{\diamondsuit_{a} \}) = \frac{\epsilon}{2}, \ \ \ \mathbb{P}(\{H_2 \}) = 1-\epsilon, \ \ \epsilon \in \ ]0,1[ .$$
Naturally associated with $G$ is the dynamical polynomial algebra, $(\mathcal{G}, [\Delta_G])$, whose pointers are labelled by $a,b,c,d$ as well and which
behaves for instance as an algebra isomorphic to $M_2(A)$ when $G$ behaves as $H_2$.
Suppose now we obtain the following dynamical sequence:
$$ H_2 \xrightarrow{} H_2 \xrightarrow{} H_2 \looparrowright \diamondsuit_{a} \looparrowright H_2 \xrightarrow{} H_2 \ldots$$
This dynamic will have repercussions on the algebraic products:
$$ [\Delta_{H_2}] \xrightarrow{} [\Delta_{H_2}] \xrightarrow{} [\Delta_{H_2}] \looparrowright [\Delta_{\diamondsuit_{a}}] \looparrowright [\Delta_{H_2}] \xrightarrow{} [\Delta_{H_2}] \ldots$$
this means that if we study the power of a 2 by 2 matrix $z$ we will get:
$$ z \xrightarrow{\textrm{matrix product}} z^2 \xrightarrow{\textrm{matrix product}} z^3
\xrightarrow{\textrm{matrix product}} z^4 \looparrowright u \xrightarrow{\textrm{matrix product}} u^2 \xrightarrow{\textrm{matrix product}} u^3 \ldots$$
where $z$ has undergone a mutation into the 2 by 2 matrix $u$.
\begin{exam}{}
This example of mutation looks like a standard projection on the pointer $a$.
\[\boldsymbol{z} = \begin{pmatrix}
 1 & 1\\
1 & 1
\end{pmatrix}
\xrightarrow{}
\boldsymbol{z^2} = \begin{pmatrix}
 2 & 2\\
2 & 2
\end{pmatrix}
\xrightarrow{}
\boldsymbol{z^3} =\begin{pmatrix}
 8 & 8\\
8 & 8
\end{pmatrix}
\xrightarrow{}
\boldsymbol{z^4} =\begin{pmatrix}
 128 & 128\\
128 & 128
\end{pmatrix}
\looparrowright
\boldsymbol{u} =\begin{pmatrix}
 (128)^2 & 0\\
0 & 0
\end{pmatrix}
\]
\\
\[
\looparrowright
\boldsymbol{u^2} =\begin{pmatrix}
 (128)^4 & 0\\
0 & 0
\end{pmatrix}
\xrightarrow{}
\boldsymbol{u^3} =\begin{pmatrix}
 (128)^6 & 0\\
0 & 0
\end{pmatrix}
\ldots
\]
\end{exam}
To explain what follows we need
to view the set $G$ as a dynamical $L$-coalgebra. We start for example in a configuration where $G$ is viewed as the
coassociative coalgebra
$H_2$. We ask
how such a graph evolves. There is three possibilities. It remains the same, contracts itself into the loop labelled by $a$ or modifies
its shape to become the oriented triangle $L$-coalgebra described above (which is a subgraph of the graph associated with $H_2$).
Fix $\epsilon >0$ and consider the set of dynamical sequences starting with $H_2$, i.e. $(H_2 \xrightarrow{}\ldots)$,
and denote by $T = \inf  \{n>0, \textrm{the} \  \triangle  \ \textrm{occurs at time n} \}$. Then, the probability for the following sequence
$$ 0 \xrightarrow{} G \xrightarrow{\Delta_G} G^{\otimes 2} \xrightarrow{\Delta_G \otimes id - id \otimes \Delta_G} G^{\otimes 3}
\ldots $$
not to be an exact complex at time $T$ is $(1-\frac{\epsilon}{2})^{(T-1)}\frac{\epsilon}{2}$. By reversing the arrows, the probability for the following sequence
$$ 0 \xleftarrow{} \mathcal{G} \xleftarrow{[\Delta_G]} \mathcal{G}^{\otimes 2} \xleftarrow{[\Delta_G] \otimes id - id \otimes [\Delta_G]} \mathcal{G}^{\otimes 3}
 \ldots $$
not to be an exact complex at time $T$ is $(1-\frac{\epsilon}{2})^{(T-1)}\frac{\epsilon}{2}$. We recall here that $\mathcal{G}$ is the polynomial
algebra over $A$, equipped with the product $[\Delta_G]$, which is isomorphic to $M_2(A)$ if $[\Delta_G]=[\Delta_{H_2}]$ and so on.
Conversely if we fix $\epsilon >0$ and consider now the set of dynamical sequences starting with the $\triangle$ $L$-coalgebra, i.e. $(\triangle \xrightarrow{}\ldots)$,
and denote by $T_1 = \inf  \{n>0, \textrm{a mutation occurs at time n} \}$. Then, the probability for the following sequence
$$ 0 \xrightarrow{} G  \xrightarrow{\Delta_G} G^{\otimes 2} \xrightarrow{\Delta_G \otimes id - id \otimes \Delta_G} G^{\otimes 3}
\ldots $$
to be an exact complex at time $T_1$ is $(1-\frac{\epsilon}{2})\frac{\epsilon}{2}^{(T_1-1)}$. By reversing the arrows, the probability for the following sequence
$$ 0 \xleftarrow{} \mathcal{G} \xleftarrow{[\Delta_G]} \mathcal{G}^{\otimes 2} \xleftarrow{[\Delta_G] \otimes id - id \otimes [\Delta_G]} \mathcal{G}^{\otimes 3}
 \ldots $$
to be an exact complex at time $T_1$ is $(1-\frac{\epsilon}{2})\frac{\epsilon}{2}^{(T_1-1)}$.
\Rk
When $G$ is represented by the loop or the graph of $Sl(2)_q$ the previous sequence are exact complex.
Yet if a mutation occurs, i.e. if $G$ becomes $\triangle$,
it will break the exactness of such a complex because the usually matrix product will undergo too
an algebraic mutation \footnote{As an example of application, imagine that some quantum measurements
are done on a quantum system which lives on a space-time represented  by the dynamical $L$-coalgebra  $G$, viewed as $H_2$.
Instead of disturbing the quantum system by the measurement, let us
suppose that we disturb $G$ and that $G$ undergoes a contraction in the loop labelled by $a$. This will
produce a change of complex for the coassociative coalgebras which will induce a mutation of algebraic product, here a projection
on the pointer labelled by $a$.

Another application of probabilistic algebraic product would be to embed the fundamental biological bricks, i.e. $A,C,G,T$ into
a random semigroup. }.
\newpage
\section{Conclusion and open problems}
Through these six parts, the author hopes to have demonstrated the interest of the
concept of $L$-coalgebra. The main idea was to unify concepts from probability and
combinatoric theory, especially
from oriented graphs
and concepts from coassociative
coalgebra theory in a new tool, called $L$-coalgebra and to show that it is necessary to consider coproducts, which allow us to manipulate
both geometry and algebra rather than the graph in its geometrical definition. Let us summarize the main ideas of the differents sections.
In affecting with each tensor product an arrow we started
to associate a graph with a coassociative coalgebra. We noticed then the non locality of the
product, (propagator) $\Delta$. To recover the locality we had to break the coassociativity
of the coproduct and to create a new coproduct $\tilde{\Delta}$, both the coproducts obeying
the coassociativity breaking equation. We gave also numerous examples
of $L$-coalgebras by showing, notably, how we could embed the quaternions and $M_2(k)$
into an oriented triangle $L$-Hopf algebra of degree 2 and a coassociative coalgebra into an Ito $L$-coalgebra.
Yet the most relevant example remains the flower graph
associated with a unital algebra.
With the notion of curvature, due to Quillen, we constructed an equation, based on the coproducts of the flower graph, which
shows a common point with the product of an associative unital algebra $A$, homomorphisms, Ito derivatives and Leibnitz derivatives
from $A$ to $A$. This equation was the started point to the realization of a di-superalgebra from the curvature of an Ito map whose the
integral calculus yields cyclic cocycles and vanishes on the Leibnitz bracket. We showed also
a bijection between homomorphisms and Ito derivatives and yielded a connection between Ito maps and the third Reidemeister
movement.
In section 6, by replacing the flower graph by another
Markov $L$-coalgebra, we showed the usefulness of virtual petals in the construction of Ito derivatives.

Thanks to the coproducts, we enlarged the notion of graph by removing the denumerability condition and enlarged
the concept of walk on a graph to walk on a $L$-coalgebra. This idea was used in section 4, to generalise the $\circ$ operation in the
context of completely positive semigroups to any coproduct. In section 7,
We also showed that we could recover the matrix product only by considering what was
called a polynomial algebra equipped with the product induced by a $L$-coalgebra, that is by
the propagator $\Delta$ of a coassociative coalgebra.
By studying the example of $Sl(2)_q$ and the example of the wedge product
we discovered the notions of probabilistic
algebraic product, of random polynomial algebra and mutation through the concept of dynamical $L$-coalgebra.

There are a lot of open questions. One of them is how to use for instance
the notion of recurrence,
ergodicity, transience of random walks on oriented graphs \cite{Petritis} \cite{Petritis2}
into the algebraic $L$-coalgebra framework.

\textbf{Acknowledgments:}
The author wishes to thank Dimitri Petritis for useful discussions and fruitful advice for the
redaction of this paper and R.L. Hudson as well for communicating him
its results prior to their publication.

\bibliographystyle{plain}
\bibliography{These}

\end{document}